**L'identité algébrique d'une pratique portée par *la discussion sur l'équation à l'aide de laquelle on détermine les inégalités séculaires des planètes* (1766-1874).**

Frédéric Brechenmacher ([1]).

---

**Résumé.** Cet article questionne l'identité algébrique d'une pratique propre à un corpus de textes publiés sur une période antérieure à l'élaboration de théories algébriques, comme la théorie des matrices, ou de disciplines, comme l'algèbre linéaire, qui donneront à cette pratique l'identité d'une méthode de transformation d'un système linéaire par la décomposition de la forme polynomiale de l'équation caractéristique associée. Dans les années 1760-1775, Lagrange élabore une pratique algébrique spécifique à la mathématisation des problèmes mécaniques des petites oscillations de cordes chargées d'un nombre quelconque de masses ou de planètes sur leurs orbites. La spécificité de cette pratique s'affirme par opposition à la méthode des coefficients indéterminés, elle consiste à exprimer les solutions des systèmes linéaires par des factorisations polynomiales d'une équation algébrique particulière*, l'équation à l'aide de laquelle on détermine les inégalités séculaires des planètes*. Elaborée en un jeu sur les primes et les indices des coefficients des systèmes linéaires, la pratique de Lagrange est à l'origine d'une caractéristique des systèmes issus de la mécanique, la disposition en miroirs des coefficients de systèmes que nous désignons aujourd'hui comme symétriques. Elle est également à l'origine d'une discussion sur la nature des racines de l'équation qui lui est associée. Nous questionnerons l'identité algébrique de cette discussion qui se développera sur plus d'un siècle en étudiant les héritages, permanences et évolutions de la pratique élaborée par Lagrange au sein de différentes méthodes élaborées dans divers cadres théoriques par des auteurs comme Laplace, Cauchy, Weierstrass, Jordan et Kronecker. Nous verrons que, préalablement à l'élaboration d'une théorie des formes dont la nature, algébrique ou arithmétique, suscitera une vive controverse entre Jordan et Kronecker en 1874, le caractère algébrique de la *discussion* renvoie davantage à l'identité historique d'un corpus qu'à une identité théorique et s'avère indissociable d'une constante revendication de généralité. Entre 1766 et 1874, la discussion sur l'équation à l'aide de laquelle on détermine les inégalités séculaires permet de mettre en évidence différentes représentations associées à une même pratique algébrique ainsi que des évolutions dans les philosophies internes portées par les auteurs du corpus sur la généralité de l'algèbre.

---

---

[1] F. BRECHENMACHER. Laboratoire de Mathématiques Lens (LML, EA2462). Fédération de Recherche Mathématique du Nord-Pas-de-Calais (CNRS, FR 2956). Université d'Artois (IUFM du Nord Pas de Calais).
Faculté des Sciences Jean Perrin, rue Jean Souvraz S.P. 18, 62 307 Lens Cedex France.
Courrier électronique : frederic.brechenmacher@euler.univ-artois.fr.

## Introduction.

Cet article est consacré à un corpus regroupant des textes mathématiques en un ensemble cohérent et fermé sur une périodisation donnée (1766-1874) et que nous désignerons sous le nom de *discussion sur l'équation à l'aide de laquelle on détermine les inégalités séculaires des planètes* ou, plus rapidement, *discussion sur l'équation des petites oscillations.* Au sens des citations et des références que les textes entretiennent entre eux, visualisées par le graphe de l'annexe 1, le corpus articule une discussion entre différents auteurs sur plus d'un siècle. Quelques nœuds apparaissent dans l'enchevêtrement des références bibliographiques. En 1766, Joseph-Louis Lagrange travaille sur l'intégration des systèmes d'équations différentielles linéaires à coefficients constants ; en 1829, la question traitée par Augustin Cauchy est la classification des surfaces du second degré ; en 1858, il s'agit pour Karl Weierstrass de caractériser les transformations de couples de fonctions homogènes du second degré par substitutions sur les variables ; en 1874 Camille Jordan et Leopold Kronecker se querellent sur les fondements et les méthodes de la théorie des formes bilinéaires. Le corpus ne se laisse donc pas aisément caractériser et, comme nous le verrons, ne s'identifie ni à une théorie, ni au développement d'une méthode ni même à la résolution d'un problème.
En envisageant ce corpus comme une *discussion*, nous voulons questionner son *identité* par l'étude des liens qu'entretiennent les textes qui le constituent. Cette question d'identité nous amène à nous démarquer de l'historiographie qui s'est intéressée au corpus dans l'objectif de dresser l'histoire d'une théorie que l'on peut désigner, pour reprendre le terme employé par Thomas Hawkins, comme la théorie spectrale des matrices ([2]). La théorie des matrices permet de formuler les différents problèmes étudiés par les auteurs du corpus comme un même problème de réduction d'un couple de matrices *(A,B)* en un couple (*D, I*) où *A* est une matrice symétrique, *B* une matrice symétrique définie, *D* une matrice diagonale et *I* la matrice identité. Une telle *identité théorique* ne s'élabore cependant qu'à la limite de la périodisation associée à la *discussion* (1766-1874), elle implique donc un regard rétrospectif sur le corpus et masque d'autres identités qui ne relèvent pas de ce cadre théorique. Par exemple, le recours à une formulation matricielle de la méthode élaborée par Lagrange en 1766 pour la résolution des systèmes différentiels linéaires à coefficients constants introduit dans le discours historique des aspects propres à la théorie des matrices comme la *transformation* d'une matrice *symétrique* qui sont absents du texte de Lagrange et masquent le caractère spécifique d'une *pratique* originale. Comme nous le détaillerons dans cet article, lorsque Lagrange ramène l'intégration d'un système de *n* équations à celle de *n* équations indépendantes, le changement de *forme* du système n'est pas envisagé comme le résultat une transformation mais s'avère indissociable d'une représentation mécanique implicite selon laquelle les petites oscillations d'un système de *n* corps se comportent comme si elles étaient composées de *n* mouvements simples ([3]). Nous verrons également que ce que nous désignerions aujourd'hui comme le caractère symétrique des systèmes linéaires de la mécanique se manifeste d'abord comme une conséquence de la *pratique algébrique* élaborée par Lagrange par opposition aux insuffisances de la méthode des coefficients indéterminés et consistant, en un jeu sur les primes et les indices des coefficients des systèmes linéaires que l'on interpréterait aujourd'hui comme une mise en œuvre de propriétés d'orthogonalité duale, à réduire l'intégration de ces derniers à la décomposition algébrique d'une équation ([4]). En envisageant l'identité de notre corpus comme « un problème et non une tautologie » pour reprendre l'expression employée par Catherine Goldstein à propos des différentes identités d'un théorème de Fermat [Goldstein 1995, p.16], notre problématique vise à questionner le caractère algébrique qui fait la spécificité d'une pratique propre au corpus sur une période antérieure à l'élaboration des théories algébriques qui donneront à cette pratique l'identité d'une méthode. Afin de restituer la « dynamique réelle des savoirs » par la « multiplicité de leurs origines » [Dhombres 2002], nous nous attacherons à mettre en évidence aussi bien la variabilité des contextes dans lesquels notre corpus se déploie que les héritages manifestant la *permanence d'une pratique* spécifique participant des *différentes méthodes* élaborées par des auteurs comme Lagrange, Cauchy, Weierstrass, Jordan et Kronecker. Il faudra donc entrer dans la technicité des textes originaux afin de saisir, comme le formulait Hourya Sinaceur dans son histoire du théorème de Sturm, « toutes les identités que le "progrès" efface : identité d'un contexte, d'un objectif, d'une perspective, d'un langage, sans parler de tout ce qui reste implicite sans manquer d'être là » [Sinaceur 1991, p.21].

---

[2] Nous proposons dans cet article une approche complémentaire aux travaux menés par l'historien T. Hawkins de l'université de Boston et auxquels nos recherches sont très largement redevables. Entre 1974 et 1977, Hawkins a consacré plusieurs articles à l'histoire de la théorie spectrale des matrices, il en a fixé la structure qui a ensuite été reprise par les traités d'histoire des mathématiques proposant une histoire de l'algèbre linéaire. Une première étape voit une « origine » de la théorie spectrale dans les travaux de mécaniques du XVIII^e^ siècle impliquant des systèmes différentiels linéaires à coefficients constants du type $AY''=BY$ où $A$ est une matrice symétrique et $B$ une matrice symétrique et définie. L'intégration de ces systèmes repose sur la recherche de valeurs propres $\lambda_i$ telles que $AX=\lambda_i BX$ permettant d'écrire le système différentiel sous la forme $DY''=Y$ où $D=(\lambda_1, \lambda_2, \dots, \lambda_n)$ est une matrice diagonale. A cette première étape succède un développement théorique initié par Cauchy (1829) et complété par Weierstrass (1858) et qui se caractérise par la démonstration de la nature réelle des valeurs propres des matrices symétriques. La troisième étape correspond à l'organisation d'une théorie, la théorie des formes bilinéaires fixée par Frobenius en 1878, autour d'un théorème énoncé par Weierstrass en 1868 permettant de caractériser les classes de similitude des matrices par un système complet d'invariants, les diviseurs élémentaires.

[3] Dans le cadre de son étude sur l'évolution des fondements de la mécanique des systèmes discrets chez Lagrange, M. Panza a mis en évidence le rôle essentiel joué par l'« interprétation mathématique des concepts mécaniques essentiels » [Panza 1992, p.205]. Nous verrons dans les troisième et quatrième parties de cet article que la pratique algébrique élaborée par Lagrange pour le problème des petites oscillations est indissociable d'une représentation mécanique consistant à envisager les oscillations d'un système de *n* corps comme la composition de *n* périodes propres.

En ne prenant pas pour appui un cadre théorique comme celui de la théorie des matrices, comment constituer le corpus auquel consacrer notre étude ? Notre méthodologie est basée sur le choix d'un moment de référence caractérisé par une date, 1874, donnant une limite à la période étudiée et un contexte, la querelle qui oppose Jordan à Kronecker ([5]). Les références de Jordan et Kronecker aux travaux de leurs prédécesseurs nous ont permis de fixer un premier groupe de textes que nous avons ensuite complété par un épuisement systématique des références bibliographiques jusqu'à obtenir le corpus fermé dont nous proposons une représentation simplifiée en annexe 1. D'autres choix de dates et de contextes de référence pourraient faire apparaître des réseaux différents ou faire jouer des rôles importants à des textes apparaissant ici comme des références mineures. Un autre point de vue, très proche puisque provenant des années 1880 et d'auteurs comme Henri Poincaré, conduirait à ouvrir le corpus à des travaux d'astronomie et à centrer le questionnement sur la stabilité du système du monde et le problème des trois corps. En prenant pour moment de référence la date de 1890 et les travaux d'Eduard Weyr sur les matrices, nous serions amenés à donner une plus grande importance aux travaux dans lesquels James Joseph Sylvester définit, en 1851, le terme « matrice » mais qui ne jouent ici qu'un rôle mineur et ne sont cités que par Charles Hermite entre 1853 et 1857 puis par Gaston Darboux en 1874 ([6]).

Les première et troisième parties de cet article sont consacrées aux deux moments particuliers qui marquent respectivement la fin et l'origine de notre corpus tandis que les deuxième et quatrième parties abordent le corpus dans son ensemble. Le choix d'adopter pour moment de référence la date de 1874 et le contexte de la querelle qui oppose Jordan et Kronecker est un point essentiel de la méthodologie de notre étude qui, bien que portant sur une période longue, n'en est pas moins locale car orientée par une problématique visant à questionner l'identité algébrique d'un corpus bien précis ; la première partie de cet article est consacrée à expliciter et justifier ce choix. Dans la deuxième partie nous questionnons l'identité du corpus par une étude des liens qu'entretiennent les textes qui le constituent. Nous verrons que cette identité est indissociable de celle d'une pratique associée à une équation algébrique spécifique dont nous étudierons d'abord, dans la troisième partie, l'origine dans les travaux de Lagrange de 1766 puis, dans la quatrième partie, les héritages et évolutions au sein de la *discussion*. Nous verrons que l'identité algébrique de cette pratique ne se réduit ni à un langage ni à un caractère formel mais que différentes représentations lui sont associées dans différents cadres théoriques. La pratique élaborée par Lagrange est notamment indissociable de représentations mécaniques et son caractère algébrique ne se réduit pas au formalisme souvent associé par l'historiographie au style de mathématisation et au modèle de théorie physique de l'auteur de la *Mécanique Analytique* où la « mise en œuvre du calcul des variations nécessite une forme d'intuition nouvelle, en ce sens qu'elle est entièrement détachée de la considération des propriétés des figures géométriques, car elle opère à l'intérieur d'un cadre de formules algébriques exclusivement. […] L'introduction du formalisme mathématique dans la dynamique signifie qu'on impose une forme esthétique algébrique à l'étude de la science du mouvement, cela veut dire que l'on abandonne tout moyen autre que celui du calcul algébrique. » [Barroso-Filho et Comte 1988, p.134]([7]).

---

[4] Comme nous le verrons dans le paragraphe I.3. et la note n°23, la pratique de Lagrange sera plus tard exprimée sous la forme d'une méthode revenant à donner une expression polynomiale générale des vecteurs propres d'une matrice comme quotients de mineurs extraits du déterminant caractéristique $|A-\lambda I|$ et d'une factorisation de l'équation caractéristique par un terme linéaire : $|A-\lambda I|/(\lambda-\lambda_i)$ (où $\lambda_i$ est une racine caractéristique de $A$). Une telle expression polynomiale s'identifie à un facteur près aux colonnes non nulles de la matrice des cofacteurs de $|A-\lambda I|$.

[5] Une étude de cette querelle est proposée dans une publication ([Brechenmacher 200?]) s'attachant au cas général de l'équivalence des couples de matrices en complément du cas symétrique détaillé dans cet article. Afin que ces deux publications, dont l'intersection est non vide, s'articulent tout en restant chacune autonome, les résultats de l'une ont été résumés dans l'autre. Ces travaux sont issus d'une thèse de doctorat, menée sous la direction de J. Dhombres à l'EHESS et au centre A. Koyré, posant plus généralement la question des différentes identités du théorème de Jordan sur la période 1870-1930 [Brechenmacher 2006a]. D'autres sources d'inspirations ont été les questionnements d'identités développés dans le travail de C. Gilain sur le théorème fondamental de l'algèbre [Gilain 1991], celui de H. Bos, C. Kers, F. Oort et D. W. Raven sur le théorème de Poncelet [Bos et al. 1987], de G. Cifoletti sur les pratiques algébriques de Pelletier et Gosselin dans le cadre de la « tradition algébrique française » de la renaissance [Cifoletti 1992] ainsi que l'étude des fluctuations d'élaborations mathématiques sur une longue période présentée dans l'ouvrage collectif coordonné par P. Benoit, K. Chemla et J. Ritter sur les fractions [Benoit et al. 1992] et le regard de C. Goldstein sur les relations entre analyse diophantienne et descente infinie [Goldstein 1993].

[6] Au sujet de Poincaré et du problème des trois corps, voir la récente thèse de doctorat d'A. Robadey [Robadey 2006] et, pour une perspective sur le long terme, [Laskar 1992]. Une histoire des pratiques associées à la représentation matricielle, notamment chez Weyr et Sylvester, est proposée dans [Brechenmacher 2006d].

# I. La controverse opposant Jordan et Kronecker en 1874, un moment de référence pour établir le corpus et organiser l'étude.

Au sein d'une série de notes et de mémoires publiés durant l'année 1874 Jordan et Kronecker opposent deux méthodes associées à deux théorèmes énoncés indépendamment et dans des théories différentes. L'un, dû à Weierstrass [1868] définit des *invariants*, les diviseurs élémentaires, pour caractériser l'équivalence des couples de formes bilinéaires ; l'autre, énoncé par Jordan en 1870, donne une réduction des substitutions linéaires à une *forme canonique* dans un contexte de travaux sur la résolubilité des équations algébriques ([8]). Comme nous allons le détailler dans cette partie, c'est au travers de la reconnaissance de leur capacité à résoudre de mêmes problèmes abordés dans le passé par des auteurs comme Lagrange ou Cauchy qu'une première identité entre ces deux théorèmes s'était manifestée de 1870 à 1873 .

[7] Le caractère « algébrique » de la mathématisation chez Lagrange a souvent été caractérisé par l'historiographie comme « formel » en raison de son détachement de l'observation physique et de la présentation de la mécanique comme une branche des mathématiques autonome de la géométrie. Le concept de formalisation associée au style de mathématisation et au modèle de théorie physique de Lagrange a notamment été étudié par C.Comte et W.Barroso-Filho qui ont mis en évidence « l'invention d'une nouvelle méthode mathématique à caractère algébrique, la plasticité du formalisme employé qui, simple et clair, donne une assez bonne description des phénomènes étudiés, la formulation logico déductive d'un seul principe unificateur, d'où toutes les propositions découlent » [Barroso-Filho et Comte 1988, p.131]. Les études sur la « formalisation de la mécanique » par Lagrange ont cependant surtout porté sur les pratiques du calcul des variations, les développements de Taylor et la méthode des coefficients indéterminés ([Fraser 1980 et 1985] parle du « style algébrique de Lagrange, [Dahan-Dalmedico 1990 et 1992] de « double mouvement de réduction de la Mécanique à l'Analyse, et de l'Analyse à l'Algèbre » et [Panza 1992 et 2003] de « réduction du calcul analytique à l'algèbre »), et non sur la pratique propre au problème des petites oscillations à laquelle est consacrée cet article.

[8] D'un point de vue qui nous est contemporain, le théorème de réduction d'une matrice à coefficients complexes à sa forme canonique de Jordan est équivalent au théorème des diviseurs élémentaires. Consulter par exemple le manuel de [Gantmacher 1959]. Ci-dessous, trois exemples de décompositions matricielles associées à la décomposition en diviseurs élémentaires d'un même polynôme caractéristique $|A-\lambda I|=(\lambda-1)^2(\lambda-2)^3(\lambda-3)$.

### 1. Une controverse opposant deux fins données à une histoire commune.

La querelle a pour origine l'ambition de Jordan de réorganiser la théorie des formes bilinéaires, développée dans les années 1860 par un petit groupe de géomètres berlinois sous la forme d'une théorie des invariants ([9]), par des méthodes de réduction à des formes canoniques:

> On sait qu'il existe une infinité de manières de ramener un polynôme bilinéaire
>
> $$P=\Sigma A_{\alpha\beta}x_\alpha y_\beta \ (\alpha=1, 2, ..., n, \beta=1, 2, ..., n)$$
>
> à la forme canonique $x_1y_1+...+x_my_m$, […] par des transformations linéaires opérées sur les deux systèmes de variables $x_1,..,x_n$, $y_1,..,y_n$. Parmi les diverses questions de ce genre que l'on peut se proposer, nous considérons les suivantes :
>
> 1. Ramener un polynôme bilinéaire $P$ à une forme canonique simple par des substitutions orthogonales opérées les unes sur $x_1,...,x_n$, les autres sur $y_1,...,y_n$.
> 2. Ramener $P$ à une forme canonique simple par des substitutions linéaires quelconques opérées simultanément sur les $x$ et les $y$.
> 3. Ramener simultanément à une forme canonique deux polynômes $P$ et $Q$ par des substitutions linéaires quelconques, opérées isolement sur chacune des deux séries de variables.
>
> Le premier de ces problèmes est nouveau, si nous ne nous trompons. Le deuxième a déjà été traité (dans le cas où $n$ est pair) par M. Kronecker, et le troisième par M. Weierstrass ; mais les solutions données par les éminents géomètres de Berlin sont incomplètes, en ce qu'ils ont laissé de côté certains cas exceptionnels qui, pourtant, ne manquent pas d'intérêt. Leur analyse est en outre assez difficile à suivre, surtout celle de M. Weierstrass. Les méthodes nouvelles que nous proposons sont, au contraire, extrêmement simples et ne comportent aucune exception. [Jordan 1873, p.1487].

Comme le font pressentir des applications récentes des formes bilinéaires à la géométrie (Klein, 1868), à l'arithmétique (Kronecker, 1874) et à divers problèmes d'intégrations de systèmes différentiels (Jordan, 1871-1872) et des équations de Fuchs (Hamburger, 1872, Jordan 1874), la théorie des formes bilinéaires jouera un rôle essentiel dans les mathématiques de la fin du XIX^e siècle ([10]). En 1874 cependant, la portée de la théorie se manifeste surtout par sa capacité à formuler dans un cadre théorique unique des problèmes issus de théories distinctes. Dans la citation de Jordan, la « forme canonique » $x_1y_1+...+x_my_m$ généralise la loi d'inertie de la théorie des formes quadratiques attachée notamment aux travaux d'Hermite des années 1850, le problème 1 fait référence à la classification des fonctions homogènes du second degré réalisée par Cauchy en 1829 dans un cadre géométrique, le problème 2 renvoie à la question arithmétique de l'équivalence des formes quadratiques dans la tradition de Carl Gauss et des *Disquitiones Aritmeticae* de 1801 et le problème 3 provient de la théorie des systèmes d'équations différentielles linéaires à coefficients

---

| | | | |
|---|---|---|---|
| *Formes de Jordan.* | $\begin{pmatrix}1&&&&&\\&1&&&&\\&&2&&&\\&&&2&&\\&&&&2&\\&&&&&3\end{pmatrix}$ | $\begin{pmatrix}1&&&&&\\&1&&&&\\&&2&1&&\\&&&2&1&\\&&&&2&\\&&&&&3\end{pmatrix}$ | $\begin{pmatrix}1&&&&&\\&1&&&&\\&&2&0&&\\&&&2&1&\\&&&&2&\\&&&&&3\end{pmatrix}$ |
| *Diviseurs élémentaires.* | *(λ-1), (λ-1), (λ-2), (λ-2),(λ-2), (λ-3)* | *(λ-1), (λ-1),* $(\lambda\text{-}2)^3$, *(λ-3)* | *(λ-1), (λ-1), (λ-2),* $(\lambda\text{-}2)^2$, *(λ-3)* |

[9] Des recherches sur la transformation des fonctions thêta de plusieurs variables sont à l'origine de la publication, en 1866, de deux mémoires de Christoffel et Kronecker qui revendiquent la création d'une théorie des formes bilinéaires. L'origine de la théorie des formes bilinéaires, à Berlin, dans les années 1860 est présentée dans [Hawkins 1977] et [Brechenmacher 2006a]. Au sujet des premiers travaux de Christoffel sur les formes consulter plus particulièrement [Mawhin 1981]. Pour des compléments mathématiques sur le problème de la réduction canonique des couples de matrices, voir [Dieudonné 1946].

[10] En des termes qui nous sont contemporains, la notion de forme bilinéaire joue pendant longtemps un rôle analogue à celui que jouera la notion de matrice dans l'algèbre linéaire du XX^e siècle.

constants remontant aux travaux de mécanique de Lagrange ([11]). La théorie des formes bilinéaires donne donc une identité nouvelle à des problèmes anciens.

Contrairement aux trois problèmes de réductions canoniques que distingue Jordan, Kronecker organise la théorie des formes bilinéaires autour d'un unique problème, celui de la caractérisation des classes d'équivalences de couples de formes bilinéaires (troisième problème dans la classification de Jordan), et de sa résolution par les invariants introduits par Weierstrass et son théorème des diviseurs élémentaires de 1868.

> Dans le Mémoire de M. *Jordan [...]*, la solution du premier problème n'est pas véritablement nouvelle ; la solution du deuxième est manquée, et celle du troisième n'est pas suffisamment établie. Ajoutons qu'en réalité ce troisième problème embrasse les deux autres comme cas particuliers, et que sa solution complète résulte du travail de M. *Weierstrass* de 1868 et se déduit aussi de mes additions à ce travail. Il y a donc, si je ne me trompe, de sérieux motifs pour contester à M. *Jordan* l'invention première de ses résultats, en tant qu'ils sont corrects [...]. [Kronecker 1874b, p.1181 (les italiques sont dans le texte original)].

Deux couples $(A,B)$ et $(A',B')$ sont dits équivalents s'il est possible les transformer l'un en l'autre par des substitutions linéaires. Comme l'ont montré Weierstrass et Kronecker dès 1868, ce problème donne une même formulation à différentes questions abordées dans le passé comme l'intégration de systèmes différentiels ($AY''=BY$), la détermination des axes principaux des coniques et quadriques ($AY=\lambda I$) ou la caractérisation des couples de formes quadratiques. Comme nous le verrons plus loin, c'est après avoir résolu à son tour ces différents problèmes entre 1871 et 1872 à l'aide de méthodes élaborées dans le cadre de la théorie des groupes, que Jordan prit connaissance de la théorie des formes bilinéaires élaborée à Berlin et démontra que « la réduction simultanée des deux fonctions $P$ et $Q$ est un problème identique à celui de la réduction d'une substitution linéaire à sa forme canonique » [Jordan 1873, p.1487]. La manifestation d'une identité entre la réduction canonique associée à *l'algèbre* des substitutions et le calcul des invariants attaché à la théorie *arithmétique* des formes, est donc indissociable d'une référence au passé. Plus précisément, la caractérisation d'un couple $(A, B)$ est abordée par l'étude de l'expression $A+sB$ - un polynôme de « polynômes bilinéaires » pour Jordan ou « faisceau de formes bilinéaires » pour Kronecker - et l'examen de l'équation polynomiale $|A+sB|=0$. Dans le cas où les racines $s_1, s_2, ..., s_n$ de l'équation caractéristique $|A+sB|=0$ sont toutes distinctes, le couple $(A, B)$ est caractérisé par les invariants que sont les racines $s_1, s_2, ..., s_n$ ou par la forme canonique $(x_1y_1+...+x_ny_n\ ,\ s_1x_1y_1+...+s_nx_my_m)$ ; l'occurrence de racines multiples nécessite la définition d'autres invariants ou formes canoniques ([12]). En donnant une résolution générale au problème, indépendamment de la multiplicité des racines, les théorèmes de Jordan et Weierstrass permettent de percevoir les travaux d'auteurs du passé comme lacunaires car limités au cas des racines distinctes. En ce sens, ces deux théorèmes se présentent comme deux fins données à une histoire commune. Pour donner corps à cette histoire à laquelle se réfèrent Jordan et Kronecker en 1874, nous avons mis en œuvre une recherche bibliographique basée sur les références données par ces deux auteurs et qui a permis, comme nous l'avons détaillé en introduction, de fixer le corpus de la discussion sur l'équation petites oscillations.

[11] La forme canonique $x_1y_1+...+x_my_m$ permet de déterminer les classes d'équivalence des matrices carrées pour la relation d'équivalence $(ARB \Leftrightarrow \exists P,Q \in GL_n(\hat{E}),\ PAQ=B)$. Le problème 1. consiste en l'étude de la relation de similitude des matrices orthogonales $(ARB \Leftrightarrow \exists P \in O(\hat{E}),\ P^{-1}AP=B)$. Le problème 2 renvoie à la congruence des matrices $(ARB \Leftrightarrow \exists P \in GL_n(\hat{E}),\ {}^tPAP=B)$. Le problème 3 à l'équivalence des couples de matrices $(A, B)$. Le problème 3 intervient pour la résolution des systèmes d'équations différentielles $AY''+BY=0$. Dans le cas particulier où $B=I$, la relation d'équivalence des couples $(A,I)$ est identique à la relation de similitude des matrices $B=P^{-1}AP$. Comme le fait remarquer Kronecker, le 3e problème suffit à déduire les deux autres : le problème 1. revenant à l'étude de la congruence du couple $(A,I)$ et le problème 2. à l'équivalence du couple $(A, {}^tA)$.

[12] Par exemple la forme bilinéaire $B(X,U)=ux+uy+vy$ n'a qu'une seule valeur propre 1 qui ne suffit pas à la caractériser car $B$ n'est pas la forme identité.

## 2. Les regards de Jordan et Yvon-Villarceau sur le corpus : une « incorrection » dans une pratique remontant à Lagrange.

Dans une note adressée à l'Académie de Paris en 1870, l'astronome Antoine Yvon-Villarceau appelle « l'attention des géomètres » sur un « point assez important de la théorie des équations linéaires », une « incorrection » dans la méthode « d'intégration des équations différentielles du mouvement de rotation d'un corps solide, soumis à l'action de la pesanteur […], présentée pour la première fois par l'illustre auteur de la *Mécanique analytique*, dans le cas des petites oscillations » d'une corde fixée en un point, lestée d'un nombre quelconque de masses et écartée de sa position d'équilibre. Le principe de conservation des forces vives permet de mathématiser le problème par un système d'équations différentielles linéaires à coefficients constants ([13]) :

> Lagrange forme trois équations différentielles du second ordre, entre lesquelles il élimine l'une des trois inconnues. Pour abréger j'écrirai le résultat de l'élimination comme il suit :
>
> $$(a)\begin{cases} g\dfrac{d^2u}{dt^2}+a\dfrac{d^2s}{dt^2}+cu=0, \\ f\dfrac{d^2s}{dt^2}+a\dfrac{d^2u}{dt^2}+cs=0, \end{cases}$$
>
> [Yvon-Villarceau 1870, p.763].

La méthode d'intégration du système *(a)* repose sur la détermination, par les méthodes d'éliminations, d'une équation algébrique :

> [Ces méthodes] fournissent l'équation caractéristique
>
> $$\frac{c^2}{\rho^4}-(f+g)\frac{c}{\rho^2}+fg-a^2=0$$

Aux racines $\rho$ et $\rho'$ de cette équation de degré 2 sont associées deux équations différentielles indépendantes, $\dfrac{d^2u}{dt^2}+\rho u=0$ et $\dfrac{d^2u}{dt^2}+\rho' u=0$, auxquelles est ramenée l'intégration du système.

> Faisant abstraction du signe des racines, et désignant leurs valeurs absolues par $\rho$ et $\rho'$, on a les expressions suivantes de $s$ et de $u$ :
>
> $$(f)\begin{cases} s=\alpha\sin(\rho t+\beta)+\alpha'\sin(\rho' t+\beta') \\ u=\dfrac{a\rho^2}{c-g\rho^2}\alpha\sin(\rho t+\beta)+\dfrac{a\rho'^2}{c-g\rho'^2}\alpha'\sin(\rho' t+\beta') \end{cases}$$

La réduction du système à deux équations indépendantes nécessite donc l'obtention de deux racines distinctes et Yvon-Villarceau critique l'interprétation donnée par Lagrange selon laquelle l'« inégalité » des racines serait garantie par l'hypothèse de la stabilité des petites oscillations. En cas de racines multiples, les oscillations ne resteraient pas bornées car le temps $t$ « sortirait du sinus » et les solutions prendraient la forme $s=t\sin(\rho t+\beta)$([14]).

> Au reste, dit Lagrange, comme cette solution est fondée sur l'hypothèse que $s$, $u$ et $\dfrac{d\theta}{dt}$ soient de très petites quantités il faudra, pour qu'elle soit légitime […] que les racines $\rho$ et $\rho'$ soient réelles et *inégales*, afin que l'angle *t soit toujours sous le signe des sinus.* C'est sur la seconde des conditions ici énoncées que je me permets d'appeler l'attention de l'Académie. Je dis qu'il n'est pas nécessaire que cette condition soit remplie, pour que les petites oscillations se maintiennent. […] Voici un cas très simple, auquel correspondent des racines égales de l'équation caractéristique : c'est celui d'un corps solide, homogène et de révolution, oscillant autour d'un point pris sur son axe de figure. Chacun comprendra sans recourir au calcul, que la petitesse des oscillations est assurée dans ce cas, si le centre de gravité est, à l'origine du mouvement, au-dessous du centre de suspension, à une petite distance de la verticale passant par ce point, et si le mouvement oscillatoire initial est suffisamment faible. [*Ibidem*, p.764-765].

[13] Dans la citation de Villarceau, $u$ et $s$ sont des fonctions de $t$, $g$, $f$ et $a$ des constantes. L'intervention de $u$ et $s$ est en miroir dans les équations, le système est donc symétrique. Nous verrons dans la troisième partie de cet article que cette propriété se dégage d'une pratique algébrique spécifique élaborée par Lagrange en 1766.

[14] On dirait aujourd'hui que le système est stable et se ramène à des « équations distinctes qui s'intègrent isolément » s'il est diagonalisable dans Ё. Une matrice dont les valeurs propres sont toutes distinctes est diagonalisable mais la condition n'est pas nécessaire, une matrice présentant des valeurs propres multiples peut être diagonalisable ou non diagonalisable. La matrice du système d'équations différentielles linéaire à coefficients constants de Lagrange est toujours diagonalisable car symétrique.

La critique d'Yvon-Villarceau est issue de préoccupations mécaniques concernant notamment l'application de la méthode des petites oscillations aux mouvements séculaires des planètes sur leurs orbites ([15]). Remettant en cause une pratique consistant à donner une représentation mécanique (les solutions ne sont pas stables) à une propriété algébrique (la multiplicité des racines) d'une équation particulière (l'équation caractéristique) ([16]), la critique de l'astronome concerne les « géomètres » et la « théorie des équations linéaires » car elle pose la question théorique de la caractérisation des systèmes se ramenant à des équations indépendantes.

> J'ai cru devoir appeler l'attention des géomètres sur un point assez important de la théorie des équations linéaires, et qui n'occupe pas une place suffisante dans les traités sur cette matière. Peut-être la question que je soulève a-t-elle déjà été résolue; mais il faut croire que la solution n'est pas généralement connue, puisque l'incorrection que je signale dans la *Mécanique analytique* a pu échapper à un géomètre aussi érudit que le savant auteur de la nouvelle édition d'un ouvrage devenu classique ([17]). [*Ibid.*, p.766].

Jordan répond à l'appel de Villarceau par la publication de deux notes aux *Comptes rendus* en 1871 et 1872. Il donne d'abord, en 1871, une méthode permettant l'intégration des systèmes différentiels à coefficients constants indépendamment de la multiplicité des racines : la résolution de Jordan, donnée en annexe 2, associe le système à une substitution linéaire qui est réduite à une forme canonique afin de donner au système une forme intégrable ([18]). Le géomètre parisien donne aussi une caractérisation des systèmes se « réduisant » à des « équations distinctes » de la forme $\frac{dy_1}{dt} = \sigma y_1$ par une relation entre le déterminant caractéristique et ses mineurs : chaque racine caractéristique $\sigma$ de multiplicité $\mu$ doit être une racine des mineurs d'ordre $\mu$-1 extraits du déterminant caractéristique ([19]). Dans une seconde note publiée en 1872 et intitulée « Sur les oscillations infiniment petites des systèmes matériels », Jordan répond plus précisément à la question posée par Villarceau en abordant le cas des systèmes intervenant en mécanique et démontre que le caractère *quadratique* de ces systèmes implique la relation entre déterminants et mineurs qui garantie la possibilité de les réduire à des « équations distinctes ». Ce faisant, Jordan dépasse le seul problème de l'intégration des systèmes symétriques en se référant au cadre plus général de la théorie des formes quadratiques en référence aux travaux de Cauchy et Hermite :

> Il est clair que la question de la réduction du système (1) à la forme canonique (7) est identique à ce problème connu : *Faire disparaître les angles des variables à la fois dans les deux formes quadratiques* T *et* U. [Jordan 1872, p.320].

Les propriétés par lesquelles Jordan caractérise, en réponse à l'appel de Villarceau, les systèmes réductibles à des équations distinctes sont en fait identiques à celles déjà énoncés par Weierstrass

---

15 Les travaux de mécanique céleste d'Yvon-Villarceau concernent la détermination des orbites de divers corps célestes. Des exemples célèbres sont le calcul de la périodicité de la comète d'Arrest (1851) et la prévision des éphémérides en tenant compte des perturbations produites par Jupiter. Pour une biographie d'Yvon-Villarceau, consulter [Baillaud 1957].

16 Pour Yvon-Villarceau, la forme $s = t\sin(pt+\beta)$ serait attribuée aux solutions en cas de racine double par une généralisation abusive de l'intégration d'une équation linéaire d'ordre $n$ à celle d'un système de $n$ équations.

> On sait que lorsque l'on a affaire à *une* équation linéaire à coefficients constants et que l'équation caractéristique présente des racines égales dont le degré de multiplicité est $m$, il faut multiplier le terme de l'intégrale correspondant à la racine multiple par un polynôme de degré $m$-1 par rapport à la variable indépendante ; or plusieurs auteurs semblent admettre la nécessité d'une modification analogue des termes correspondants à une racine multiple, dans le cas *d'un système* d'équations linéaires ; ces auteurs se bornent à renvoyer aux explications fournies à l'occasion d'*une* équation unique. [1870, p.766].

17 Yvon-Villarceau fait ici allusion à la nouvelle édition de la *Mécanique Analytique* de Lagrange par J. Bertrand en 1853.

18 Les années 1870-1880 sont une période charnière dans la carrière de Jordan marquée par une entreprise d'application à des domaines variés des notions et méthodes développées pour la théorie des substitutions dans les années 1860. Voir à ce sujet [Brechenmacher 2006b].

19 Cette propriété caractérise les endomorphismes diagonalisables par la condition que leurs diviseurs élémentaires sont simples. Le degré de multiplicité $\mu$ d'une valeur propre doit coïncider avec la dimension du sous espace propre associé.

en 1858 (cas symétrique des couples de formes quadratiques) et 1868 (cas général des couples de formes bilinéaires) ([20]). Cette identité est d'abord relevée par Meyer Hamburger qui, en 1873, s'inspire des méthodes de Jordan et Weierstrass pour résoudre le problème posé par l'occurrence de racines caractéristiques multiples dans l'intégration des équations différentielles linéaires à coefficients non constants, dites équations de Fuchs. La parution du mémoire d'Hamburger attire l'attention de Jordan sur la théorie des formes bilinéaires développée à Berlin et sur laquelle ce dernier va publier une série de notes entre 1873 et 1874. Comme nous l'avons vu dans le premier paragraphe, Jordan démontre que la transformation des formes peut s'interpréter comme l'action de groupes de substitutions (linéaires, orthogonales etc.) et, en élaborant une identité mathématique entre ses méthodes de réductions canoniques des substitutions et les calculs d'invariants de Kronecker et Weierstrass, il propose une nouvelle organisation de la théorie des formes bilinéaires en la présentant comme une application de la théorie des groupes structurée dans son *Traité des substitutions et des équations algébriques* de 1870.

## 3. Le regard de Kronecker sur le corpus : le théorème de Weierstrass perçu comme une rupture dans l'histoire des raisonnements algébriques.

Comme nous l'avons vu, de premières identités entre les théorèmes de Jordan et de Weierstrass se sont manifestées par la capacité de chacun des théorèmes à donner une résolution générale à de mêmes problèmes dont le traitement par des auteurs du passé apparaît par conséquent lacunaire. En 1874, Kronecker associe au théorème de Weierstrass un idéal de généralité qu'il va expliciter par une critique des pratiques algébriques du passé et sur lequel il va s'appuyer pour dénoncer le caractère formel qu'il prête à la réduction canonique de Jordan. A partir d'un « défaut » relevé dans un calcul de Jordan, la mise au dénominateur d'une expression algébrique susceptible de s'annuler ([21]), Kronecker dénonce un certain type de pratiques algébriques qui s'appuient sur des expressions « prétendument générales » mais en réalité formelles car perdant toute signification dans certains « cas singuliers » au contraire de la « vraie généralité » qu'il dépeint comme un océan dont il faut percer la « surface » (la prétendue généralité) pour pénétrer les « profondeurs » en résolvant une « question algébrique dans ses moindres détails ». Comme le montre l'extrait ci-dessous, cette opposition entre deux types de « généralité » va orienter un regard porté sur l'histoire de « la théorie des faisceaux des formes quadratiques »:

> Ceci se confirme partout dans les rares questions algébriques qui sont mises en œuvre complètement jusqu'à leurs moindres détails, notamment dans la théorie des faisceaux des formes quadratiques […]. Parce que, pendant si longtemps, on n'osait pas faire tomber la condition que le déterminant ne contient que des facteurs inégaux, on est arrivé avec cette question connue de la transformation simultanée de deux formes quadratiques; qui a été si souvent traitée depuis un siècle, mais de manière sporadique, à des résultats très insuffisants et les vrais aspects de l'étude ont été ignorés. Avec l'abandon de cette condition, le travail de Weierstrass de l'année 1858 à conduit à un aperçu plus élevé et notamment à un règlement complet du cas, dans lequel n'existent que des diviseurs élémentaires simples. Mais l'introduction générale de cette notion de diviseur élémentaire, dont seule une étape provisoire était alors accomplie, intervient seulement dans le mémoire de Weierstrass de l'année 1868, et une lumière tout à fait nouvelle est ainsi faite sur la théorie des faisceaux pour n'importe quel cas, avec la seule condition que le déterminant soit différent de zéro. Quand j'ai aussi dépouillé cette dernière restriction et l'ai développé à partir de la notion de diviseur élémentaire des faisceaux élémentaires généraux, la clarté la plus pleine s'est répandue sur une quantité de nouvelles formes algébriques, et par ce traitement complet de l'objet des vues plus élevées ont été acquises sur une théorie des invariants comprise dans sa vraie généralité [Kronecker 1874c, p. 404, *traduction F.B.*].

[20] Il semble que les travaux de Weierstrass de 1858 aient influencé Jordan de manière indirecte, par l'intermédiaire d'un mémoire publié à l'Académie de Saint Petersbourg par Somof [1859]. Ce mémoire, porté à l'attention de Jordan par Yvon-Villarceau en 1872, contient un exposé de la méthode de Weierstrass appliquée au problème des petites oscillations.

[21] Ce défaut, immédiatement corrigé par Jordan, est sans conséquence sur l'organisation théorique que celui-ci propose. Voir à ce sujet [Brechenmacher 2006a].

Par son traitement général et homogène de la caractérisation des couples de formes bilinéaires, le théorème de Weierstrass vient sanctionner les « résultats très insuffisants » de pratiques développées « durant tout un siècle » mais « de manière sporadique » et négligeant d'aborder les singularités en « n'osant pas faire tomber » la condition que le déterminant $S=|A+sB|$ contient des facteurs inégaux ([22]). Faisant implicitement référence aux travaux d'auteurs comme Lagrange, Charles Sturm, Cauchy ou Carl Gustav Jacobi, ces pratiques, que nous préciserons dans les troisième et quatrième parties de cet article, avaient été développées pour le cas symétrique des « systèmes quadratiques » dont elles exprimaient les solutions $x_i$ par des expressions polynomiales « générales » données en fonction de factorisations de $S$ et de ses sous déterminants $P_{1i}$ obtenus par développements par rapport à la $1^{re}$ ligne et $i^{e}$ colonne. L'expression générale (*) $\dfrac{P_{1i}}{\dfrac{S}{x-s_j}}(x)$ donnant $x_i^{s_j}=\dfrac{P_{1i}}{\dfrac{S}{x-s_j}}(s_j)$ où $(x_i^{s_j})_{1\leq i\leq n}$ désigne le système de solutions associé à la racine $s_j$ ([23]). Kronecker condamne l'utilisation qui avait été faite de telles expressions, susceptibles de prendre une valeur du type $\dfrac{0}{0}$ en cas d'occurrence de racines communes entre les équations obtenues par les sous déterminants successifs de $S=0$, avant que Weierstrass ne démontre en 1858,

---

[22] C'est-à-dire la condition imposant des racines caractéristiques distinctes ou, en termes contemporains, des matrices diagonalisables.

[23] On dirait aujourd'hui que les coordonnées des vecteurs propres d'une matrice $A$ sont données par les colonnes non nulles de la matrice adjointe de la matrice caractéristique $A-xI$, c'est-à-dire la matrice des cofacteurs. Par exemple, pour la matrice

$$A=\begin{pmatrix}1 & -1 & 0\\ -1 & 2 & 1\\ 0 & 1 & 1\end{pmatrix}$$

d'équation caractéristique et de mineurs : $S=-x(3-x)(1-x)$, $P_{11}(x)=(1-x)(2-x)-1$, $P_{12}(x)=(1-x)$ et $P_{13}=-1$.
L'expression $((1-x)(2-x)-1,\ 1-x,\ -1)$ donne une écriture polynomiale des coordonnées d'un vecteur propre dont $|S|$ donne le carré de la norme. Pour la valeur propre $s_1=1$, $\dfrac{S}{x-1}=x(3-x)$. Par conséquent,

$$x_1^{s_1}=\frac{P_{11}}{\dfrac{S}{x-1}}(1)=-\frac{1}{2},\ x_2^{s_1}=\frac{P_{12}}{\dfrac{S}{x-1}}(1)=0,\ x_3^{s_1}=\frac{P_{11}}{\dfrac{S}{x-1}}(1)=-\frac{1}{2},$$

On en déduit les coordonnées d'un vecteur propre normé associé à la valeur propre 1 : $(1/\sqrt{2},\ 0,\ 1/\sqrt{2})$.
En procédant de manière analogue pour les valeurs propres $s_2=0$,

$$x_1^{s_2}=\frac{1}{3},\ x_2^{s_2}=\frac{1}{3},\ x_3^{s_2}=-\frac{1}{3}$$

et $s_3=3$,

$$x_1^{s_3}=\frac{1}{6},\ x_2^{s_3}=\frac{2}{6},\ x_3^{s_3}=-\frac{1}{6}$$

on obtient respectivement pour vecteurs propres normés, $(1/\sqrt{3},\ 1/\sqrt{3},\ -1/\sqrt{3})$ et $(1/\sqrt{6},\ 2/\sqrt{6},\ -1/\sqrt{6})$.
Interprétation dans le cadre des formes quadratiques : la forme associée à $A$ dans la base canonique de Ë$^3$,

$$A(x_1,x_2,x_3)=x_1^2-2x_1x_2+x_2^2+2x_2x_3+x_3^2=1.X_1^2+0.X_2^2+3.X_3^2$$

avec $(X_1,X_2,X_3)=(x_i^{s_j})_{i,j=1,2,3}A(x_i^{s_j})^{-1}_{i,j=1,2,3}$ où $(x_i^{s_j})_{i,j=1,2,3}$ est la matrice orthogonale de passage (ou de changement de base orthonormées) :

$$(x_i^{s_j})_{i,j=1,2,3}=\begin{pmatrix}\frac{1}{\sqrt{2}} & \frac{1}{\sqrt{3}} & \frac{1}{\sqrt{6}}\\ 0 & \frac{1}{\sqrt{3}} & \frac{2}{\sqrt{6}}\\ \frac{1}{\sqrt{2}} & -\frac{1}{\sqrt{3}} & -\frac{1}{\sqrt{6}}\end{pmatrix}$$

Interprétation dans le cadre des formes bilinéaires (symétriques) :

$$A(x,y)=x_1y_1-x_1y_2-y_1x_2+x_2y_2+x_2y_3+y_3x_2+x_3y_3=1.X_1Y_1+0.X_2Y_2+3.X_3Y_3.$$

comme nous le détaillerons plus loin, que tout couple ( $A=\sum_{i=1}^{n} A_{ij}x_i x_j$ $B=\sum_{i=1}^{n} B_{ij}x_i x_j$ ) de formes quadratiques - avec $A$ définie positive -, peut être écrit sous la forme ( $\sum_{i=1}^{n} X_i^2, \sum_{i=1}^{n} s_i X_i^2$ ) où les racines $s_1, s_2, ..., s_n$ de $S=0$ sont toujours réelles, quelle que soit leurs multiplicités, et ne prouve par la même occasion que toute racine d'ordre $p$ de $S$ est une racine d'ordre supérieur ou égal à $p-1$ pour $P_{1i}$, l'expression (*) étant par conséquent toujours parfaitement définie ([24]). Le théorème de 1858 généralise la loi d'inertie selon laquelle, comme l'énonçait Hermite, une forme quadratique homogène à $n$ variables $A$ peut s'écrire comme une somme de carrés de la forme

$$A=\Delta_{n-1}X_1^2+\frac{\Delta_{n-2}}{\Delta_{n-1}}X_2^2+...+\frac{\Delta}{\Delta_1}X_n^2,$$

où $\Delta, \Delta_1, \Delta_2, ..., \Delta_{n-1}, 1$ sont les mineurs principaux successifs du déterminant de $A$. Comme le formule explicitement Darboux en 1874, cette loi se généralise aux couples $(A, I)$ envisagés comme des formes quadratiques polynomiales $A+sI$, les quotients de mineurs successifs qu'elle met en œuvre sont alors toujours définis car « si elle [l'équation caractéristique] a des racines multiples, une racine d'ordre $p$ devra annuler tous les mineurs d'ordre $p-1$ de l'invariant » [Darboux 1874, p.367]. Le mémoire de Darboux, dont un extrait est donné en annexe 3, est publié dans le *Journal de Mathématiques* à la suite du mémoire de Jordan [1874] qui est au cœur de la controverse avec Kronecker ; les deux mémoires se complètent et tous deux concernent des théorèmes de Weierstrass, d'une part Jordan propose une nouvelle théorie des faisceaux de formes bilinéaires abordés par Weierstrass en 1868, d'autre part Darboux propose une nouvelle démonstration du résultat de 1858 sur les faisceaux de formes quadratiques ([25]).

Déjà en germe dans le résultat de 1858 sur les formes quadratiques qui montre que la question essentielle n'est pas la multiplicité des racines mais les relations de divisibilité entre le déterminant caractéristique et ses mineurs successifs, le théorème des diviseurs élémentaires est considéré comme représentant une *rupture* dans l'histoire des raisonnements algébriques. Ce regard historique de Kronecker n'est pas indépendant de celui porté en 1977 par Hawkins sur l'histoire de la théorie spectrale des matrices. Dans ses travaux consacrés à l'origine de cette théorie, Hawkins a mis en évidence le rôle des catégories de généralité et de rigueur dans l'évolution de l'algèbre aux XVIII^e^ et XIX^e^ siècles.

> Dès l'origine des méthodes de l'analyse [chez Viète], l'importance de la valeur attribuée à ces méthodes était indissociable de leur généralité. Ainsi Viète présentait-il la nouvelle méthode de l'analyse comme « un art particulier, qui n'exerce plus sa raison par les nombres, ce qui se faisait auparavant par la nonchalance des anciens analystes, mais par une logistique des espèces […] » [Viète 1591, p.321-322]. […]. L'analyse devint une méthode pour raisonner en manipulant des expressions comportant des symboles prenant des valeurs « générales » et une tendance se développa à penser exclusivement en termes de cas « généraux » avec peu, ou aucune, attention pour les difficultés potentielles ou inexactitudes pouvant apparaître en assignant des valeurs particulières aux symboles. Je qualifierai de tels raisonnements basés sur des expressions « générales » de « raisonnements génériques ». Le traitement par Lagrange du théorème des axes principaux de la mécanique [1788] donne un exemple simple de raisonnement générique […]. Durant la première partie du XIX^e^ siècle, certains mathématiciens, préoccupés par une élévation des critères de rigueur en analyse, remirent en cause la légitimité du raisonnement générique. Ainsi Cauchy, mettait-il en avant, dans sa préface des *Cours d'analyse* [1821], que le raisonnement mathématique ne devait pas se baser sur des arguments « tirés de la généralité de l'algèbre », arguments qui tendent « à faire attribuer aux formules algébriques une étendue indéfinie, tandis que

[24] Contrairement au cas bilinéaire, dans le cas quadratique étudié par Weierstrass en 1858 les matrices sont symétriques donc toujours diagonalisables sur Ë. Une matrice est diagonalisable si et seulement si ses diviseurs élémentaires sont simples.

[25] Le mémoire de Darboux concerne avant tout les surfaces, et sa méthode est insérée par Gundelfinger dans la troisième édition de la géométrie analytique de Hesse. Pour une description détaillée du travail de Darboux, voir l'article de Drach et Meyer dans *L'encyclopédie des sciences mathématiques* de 1907.

> […] la plupart de ces formules subsistent uniquement sous certaines conditions ». Weierstrass endossa le programme de Cauchy visant à rendre l'analyse plus rigoureuse. Les deux contributions à la théorie des matrices qu'il publia [1858, 1868] étaient en effet une conséquence de ce programme et des insuffisances du raisonnement générique. [Hawkins 1977, p.122, *traduction F.B.*].

Le qualificatif de « générique », donné par Hawkins à un certain type de raisonnement qui ne se préoccupe pas de la signification des symboles, structure une histoire à long terme de l'algèbre, de son origine dans l'analyse de François Viète au développement d'exigences de rigueur dans les « raisonnements mathématiques » caractéristiques du XIX^e^ siècle et personnifiés par des auteurs comme Cauchy et Weierstrass. Les travaux des auteurs du corpus de l*a discussion* sont alors envisagés sous l'angle d'une opposition entre « raisonnement générique » et « rigueur » dans la manière d'appréhender la question des racines caractéristiques multiples.

> En raison de la nature bornée des solutions [des problèmes de petites oscillations], Lagrange et Laplace affirmaient que les racines caractéristiques $\lambda$ devaient être réelles. Aucun des deux savants n'avait cependant poursuivi l'étude des solutions des systèmes d'équations différentielles avec suffisamment de soin pour justifier ses affirmations. S'ils n'avaient aucune difficulté à étudier les systèmes lorsque les racines caractéristiques étaient distinctes, leurs analyses du cas des racines multiples étaient inadéquates. Compte tenu de la tendance générique de leurs méthodes analytiques, il est déjà remarquable qu'ils aient mentionné ce cas de figure. […] La reconnaissance par Weierstrass du caractère discutable de leurs affirmations était à l'origine de la recherche qui allait culminer avec sa théorie des diviseurs élémentaires. [….] En [1868] et dans un mémoire préliminaire [1858], Weierstrass démontra plus que des théorèmes. Il démontra également la possibilité d'une approche plus rigoureuse de l'analyse algébrique qui ne se satisferait pas de la tendance dominante à raisonner vaguement dans le cadre du cas « général ». [*Ibidem*, p.123-124, *traduction F.B.*].

Le théorème des diviseurs élémentaires de Weierstrass est considéré comme représentant une rupture dans l'histoire des raisonnements algébriques par son caractère à la fois rigoureux, par opposition aux raisonnements génériques du passé, et homogène, par opposition à une perception des différents cas de multiplicité de racines comme autant de cas singuliers nécessitant des traitements spécifiques. Comme l'a montré Hawkins, cette rupture se manifeste notamment dans l'évolution des travaux de Kronecker. Dans un mémoire consacré aux fonctions thêta en 1866, ce dernier avait été amené à considérer l'occurrence de racines caractéristiques multiples comme un cas singulier nécessitant un traitement spécifique distinct de la méthode algébrique développée pour le « cas général ». Par son caractère homogène, le système complet d'invariant introduit par Weierstrass en 1868 démontre la possibilité de concilier une approche algébrique et un traitement général du problème. Il est à l'origine de nouveaux idéaux de généralité dans les raisonnements algébriques qui se manifestent notamment dans l'élaboration de la théorie des formes bilinéaires à Berlin et dans lesquels Hawkins voit une origine de la « théorie des matrices ».

> J'aimerais suggérer que, si un seul auteur devait mériter le titre de fondateur de la théorie des matrices, cet individu serait Weierstrass. […] Sa théorie des diviseurs élémentaires constitua un cœur théorique, une fondation sur laquelle bâtir. Son œuvre démontra la possibilité d'employer les méthodes de l'analyse de manière non générique, ouvrant par là un nouveau monde à la recherche mathématique, un monde dont ses collèges et étudiants entreprirent l'exploration. […] L'émergence de la théorie [spectrale] des matrices était directement liée à la chute des approches génériques de l'analyse algébrique. [*Ibid.*, p.157-159, *traduction F.B.*].

Nous avons vu que les premières lignes de l'histoire qui qualifie de « génériques » des « raisonnements algébriques » des XVIII^e^ et XIX^e^ siècles sont écrites par Kronecker en 1874. Cette perception de l'histoire s'avère donc indissociable de l'énoncé du théorème des diviseurs élémentaires qu'elle présente comme une rupture et que Kronecker oppose au théorème de réduction canonique de Jordan ([26]). Lorsqu'il attirait l'attention de l'Académie sur une « lacune » dans la résolution du problème des petites oscillations par Lagrange, Yvon-Villarceau ne critiquait pas le caractère générique d'un raisonnement mais une pratique donnant une représentation mécanique erronée à la nature algébrique des racines de l'équation caractéristique. La reconnaissance du contexte très particulier dans lequel Kronecker sanctionne les pratiques algébriques du passé ainsi que des regards différents portés par Yvon-Villarceau et Jordan sur ces mêmes pratiques nous amène à développer une approche complémentaire à celle suivie par Hawkins lorsqu'il avait mis en évidence le caractère générique de *raisonnements* algébriques en nous concentrant sur les différentes représentations associées à une même *pratique* algébrique dans divers cadres théoriques comme la mécanique, la géométrie ou l'arithmétique avant que cette pratique ne revête l'identité d'une *méthode* dans un cadre théorique unique comme celui de la théorie des formes bilinéaires. En 1874, la nature, algébrique ou arithmétique, d'une telle théorie oppose Jordan et Kronecker qui se réfèrent pourtant tous deux à un même corpus de travaux du passé dans lequel ils identifient une pratique spécifique qui leur est commune et consistant à aborder les transformations des couples de formes bilinéaires $(A,B)$ par la décomposition – algébrique ou arithmétique- de la *forme* polynomiale de l'équation $|A+sB|=0$. Le fait que, à la suite de l'appel d'Yvon-Villarceau, Jordan résolve entre 1871 et 1872 des problèmes différents attachés à des auteurs du corpus comme Lagrange, Cauchy et Hermite indique qu'aux yeux de ce dernier le corpus manifeste une cohérence interne indépendante de l'identité théorique donnée par Weierstrass qui, dès 1858 concevait l'ensemble des problèmes du corpus comme relevant de la « théorie des couples de fonctions homogènes ». La deuxième partie de notre article sera consacrée à expliciter de tels éléments de cohérence interne qui font l'identité de notre corpus au-delà des cadres théoriques au travers desquels celui-ci se déploie sur plus d'un siècle, de son origine chez Lagrange aux deux fins qu'en donnent Weierstrass en 1858 et 1868 et Jordan en 1871-1872.

[26] Sur la question de la fabrication de l'histoire par les textes mathématiques, voir la discussion de référence de [Goldstein 1995], celle de [Dhombres 1998] en relation avec le concept de postérité, les exemples donnés par [Cifoletti 1992 et 1995] et, pour notre propos, le développement proposé dans [Brechenmacher 2006c].

# II. Les identités de la *discussion* sur l'équation des petites oscillations.

La recherche bibliographique menée à partir du moment de référence que nous avons détaillé dans la partie précédente permet d'obtenir un corpus qui, bien que s'étendant sur une longue période, reste fermé et limité à un nombre restreint d'auteurs. Ce caractère fermé donne au corpus l'identité d'une discussion présentant des caractéristiques propres que nous allons préciser en examinant les principaux nœuds qui apparaissent dans l'enchevêtrement des références bibliographiques et sont attachés aux travaux de Lagrange, Pierre Simon Laplace, Cauchy et Weierstrass. En portant un premier regard d'ensemble sur notre corpus, nous nous attacherons plus particulièrement à la manière dont les auteurs se font référence les uns aux autres.

## 1. Une référence commune à la *Mécanique Analytique* de Lagrange.

La référence systématique à la résolution par Lagrange du « problème des oscillations très petites d'un système quelconque de corps » est une première caractéristique du corpus. Il s'agit d'étudier le déplacement d'un système de corps *m, m', m'', &c,* de coordonnées *x, y, z, ...*, dont on suppose qu'il s'éloigne très peu d'une position d'équilibre de coordonnées *a, b, c, ...* :

> Au reste, on pourra souvent aussi, en ayant égard aux conditions du problême, réduire les coordonnées immédiatement par des substitutions, en fonctions rationnelles & entières d'autres variables indépendantes entr'elles, & très petites, dont la valeur soit nulle dans l'état d'équilibre.
> Ainsi nous supposerons en général que l'on ait :
> $$x = a + a_1\xi + a_2\psi + a_3\varphi + \&c + a'_1\xi^2 + \&c,$$
> $$y = b + b_1\xi + b_2\psi + b_3\varphi + \&c + b'_1\xi^2 + \&c,$$
> $$z = c + c_1\xi + c_2\psi + c_3\varphi + \&c + c'_1\xi^2 + \&c,$$
> & ainsi des autres coordonnées *x', y',* &c, les quantités $a, b, c, a_1, b_1$ & c, sont constantes, & les quantités ξ, ψ, φ, &c, sont variables, très petites, & nulles dans l'équilibre. [Lagrange 1788, p.243].

Les équations du mouvement sont déduites des conditions initiales et du principe de conservation des forces vives ([27]). Elles se présentent sous « une forme linéaire avec des coefficients constants » :

$$0 = (I)\frac{d^2\xi}{dt^2} + (I,2)\frac{d^2\psi}{dt^2} + (I,3)\frac{d^2\varphi}{dt^2} + \&c + [I]\xi + [I,2]\psi + [I,3]\varphi + \&c.$$

$$0 = (2)\frac{d^2\psi}{dt^2} + (I,2)\frac{d^2\xi}{dt^2} + (2,3)\frac{d^2\varphi}{dt^2} + \&c + [2]\psi + [I,2]\xi + [2,3]\varphi + \&c.$$

$$0 = (3)\frac{d^2\varphi}{dt^2} + (I,3)\frac{d^2\xi}{dt^2} + (2,3)\frac{d^2\psi}{dt^2} + \&c + [3]\varphi + [I,3]\xi + [2,3]\psi + \&c.$$

*&c.*

> équations qui étans sous une forme linéaire avec des cœfficients constants, peuvent être intégrées rigoureusement & généralement par les méthodes connues [*Ibidem*].

L'intégration du système est basée sur le principe de réduction de l'ordre d'une équation différentielle par la connaissance de solutions particulières. Trouver *n* solutions particulières $\xi = E\sin(t\sqrt{K} + \varepsilon)$ de *n* équations indépendantes $\frac{d^2\xi}{dt^2} + K\xi = 0$ permet, « en les joignant ensemble », d'exprimer les solutions générales ([28]) :

> [En supposant que] les variables dans ces sortes d'équations ayent entr'elles des rapports constants; c'est-à-dire que l'on ait $\psi = f\xi$, $\phi = g\xi$, *&c*, le système s'écrit alors :
>
> $$((I) + (I,2)f + (I,3)g + \&c)\frac{d^2\xi}{dt^2} + ([I] + [I,2]f + [I,3]g + \&c)\xi = 0,$$
>
> $$((2)f + (I,2) + (2,3)g + \&c)\frac{d^2\xi}{dt^2} + ([2]f + [I,2] + [2,3]g + \&c)\xi = 0,$$
>
> $$((3)g + (I,3) + (2,3)f + \&c)\frac{d^2\xi}{dt^2} + ([3]g + [I,3] + [2,3]f + \&c)\xi = 0,$$
>
> *&c,*
>
> lesquelles donnent $\frac{d^2\xi}{dt^2} + K\xi = 0$, en faisant
>
> $$K = \frac{[I]+[I,2]f+[I,3]g+\&c}{(I)+(I,2)f+(I,3)g+\&c} = \frac{[2]f+[I,2]+[2,3]g+\&c}{(2)f+(I,2)+(2,3)g+\&c} = \frac{[3]g+[I,3]+[2,3]f+\&c}{(3)g+(I,3)+(2,3)f+\&c}$$
>
> Maintenant l'équation $\frac{d^2\xi}{dt^2} + K\xi = 0$, donne par l'intégration
>
> $$\xi = E\sin(t\sqrt{K} + \varepsilon).$$
>
> *E*, *ε* étant des constantes arbitraires ; ainsi comme on a supposé $\psi = f\xi$, $\varphi = g\,\xi$, &c, on a aussi les valeurs de *ψ*, *φ*, &c. Cette solution n'est que particulière, mais elle est en même temps double, triple, &c, selon le nombre des valeurs de *K* ; par conséquent en les joignant ensemble, on aura la solution générale […]. Dénotant par *K', K'', K''', &c,* les différentes valeurs de *K*, c'est-à-dire, les racines de l'équation en *K* […]
>
> $$\xi = E'\sin(t\sqrt{K'}+\varepsilon') + E''\sin(t\sqrt{K''}+\varepsilon'') + E'''\sin(t\sqrt{K'''}+\varepsilon''') + \&c,$$
>
> $$\psi = f'E'\sin(t\sqrt{K'}+\varepsilon') + f''E''\sin(t\sqrt{K''}+\varepsilon'') + f'''E'''\sin(t\sqrt{K'''}+\varepsilon''') + \&c,$$
>
> $$\varphi = g'E'\sin(t\sqrt{K'}+\varepsilon') + g''E''\sin(t\sqrt{K''}+\varepsilon'') + g'''E'''\sin(t\sqrt{K'''}+\varepsilon''') + \&c,$$
>
> *&c,*
>
> [*Ibid.*, p.244].

Les valeurs de *K* apparaissent, par usage des méthodes d'élimination, comme les racines d'une équation résultante de degré *n* :

---

[27] Il s'agit du calcul du Lagrangien $L=T-V$ où *T* représente l'énergie cinétique et *V* l'énergie potentielle. Le système est obtenu après prise en compte des conditions initiales et développements de Taylor de *V* dont on néglige les termes supérieurs à l'ordre 2 en *t* selon le principe des petites oscillations.

[28] Dire que les solutions ont entre elles des rapports constants ($\psi = f\xi$, $\varphi = g\,\xi$,…) et sont exprimables à l'aide de *n* solutions particulières indépendantes s'interprète, dans le cadre des mathématiques contemporaines, par la structure linéaire de l'ensemble des solutions formant un espace vectoriel de dimension *n*.

Le nombre de ces équations est, comme l'on voit, égal à celui des inconnues *f, g, &c, K* ; par conséquent elles déterminent exactement ces inconnues ; & comme en retenant pour premier membre le terme *K*, le multipliant respectivement par le dénominateur du second, on a des équations linéaires en *f, g, &c*, il sera facile de les éliminer par les méthodes connues, & il n'est pas difficile de voir par les formules générales d'élimination, que la résultante en *K* sera d'un degré égal à celui des équations, & par conséquent égal à celui des équations différentielles proposées ; de sorte que l'on aura pour *K* un pareil nombre de différentes valeurs, dont chacune étant substituée dans les expressions de *f, g, &c*, donnera les valeurs correspondantes de ces quantités [*Ibidem*, p.245].

Le problème se ramène donc à la résolution d'une équation algébrique caractéristique du problème des petites oscillations. Pour Lagrange, la *nature mécanique des oscillations*, supposées petites dans le problème posé, implique en effet un caractère spécifique de la *nature algébrique des racines* de l'équation, celles-ci doivent être réelles et inégales :

De là on tire une méthode générale pour voir si l'état d'équilibre d'un système quelconque donné de corps est stable, c'est-à-dire si, les corps étant infiniment peu dérangés de cet état, ils y reviendront d'eux-mêmes, ou au moins tendront à y revenir. […]

1° Si toutes les racines de cette équation sont réelles négatives et inégales, l'état d'équilibre sera stable en général, quel que soit le dérangement initial du système.

2° Si ces racines sont toutes réelles, positives ou toutes imaginaires ou en partie positives, et en partie imaginaires, l'état d'équilibre n'aura aucune stabilité, et le système une fois dérangé de cet état ne pourra le reprendre ;

3° Enfin, si les racines sont en partie réelles négatives et inégales, et en partie réelles négatives et égales ou réelles et positives, ou imaginaires, l'état d'équilibre aura seulement une stabilité restreinte et conditionnelle. [Lagrange 1766, p.532].

## 2. Un problème considéré comme résolu par Lagrange pour la durée de la discussion.

La « méthode générale » donnée par Lagrange pour caractériser la stabilité des systèmes mécaniques ne sera pas remise en cause avant la publication par Weierstrass, en 1858, d'un mémoire intitulé « Ueber ein die homogenen functionen zweiten grades betreffendes Theorem, nebst anwendung desselben auf die theorie der kleinen schwingungen ».

Après avoir indiqué et énoncé la forme des intégrales, Lagrange a conclu que, comme les oscillations $x_1$, $\frac{dx_1}{dt}$ restent toujours petites si elles le sont à l'origine, l'équation ne peut pas avoir de racines égales car les intégrales pourraient devenir arbitrairement grandes avec le temps. La même affirmation se trouve répétée chez Laplace lorsqu'il traite dans la *Mécanique céleste* des variations séculaires des planètes. Beaucoup d'autres auteurs, comme, par exemple, Poisson, mentionnent cette même conclusion. Mais cette conclusion n'est pas fondée […] et, si la fonction $\Psi$ reste négative et de déterminant non nul, on peut énoncer le même résultat, que les racines de l'équation $f(s)=0$ soient ou non toutes distinctes ; l'homogénéité de cette conclusion n'a pu être découverte dans la passé car on a toujours envisagé ce cas [des racines multiples] par des approches particulières. [Weierstrass 1858, p. 244, *traduction F.B.*]

La « conclusion homogène » de Weierstrass est fondée sur l'énoncé d'un théorème portant sur les couples $(\Phi,\Psi)$ de fonctions homogènes du second degré et s'inscrivant dans la théorie des formes quadratiques développée dans les années 1850 par Jacobi, Hermite, Dirichlet et Sylvester ([29]). Dans le cas où la forme $\Phi$ est définie, positive et où les coefficients de $\Phi$ et $\Psi$ sont réels, le théorème de Weierstrass énonce que les racines $s_i$ de l'équation $f(s)=0$ sont toujours réelles indépendamment de leur multiplicité. La stabilité des systèmes mécaniques est donc assurée par leur nature quadratique. Mais comme le montre l'appel d'Yvon-Villarceau qui, comme nous l'avons vu dans la première partie, attire en 1870 l'attention des géomètres de l'Académie sur une

[29] La forme $\Phi$ étant supposée de déterminant non nul. La relation entre les travaux de Weierstrass et l'arithmétique des formes quadratiques sera développée en conclusion de cet article.

« lacune » dans la relation établie par Lagrange entre stabilité mécanique et multiplicité des racines caractéristiques, le problème des petites oscillations est considéré comme résolu dans la période qui sépare Lagrange des théorèmes de Weierstrass [1858, 1868] et de Jordan [1871, 1872]. Notre corpus ne peut donc pas se caractériser comme la résolution d'un problème. Comment un problème considéré comme résolu est-il à l'origine de la discussion mathématique que forme notre corpus ?

### 3. Le caractère spécifique d'une équation algébrique.

Comme nous le verrons plus en détail dans la troisième partie de cet article, le problème des petites oscillations apparait dans l'œuvre de Lagrange en 1766 pour l'étude des oscillations d'un fil chargé d'un « nombre quelconque de masses ». Mais le problème tel qu'il est posé dans la *Mécanique analytique* dépasse le seul cas des cordes vibrantes. Entre 1766 et 1788, l'efficacité du procédé élaboré pour l'intégration des systèmes différentiels à coefficients constants a conduit Lagrange à donner, dans un mémoire adressé depuis Berlin à l'Académie de Paris en 1774, une même mathématisation au problème des petites oscillations d'un fil et au problème des petites oscillations des planètes sur leurs orbites :

> Si les Planètes étaient simplement attirées par le Soleil, et n'agissaient point les unes sur les autres, elles décriraient autour de cet astre, des ellipses variables suivant les lois de Kepler, comme Newton l'a démontré le premier, et une foule d'Auteurs après lui. Mais les observations ont prouvé que le mouvement elliptique des Planètes est sujet à de petites oscillations, et le calcul a démontré que leur attraction mutuelle peut en être la cause. Ces variations sont de deux espèces : les unes périodiques et qui ne dépendent que de la configuration des Planètes entre elles ; celles-ci sont les plus sensibles, et le calcul en a déjà été donné par différents Auteurs ; les autres séculaires et qui paraissent aller toujours en augmentant, ce sont les plus difficiles à déterminer tant par les observations que par la Théorie. Les premières ne dérangent point l'orbite primitive de la Planète ; ce ne sont, pour ainsi dire, que des écarts passagers qu'elle fait dans sa course régulière, et il suffit d'appliquer ces variations au lieu de la Planète calculé par les Tables ordinaires du mouvement elliptique. Il n'en est pas de même des variations séculaires. Ces dernières altèrent les éléments mêmes de l'orbite, c'est-à-dire la position et la dimension de l'ellipse décrite par la planète ; et quoique leur effet soit insensible dans un court espace de temps, il peut néanmoins devenir à la longue très considérable. [Lagrange 1781, p.125].

Malgré les approximations nécessaires à la représentation par un système différentiel à coefficients constants des « équations séculaires des mouvements des nœuds et des inclinaisons des planètes », la mathématisation élaborée par Lagrange est adoptée par Laplace dès 1775 dans le *Mémoire sur les solutions particulières des équations différentielles et sur les inégalités séculaires des planètes* :

> J'ai donné dans un autre Mémoire (Tome VII des Savants étrangers) les expressions des inégalités séculaires des planètes […]. Celles que j'ai donné n'en sont que les différentielles ; je m'étais proposé depuis longtemps de les intégrer ; mais le peu d'utilité de ce calcul pour les besoins de l'Astronomie, joint aux difficultés qu'il présentait, m'avait fait abandonner cette idée , et j'avoue que je ne l'aurais pas reprise, sans la lecture d'un excellent Mémoire *Sur les inégalités séculaires du mouvement des nœuds et de l'inclinaison des orbites des planètes* que M. de Lagrange vient d'envoyer à l'Académie, et qui paraîtra dans un des Volumes suivants. Cet illustre géomètre, au moyen d'une transformation heureuse, réduit le problème à l'intégration d'autant d'équations linéaires du premier ordre qu'il y a d'inconnues ; il donne ensuite une méthode fort ingénieuse pour les intégrer, et pour déterminer les constantes que renferme l'intégrale, quelque soit le nombre de planètes. En employant la même transformation, j'ai tiré les mêmes équations de mes formules ; j'ai de plus cherché si l'on ne pourrait pas déterminer d'une manière analogue les inégalités séculaires de l'excentricité et du mouvement de l'aphélie […] on aura ainsi une théorie complète et rigoureuse de toutes les inégalités séculaires des orbites des planètes. [Laplace 1775, p.354].

La référence de Laplace à Lagrange manifeste la reconnaissance de l'efficacité d'une « méthode fort ingénieuse » ramenant la résolution de problèmes mécaniques à la considération d'une équation algébrique. C'est donc l'efficacité d'une méthode qui justifie la mathématisation du problème par « autant d'équations linéaires du premier ordre qu'il y a d'inconnues » et par là, de la stabilité mécanique par la nature des racines de l'équation algébrique associée au système différentiel. Dans deux publications successives, le *Mémoire sur les inégalités séculaires des planètes et des satellites* [1787] et le *Mémoire sur les variations séculaires des orbites des planètes* [1789], Laplace entreprend de démontrer mathématiquement, en « dehors toute hypothèse » numérique sur la masse des planètes, la stabilité du système solaire c'est-à-dire la nature réelles, négatives et inégales des racines de l'équation algébrique caractéristique de la méthode élaborée par Lagrange. La démonstration de 1789 repose sur un « artifice particulier » s'appuyant sur la présence de « rapports remarquables » [Laplace 1789, p.297], $(1,2)=(2,1)$ et $[1,2]=[2,1]$ dans le système différentiel :

$$0 = (1)\ \frac{d^2\xi}{dt^2} + (1,2)\ \frac{d^2\psi}{dt^2} + (1,3)\ \frac{d^2\varphi}{dt^2} + \&c + [1]\ \xi + [1,2]\ \psi + [1,3]\varphi + \&c,$$
$$0 = (2)\ \frac{d^2\psi}{dt^2} + (1,2)\ \frac{d^2\xi}{dt^2} + (2,3)\ \frac{d^2\varphi}{dt^2} + \&c + [2]\psi + [1,2]\xi + [2,3]\varphi + \&c,$$
$$0 = (3)\ \frac{d^2\varphi}{dt^2} + (1,3)\ \frac{d^2\xi}{dt^2} + (2,3)\ \frac{d^2\psi}{dt^2} + \&c + [3]\varphi + [1,3]\xi + [2,3]\psi + \&c,$$

Comme le manifestent les titres des travaux d'auteurs du corpus comme le mémoire de Cauchy de 1829 intitulé « Sur l'équation à l'aide de laquelle on détermine les inégalités séculaires des planètes », celui de Sylvester de 1852 « Sur une propriété nouvelle de l'équation qui sert à déterminer les inégalités séculaires des planètes » ou encore la formulation abrégée utilisée par Hermite en 1857, « Mémoire sur l'équation à l'aide de laquelle, etc. », les travaux de mécanique de Lagrange et Laplace confèrent un caractère spécifique à une équation algébrique qui se caractérise par la nature de ses racines (réelles, inégales) et les « rapports remarquables » du système dont elle issue.

Nous allons à présent préciser, en évoquant les travaux de Cauchy, la manière dont le caractère spécifique de l'équation des petites oscillations donne au corpus son principal élément de cohérence interne. En 1829, Cauchy formule les problèmes des petites oscillations des systèmes mécaniques, des axes principaux d'un solide en rotation et de la classification des coniques et quadriques comme relevant d'une même question de recherche d'extremum d'une fonction homogène du second degré ([30]). Comme le manifestent les deux citations ci-dessous, la reconnaissance par Cauchy d'une similitude entre ces différents problèmes tient à une analogie formelle des caractéristiques « dignes de remarque » des équations algébriques que leurs résolutions mettent en œuvre ([31]).

> Soit *(1)* $s = f(x\ y,z,...)$ une fonction homogène et du second degré. Soient de plus *(2)* $\varphi(x,y,z,...)$, $\chi(x,y,z,...)$, $\psi(x,y,z,...)$, ... les dérivées partielles de $f(x,y,z,...)$ prises par rapport aux variables $x,y,z,...$ Si l'on assujettit ces variables à l'équation de condition (3) $x^2+y^2+z^2 + ... =1$ les maxima et minima de la fonction $s$ seront déterminées (voir les Leçons sur le Calcul infinitésimal, p.252) par la formule […]
>
> $$(6)\ \frac{1}{2}\varphi(x,y,z,....)=sx\ ,\ \frac{1}{2}\chi(x,y,z,....)= sy\ ,\ \frac{1}{2}\psi(x,y,z,...)= sz$$
>
> Soit maintenant *(7)* $S=0$ l'équation que fournira l'élimination des variables $x,y,z,$ … entre les formules *(6)*. Les *maxima et minima* de la fonction $s=f(x,y,z\ ,..)$ ne pourront être que des racines de l'équation *(7)*. D'ailleurs cette équation sera semblable à celle que l'on rencontre dans la théorie des inégalités séculaires des mouvements des planètes, et dont les racines, toutes réelles, jouissent de propriétés dignes de remarque. Quelques unes de ces propriétés étaient déjà connues : nous allons les rappeler ici, et en indiquer des nouvelles. [Cauchy 1829, p.173].

[30] Si $q$ est une forme quadratique et $\varphi$ sa forme polaire alors $dq_x(y)=2\varphi(x,y)$. On peut déterminer un vecteur propre $e$ de $q$ par une recherche d'extremum sur la sphère unité $S$. Il suffit de démontrer qu'il existe $e \in E$ tel que $\|e\|=1,\ <e,x>=0 \Rightarrow \varphi(e,x)=0$ (< > désigne le produit scalaire, $\|\ \|$ la norme de l'espace euclidien). Soit la fonction $f: E\text{-}\{0\}\rightarrow R,\ x\rightarrow q(x)/\|x\|$ . Comme $S$ est compact, $f/S$ admet un maximum atteint en un point $e$ tel que $q(e)=1$ et $df_{e1}=0$. Mais $dq_x(y)=2\varphi(x,y)$ donc $df_{e1}(x)=0 \Leftrightarrow 2<e,x>-2q(e)\varphi(e,x)=0$. On a bien $<e,x>=0 \Rightarrow \varphi(e,x)=0$.

[31] Hawkins a mis en évidence le rôle joué par Sturm dans la reconnaissance par Cauchy, entre 1826 et 1829, des analogies formelles portées par l'équation des petites oscillations. Les problèmes de conduction de la chaleur abordés par Sturm dans la tradition de Fourier conduisent à l'étude de systèmes d'équations différentielles linéaires semblables à ceux des petites oscillations. Sturm aborde la question de la multiplicité des racines des équations caractéristiques de ces systèmes en appliquant son théorème de dénombrement des racines réelles des équations. Voir à ce sujet, [Hawkins 1975, p.21-22].

Parmi les méthodes employées par les géomètres pour discuter les surfaces représentées par des équations du second degré, l'une des plus simples est celle qui consiste à couper ces surfaces par des droites parallèles. En suivant cette méthode, on peut facilement déterminer la nature des surfaces dont il s'agit, leurs centres, s'il en existe, leurs axes principaux etc. ; et l'on reconnaît, en particulier, que pour fixer la direction de ces axes, il suffit de résoudre une équation du troisième degré. Cette équation qui se représente dans diverses questions de Géométrie ou de Mécanique, et, en particulier, dans la théorie des moments d'inertie, a cela de remarquable que ses trois racines sont toujours réelles. [Cauchy 1828, p.9].

Dans la période qui sépare les travaux de Lagrange de 1766 et les organisations théoriques élaborées par Jordan et Kronecker en 1874, l'équation à l'aide de laquelle on détermine les inégalités séculaires des planètes revêt une identité qui ne se limite pas à un cadre théorique mais l'identifie à un corpus d'auteur*s* que l'on sollicite lorsque cette équation se manifeste et qui s'étoffe en conséquence. Nous proposons à présent d'étudier l'origine du caractère spécifique de cette équation dans la pratique élaborée par Lagrange pour l'intégration des systèmes différentiels linéaires à coefficients constants. Si une formulation contemporaine de cette pratique polynomiale a déjà été présentée par l'expression (*) donnée dans le troisième paragraphe de la partie I et détaillée dans la note n°23, les deux premiers paragraphes de la troisième partie de cet article proposent d'entrer dans la technicité mathématique des textes originaux, une synthèse est proposée dans le troisième paragraphe dans un cadre mathématique plus familier au lecteur contemporain.

# III. La spécificité de la pratique algébrique élaborée par Lagrange pour le « problème des oscillations très petites d'un système quelconque de corps ».

## 1. La portée générale donnée au problème par la *Mécanique Analytique* de 1788.

Le problème « des oscillations très petites d'un système quelconque de corps » ouvre la section V de la *Mécanique Analytique* intitulée « Solutions de différens problêmes de Dynamique ». Relevant d'une méthode « générale » applicable à un système « quelconque » de corps, le problème des petites oscillations illustre la portée des « principes généraux » exposés dans le traité et porte des enjeux historiques. L'architecture de la *Mécanique Analytique* ne peut en effet être dissociée d'une certaine forme d'histoire, faite d'un jeu de postérité, portée par une idée de « révolution » de la science par les modernes et présentée dans l'exposé préliminaire à l'ouvrage :

La dynamique est la science des forces accélératrices ou retardatrices, & des mouvements variés qu'elles peuvent produire. Cette science est due entièrement aux Modernes, & Galilée est celui qui en a jeté les premiers fondements […]. La Mécanique devint une science nouvelle entre les mains de Newton, & ses *Principes Mathématiques* qui parurent pour la première fois en 1687, furent l'époque de cette révolution […]. Enfin l'invention du calcul infinitésimal mit les Géomètres en état de réduire à des équations analytiques les lois du mouvement des corps; & la recherche des forces & des mouvements qui en résultent est devenue depuis le principal objet de leurs travaux. Je me suis proposé ici de leur offrir un nouveau moyen de faciliter cette recherche […]. [Lagrange 1788, p.158-159].

La *Mécanique analytique,* qui se présente comme une œuvre de synthèse du siècle écoulé, a prétention à participer de cette « révolution » par sa « facilité ». En faisant référence à la querelle des cordes vibrantes, Lagrange oppose la « simplicité » des principes exposés dans son traité aux difficultés et confusions suscitées par l'étude des systèmes de corps aux XVII$^{e}$ et XVIII$^{e}$ siècles.

La résolution du « problème des petites oscillations d'un système quelconque de corps », que Daniel Bernouilli pensait trop irrégulier pour les méthodes analytiques et dont Jean le Rond d'Alembert n'était parvenu à traiter que le cas particulier d'un système restreint à deux ou trois masses, illustre à la fois la simplicité et la généralité qui font la nouveauté de la *Mécanique analytique* de 1788 ([32]). La mise en équation du mouvement vibratoire d'une corde sans poids, chargée de deux masses, illustrait déjà, dans le *Traité de dynamique* de 1743, la portée générale du principe des forces vives sur lequel d'Alembert avait unifié la mécanique des corps solides ([33]).

> Un fil *CmM* fixe en *C*, & chargé de deux poids *m, M,* étant infiniment peu éloigné de la verticale *CO*, trouver la durée des oscillations de ce fil [...].
>
> $$(K) - ddx = [\ \frac{px}{l} - \frac{M.P}{m}\ (\ \frac{y}{L} - \frac{x}{l}\ )]\ dt^2$$
>
> $$\& \ -\ ddy = [\ \frac{yP}{L} . \frac{M+m}{m} - \frac{x}{L}.(p+ \frac{M.P}{m}\ )]dt^2\ (N).$$
>
> [d'Alembert 1758, p.139-142].

Sous la simplification *M=m, l=L,* notant *T* le temps pendant lequel la force accélératrice fait parcourir au corps un espace d'une unité, les équations du mouvement sont les suivantes :

$$(P) - \quad ddx = (2x\text{-}y).\ \frac{2dt^2}{T^2}\ ,$$

$$\&\ (Q)- \quad ddy= (2y\text{-}2x)\ \frac{2dt^2}{T^2}\ .$$

d'Alembert recourt à la méthode des coefficients indéterminées pour ramener le système des équations du mouvement à deux équations pouvant s'intégrer séparément. La méthode consiste à introduire un coefficient indéterminé ν permettant de combiner les deux équations en une équation algébrique du second degré dont les deux solutions, $\nu=\frac{\pm 1}{\sqrt{2}}$, donnent les variables indépendantes recherchées, *u=x+vy et u'=x+vy* :

> Je multiplie la seconde par un coefficient indéterminé *v*, & ensuite je les ajoute ensemble, ce qui donne
>
> $$-ddx\text{-}vddy = \frac{2dt^2}{T^2} \times (\overline{2\text{-}2v}\ .x+ \overline{2v\text{-}1}\ .y)\ldots(R).$$
>
> Je fais en sorte *que (2-2v)x+(2v-1)y* soit un multiple de *–x-vy*, ce qui donne
>
> $$2\text{-}2v= \frac{2v-1}{v}\ ;\ \&\ v= \frac{\pm 1}{\sqrt{2}}\ ;$$
>
> donc faisant *x+vy=u*, ou plutôt
>
> $$x+ \frac{y}{\sqrt{2}} =u,\ \&\ x\text{-} \frac{y}{\sqrt{2}} =u',$$
>
> on aura les deux équations
>
> $$-ddu=(2\text{-}\sqrt{2})\ 2u\ \frac{dt^2}{T^2},\quad \&\ \text{-}ddu'=(2+\sqrt{2}).\ \frac{2u'dt^2}{T^2}\ .$$

L'intégration de ces équations linéaires à coefficients constants, basée sur la solution donnée par Euler en 1743 ([34]), est présentée par d'Alembert dans un mémoire adressé à l'académie de Berlin en 1747 [d'Alembert 1750] ([35]) :

---

[32] Au sujet du modèle de théorie physique et du style de mathématisation associés au formalisme de Lagrange, consistant à réduire et ordonner les différents principes connus en mécanique autour d'un unique principe général et d'un formalisme « algébrique », voir [Fraser 1980 et 1985], [Barroso-Filho et Comte 1988], [Dahan Dalmedico 1990 et 1992], [Galleto 1991] et [Panza 1992 et 2003].

[33] Au sujet des principes variationnels chez Bernouilli, d'Alembert et Lagrange, voir [Panza 2003, p.147]. Les citations seront extraites de la deuxième édition de 1758. La notation *x* désigne le déplacement vertical du corps *m*, *y* celui du corps *M*, *l* le déplacement parcouru pendant le premier instant par *m* et *L* par *M*. *P* est la pesanteur du corps *M* et *p* celle de *m*.

[34] Au sujet des travaux d'Euler et de D. Bernouilli sur les petites oscillations d'une corde et à propos de la solution d'Euler de 1743, consulter [Gilain 2003, p.443].

[35] Ci-dessous, un exemple avec des notations contemporaines : si on considère des oscillateurs harmoniques couplés, de même masse *m*, reliés à des points *A* et *B* par des ressorts de même raideur $k_0$ et reliés entre eux par un ressort de raideur *k*, glissant sans frottement sur *(AB)*, les équations du mouvement s'écrivent :

$$md^2x_1 /dt^2 = -\ k_0x_1 - k(x_1\text{-}x_2)$$
$$md^2x_2 /dt^2 = -\ k_0x_2 - k(x_2\text{-}x_1)$$

[…] multipliant la 1[ere] par $du$, on a l'intégrale $\frac{du^2}{\sqrt{AA-uu}} = \frac{dt}{T}\sqrt{4-2\sqrt{2}}$, parceque $t$ croissant, $u$ diminue ; donc

$$u=A\cos\frac{t\sqrt{4-2\sqrt{2}}}{T}, \;\&\; u'=B\cos.\frac{t\sqrt{4+2\sqrt{2}}}{T},$$

intégrales qui sont complètes, […]; de là on tirera les valeurs de $x$ & de $y$, & on déterminera les constantes $A$&$B$ par les valeurs connues & données de $x$ & de $y$ lorsque $t=0$. [*Ibidem*, p.145].

Bien que revendiquant la généralité de sa méthode, d'Alembert ne mettra celle-ci en œuvre que pour les cas de deux masses et de trois masses égales. Comme nous allons le voir, lorsque Lagrange entreprendra de résoudre, en 1766, le problème dans le cas d'un « système quelconque de corps », cette *généralisation* s'accompagnera de l'élaboration d'une pratique *polynomiale* spécifique s'affranchissant de la méthode des coefficients indéterminés.

## 2. Une pratique élaborée par Lagrange en un *jeu* sur les *primes* et les *indices* des systèmes linéaires.

Lagrange aborde pour la première fois le problème des petites oscillations d'une corde dans ses « Solutions de différents problèmes de calcul intégral » [1766], écrites entre 1762 et 1765. La résolution se présente comme une application du principe, exposé dans le même mémoire, de réduction de l'ordre d'une équation différentielle dont on connaît des solutions particulières :

> Les équations proposées seront intégrables algébriquement, si l'on peut trouver, […] autant de valeurs particulières de chacune des quantités *y, y', y'',...* qu'il y a d'unités dans la somme des exposants des plus hautes différences de ces variables. […] Méthode générale pour déterminer le mouvement d'un système quelconque de corps qui agissent les uns sur les autres, en supposant que ces corps ne fassent que des oscillations infiniment petites autour de leurs points d'équilibre.
> Soit $n$ le nombre des corps qui composent le système, et nommons *y', y'', y''',...* les espaces infiniments petits que ces corps décrivent dans leurs oscillations pendant le temps $t$ ; on aura, en négligeant les quantités infiniment petites du second ordre et des ordres plus élevés des équations de cette forme
>
> $$(a)\begin{cases}\frac{d^2y'}{dt^2}+A'y'+B'y''+C'y'''+\ldots+N'y^{(n)}=0\\ \frac{d^2y''}{dt^2}+A''y'+B''y''+C''y'''+\ldots+N''y^{(n)}=0\\ \frac{d^2y'''}{dt^2}+A'''y'+B'''y''+C'''y'''+\ldots+N'''y^{(n)}=0\\ \ldots\\ \frac{d^2y^{(n)}}{dt^2}+A^{(n)}y'+B^{(n)}y''+C^{(n)}y'''+\ldots+N^{(n)}y^{(n)}=0\end{cases}$$
>
> *A', B', C',..., A'', B'', C'',...* étant des constantes données par la nature du problème. [Lagrange 1766, p.519].

- **Les coefficients en miroirs des systèmes obtenus par un emploi de la méthode des coefficients indéterminés.**

Le système différentiel des petites oscillations est d'ordre 2. Lagrange considère un système de $n$ solutions particulières de la forme $\lambda e^{\rho t}dt$, où $\lambda$ et $\rho$ sont des « constantes indéterminées » permettant de réduire l'ordre du système. Les $n$ équations sont ramenées en une seule par

---

Le changement de variable : $u = x_1+x_2$ et $v = x_1-x_2$ donne les deux équations indépendantes :

$$md^2u/dt^2 = -k_0u$$
$$md^2v/dt^2 = -k_0 v-2kv$$

Les solutions $u=a\cos\omega' t$ et $v=a\cos\omega'' t$ peuvent s'interpréter comme deux oscillateurs indépendants, tandis que les solutions du système feront apparaître un phénomène de battement :

$$x_1 = (a/2)(\cos\omega' t + \cos\omega'' t)$$
$$x_2 = (a/2)(\cos\omega' t - \cos\omega'' t)$$

sommation après multiplication par les coefficients indéterminés $\lambda e^{\rho t}d$. L'unique équation obtenue peut alors s'intégrer par partie :

> Pour intégrer ces équations suivant la méthode expliquée ci-dessus, on multipliera la première par $\lambda' e^{\rho t}dt$, la seconde *par* $\lambda'' e^{\rho t}dt$, et ainsi de suite, $\lambda, \lambda', \lambda''$, ... étant, ainsi que $\rho$, des constantes indéterminées ; ensuite on les ajoutera ensemble, et on en prendra l'intégrale en faisant disparaître de dessous le signe $\int$ les différences des variables $y', y'', y'''$, ... ; après quoi on fera les coefficients des quantités $\int y' e^{\rho t}dt$, $\int y'' e^{\rho t}dt$, $\int y''' e^{\rho t}dt$,... égaux à zéro ; de cette manière on aura d'abord l'équation intégrale
>
> $$(b)\ [\lambda'(\frac{dy'}{dt}-\rho y') + \lambda''(\frac{dy''}{dt}-\rho y'')+ \lambda'''(\frac{dy'''}{dt}-\rho y''')+...+\lambda^{(n)}(\frac{dy^{(n)}}{dt}-\rho y^{(n)})]e^{\rho t}=const.$$ [*Ibidem*].

Les termes $\frac{dy'}{dt}-\rho y'$ expriment la qualité de solutions particulières des expressions $\lambda' e^{\rho t}dt$ réduisant l'ordre du système différentiel d'une unité. La méthode des coefficients indéterminés donne les conditions sur les $\lambda$ et $\rho$ pour l'obtention de l'équation *(b)* : « on fera les coefficients des quantités $\int y' e^{\rho t}dt$, $\int y'' e^{\rho t}dt$, $\int y''' e^{\rho t}dt$,... égaux à zéro ». Ces conditions sont exprimées par un système d'équations linéaires ([36]) :

$$(c)\begin{cases} \rho^2\lambda'+A'\lambda'+A''\lambda''+A'''\lambda'''+...+A^{(n)}\lambda^{(n)}=0, \\ \rho^2\lambda''+B'\lambda'+B''\lambda''+B'''\lambda'''+...+B^{(n)}\lambda^{(n)}=0, \\ \rho^2\lambda'''+C'\lambda'+C''\lambda''+C'''\lambda'''+...+C^{(n)}\lambda^{(n)}=0, \\ ... \\ \rho^2\lambda^{(n)}+N'\lambda'+N''\lambda''+N'''\lambda'''+...+N^{(n)}\lambda^{(n)}=0, \end{cases}$$

> lesquelles [équations] serviront à déterminer les quantités $\rho, \lambda', \lambda'', \lambda'''$,... . [*Ibid.*, p.520]

L'« élimination des équations » *(c)* donne une équation de degré $n$ en $\rho$ notée $P=0$ et permettant de déterminer les $n$ valeurs $\rho_1, \rho_2,...,\rho_n$ des indéterminés $\rho$. Les solutions particulières $\lambda e^{\rho t}dt$ obtenues, la résolution de l'équation intégrale *(b)* permettraient alors de trouver les solutions aux conditions initiales, Lagrange élabore cependant une méthode permettant d'obtenir l'expression générale, c'est-à-dire polynomiale, de ces solutions par la seule donnée de l'équation algébrique $P=0$. Comme nous allons le voir, la méthode développée par Lagrange pour la résolution du système aux conditions initiales est propre au problème des petites oscillations et met en œuvre un jeu sur les primes et les indices à partir d'une propriété des systèmes obtenus par la méthode des coefficients indéterminés : les coefficients du système *(c)* sont en *miroirs* de ceux du système *(a)*, par exemple le coefficient $A''$, deuxième coefficient de la première équation de *(a)*, est identique au premier coefficient de la deuxième équation de *(c)*.

- **Expression des solutions en fonction des conditions initiales.**

Lagrange note $\rho_1^2, \rho_2^2, \rho_3^2$,... les racines, supposées réelles positives et distinctes, de l'équation $P=0$. Si les conditions initiales pour $t=0$ sont $y'=Y'$, $y''= Y''$,..., $dy'/dt=V'$, $dy''/dt=V''$, ..., l'équation *(b)* s'écrit :

$$\lambda'\frac{dy'}{dt} + \lambda''\frac{dy''}{dt} + \lambda'''\frac{dy'''}{dt} + ... + \lambda^{(n)}\frac{dy^{(n)}}{dt} \quad -\rho\ [\lambda' y' - \lambda'' y'' - ... - \lambda^{(n)}y^{(n)}]$$
$$= [\lambda'V' + \lambda''V'' + \lambda'''V''' + .... + \lambda^{(n)}V^{(n)} - \rho\,(\lambda'Y' + \lambda''Y'' + ... + \lambda^{(n)}Y^{(n)})]\ e^{-\rho t}$$

Cette équation met en relation les conditions initiales ($Y$) et ($V$) et le système des coefficients ($\lambda$) des solutions particulières $\lambda e^{\rho t}dt$. Les racines de $P=0$ s'écrivant $\rho_i^2$ et l'équation intégrale *(b)* étant

---

[36] Par exemple la méthode des coefficients indéterminés sur le système à deux inconnues $\begin{cases}\frac{d^2y'}{dt^2}+A'y'+B'y''=0 \\ \frac{d^2y''}{dt^2}+A''y'+B''y''=0\end{cases}$

donne l'équation $\lambda' e^{\rho t}(\frac{d^2y'}{dt^2}+A'y'+B'y'')+\lambda'' e^{\rho t}(\frac{d^2y''}{dt^2}+A''y'+B''y'')=0$.

L'équation intégrale est obtenue par intégration par partie :

$$\int\frac{d^2y}{dt^2}e^{\rho t}dt = \int\rho\frac{dy}{dt}e^{\rho t}dt - [\frac{dy}{dt}e^{\rho t}] = \int\rho^2 ye^{\rho t}dt + [\rho ye^{\rho t}] - [\frac{dy}{dt}e^{\rho t}]$$

à la condition que les termes en $\int y'' e^{\rho t}dt$ s'annulent, c'est-à-dire que $\rho^2\lambda''+A''\lambda'+A''\lambda''=0$.

vérifiée pour $\rho$ et $-\rho$, sommer les équations intégrales obtenues pour ces valeurs opposées permet d'écrire la relation :

$$\lambda' y' + \lambda'' y'' + \ldots + \lambda^{(n)} y^{(n)} = [\lambda' Y' + \lambda'' Y'' + \ldots + \lambda^{(n)} Y^{(n)}]\ \frac{e^{\rho t}+e^{-\rho t}}{2}$$
$$+ [\lambda' V' + \lambda'' V'' + \lambda''' V''' + \ldots + \lambda^{(n)} V^{(n)}]\ \frac{e^{\rho t}-e^{-\rho t}}{2\rho}$$

Lagrange note $\theta$ le membre de droite qui exprime les conditions initiales (*Y*) et (*V*) dans le système de variables *(λ),*

$$\theta = [\lambda' Y' + \lambda'' Y'' + \ldots + \lambda^{(n)} Y^{(n)}]\ \frac{e^{\rho t}+e^{-\rho t}}{2} + [\lambda' V' + \lambda'' V'' + \lambda''' V''' + \ldots + \lambda^{(n)} V^{(n)}]\ \frac{e^{\rho t}-e^{-\rho t}}{2\rho}$$

Le système *(c)* associe aux *n* valeurs de l'indéterminée $\rho^2$, c'est-à-dire aux *n* racines $\rho_1^2, \rho_2^2, \rho_3^2, \ldots$, *n* systèmes de valeurs $(\lambda'_1, \lambda''_1, \ldots, \lambda^{(n)}_1)$, $(\lambda'_2, \lambda''_2, \ldots, \lambda^{(n)}_2)$, ...., $(\lambda'_n, \lambda''_n, \ldots, \lambda^{(n)}_n)$ des coefficients des solutions particulières et *n* valeurs correspondantes des conditions initiales $\theta_1, \theta_2, \ldots$ ([37]). Déterminer les solutions aux conditions initiales revient alors à résoudre le système suivant afin d'exprimer les *y* en fonction des $(\lambda_i)$ et $(\theta_i)$:

$$\lambda'_1 y' + \lambda''_1 y'' + \ldots + \lambda^{(n)}_1 y^{(n)} = \theta_1,$$
$$\lambda'_2 y' + \lambda''_2 y'' + \ldots + \lambda^{(n)}_2 y^{(n)} = \theta_2,$$
$$\ldots\ldots$$
$$\lambda'_n y' + \lambda''_n y'' + \ldots + \lambda^{(n)}_n y^{(n)} = \theta_n,$$

- **Une première application de la méthode des coefficients indéterminés.**

L'introduction d'un premier jeu d'indéterminées *μ', μ'', ...* donne une première expression des solutions $y^{(s)}$ ([38]) :

> Pour en venir à bout, je multiplie la première de ces équations par *μ',* la seconde par *μ''*, la troisième par *μ'''*, et ainsi de suite, *μ', μ'', μ''',...* étant des coefficients indéterminés, puis je les ajoute ensemble, ce qui me donne [...] la valeur d'une *y* quelconque, comme $y^{(s)}$, en égalant à zéro chacun des coefficients des autres *y* ; ainsi l'on aura
>
> $$(e)\quad y^{(s)} = \frac{\mu'\theta_1 + \mu''\theta_2 + \mu'''\theta_3 + \ldots + \mu^{(n)}\theta_n}{\mu'\lambda_1^{s} + \mu''\lambda_2^{(s)} + \mu'''\lambda_3^{(s)} + \ldots + \mu^{(n)}\lambda_n^{(s)}},$$
>
> et ensuite ces équations de condition :
>
> $$(f)\ \mu'\lambda'_1 + \mu''\lambda'_2 + \ldots + \mu^{(n)}\lambda'_n = 0,$$
> $$\mu'\lambda''_1 + \mu''\lambda''_2 + \ldots + \mu^{(n)}\lambda'_n = 0,$$
> $$\ldots.$$
> $$\mu'\lambda_1^{(n)} + \mu''\lambda^{(n)}{}_2 + \ldots + \mu^{(n)}\lambda'^{(n)}_n = 0,$$
>
> à l'exception seulement de celle qui répondrait à l'exposant *s*. [*Ibid.,* p.523].

L'expression donnée à la solution $y^{(s)}$ par la formule (*e*) ne satisfait pas Lagrange car les *primes* et les *indices* se trouvent mélangés et les coefficients présents au dénominateur n'ont donc pas de signification par rapport au problème considéré. Comme l'illustre le tableau suivant qui représente l'affectation des $\lambda_i^{(j)}$ aux racines $\rho_i$, le dénominateur de la formule (*e*) mêle $\lambda_1^s$, le $s^e$ coefficient associé à $\rho_1$ à $\lambda_2^s$ le $s^e$ coefficient associé à $\rho_2$ etc. Les *indices* de $\lambda^{(s)}$ sont à l'endroit où l'on attendrait des *primes* de $\lambda_s$ ([39]).

---

[37] En des termes qui nous sont contemporains : les coefficients θ expriment les vecteurs conditions initiales *Y* dans la base de vecteurs propres *(λ)* de la matrice associée au système (a) : $\Theta = \Lambda Y$.

[38] La méthode des coefficients indéterminés s'interprèterait aujourd'hui comme une recherche de la solution $y^{(s)}$ par des combinaisons sur les lignes annulant les coefficients des autres inconnues. Le système *(f)* donne les conditions sur les coefficients indéterminées *(μ),* c'est-à-dire les combinaisons à effectuer sur les lignes afin d'isoler $y^{(s)}$. La solution $(y^{(s)})$ donnée par l'équation *(e)* est analogue à celle exprimée par la règle de Cramer comme quotient de deux déterminants. L'emploi de déterminants pour exprimer la méthode de Lagrange sera d'abord effectué par [Laplace 1775], le rôle joué par les déterminant sera discuté dans la suite de ce travail lorsque sera abordée la méthode de [Cauchy 1829].

[39] $(\lambda_1', \ldots, \lambda_1^{(n)})$ est un système de cordonnées du vecteur propre de $\rho_1^2$. $\lambda_1^{(i)}$ est la $i^e$ coordonnée du vecteur propre. La méthode de l'élimination appliquée au système *(c)* conduit à des relations que l'on interprèterait aujourd'hui comme des conditions d'orthogonalités duales. Si les $(\lambda_i)$ sont les vecteurs propres de *A* alors les $(v_m)$ sont construit comme appartenant à l'orthogonal de l'hyperplan $(\lambda_1, \lambda_2, \ldots, \lambda_{m-1}, \lambda_{m+1}, \ldots, \lambda_n)$ et sont dont les vecteurs propres de la matrice ${}^tA$ .

| | | | |
|---|---|---|---|
| $\rho_1^2$ | $\rho_2^2$ | ... | $\rho_n^2$ |
| $(\lambda'_1, \lambda''_1, ..., \lambda^{(n)}_1)$ | $(\lambda'_2, \lambda''_2, ..., \lambda^{(n)}_2)$ | ... | $(\lambda'_n, \lambda''_n, ..., \lambda^{(n)}_n)$ |
| $\theta_1$ | $\theta_2$ | ... | $\theta_n$ |

La méthode des coefficients indéterminés inverse les primes et les indices, Lagrange recourt donc à un second emploi de cette méthode sur le système *(f)* afin de rétablir l'ordre désiré entre primes et indices.

- **Une seconde application de la méthode des coefficients indéterminés.**

Le système comprend $n$ équations, dont $n-1$ de seconds membres nuls (les équations *f)*, et une de second membre non nul *(e)*.

Supposons que l'on ait en général

$$\mu'\lambda'_1 + \mu''\lambda'_2 + .... + \mu^{(n)}\lambda'_n = \Delta'$$
$$\mu'\lambda''_1 + \mu''\lambda''_2 + .... + \mu^{(n)}\lambda'_n = \Delta''$$
$$....$$
$$\mu'\lambda_1^{(n)} + \mu''\lambda^{(n)}{}_2 + .... + \mu^{(n)}\lambda'_n{}^{(n)} = \Delta^{(n)},$$

et qu'il faille trouver la valeur d'une $\mu$ quelconque comme $\mu^{(m)}$. On multipliera ces équations par des coefficients indéterminés $v', v'', v''',...,v^{(n)}$, et, après les avoir ajoutées ensemble, on fera les coefficients des quantités $\mu', \mu'', \mu''',$... chacun égal à zéro, excepté celui de la quantité $\mu^{(m)}$ ; de cette manière on aura

$$(g)\quad \mu^{(m)} = \frac{v'\Delta' + v''\Delta'' + v'''\Delta''' + ... + v^{(n)}\Delta^{(n)}}{v'\lambda'_m + v''\lambda''_m + v'''\lambda'''_m + ... + v^{(n)}\lambda_m^{(n)}},$$

et la détermination des quantités $v', v'', v''',$... dépendra de cette condition que

$$(h)\ v'_m\lambda' + v'_m{}'\lambda'' + ... + v(n)_m\,\lambda(n) = 0$$

lorsque $\rho = \rho_1, \rho_2, \rho_3, ...,\rho_n$ excepté $\rho_m$. [*Ibid.*, p.525].

La double inversion du rôle des primes et des indices permet un retour à la configuration d'origine. L'expression *(g)* obtenue par l'introduction des nouvelles indéterminées *(v)* exprime en effet les solutions $\mu^{(m)}$ à l'aide des coefficients $(\lambda'_m, \lambda''_m, ...., \lambda^{(n)}_m)$ associés à une même racine $\rho_m$. La formule *(g)* permet alors de déterminer $\mu^m$ puis, en remontant les calculs, la solution $y^{(s)}$ cherchée.

- **L'identité des équations algébriques associées aux systèmes obtenus par la méthode des coefficients indéterminés.**

Le double emploi de la méthode des coefficients indéterminés que nous avons vu à l'œuvre chez Lagrange donne un procédé d'obtention des solutions aux conditions initiales $y^{(s)}$ sans pour autant permettre une expression directe de ces solutions. Afin d'obtenir une telle expression, Lagrange va tirer partie du jeu sur les primes et les indices qu'implique chaque introduction de coefficients indéterminés. Si l'on applique les coefficients $v', v''$, ... au système *(c)*, on obtient une équation *(i)* et des conditions sur les *(v)* données par un nouveau système (*k)*:

$$(c)\begin{cases} \rho^2\lambda' + A'\lambda' + A''\lambda'' + A'''\lambda''' + ... + A^{(n)}\lambda^{(n)} = 0, \\ \rho^2\lambda'' + B'\lambda' + B''\lambda'' + B'''\lambda''' + ... + B^{(n)}\lambda^{(n)} = 0, \\ \rho^2\lambda''' + C'\lambda' + C''\lambda'' + C'''\lambda''' + ... + C^{(n)}\lambda^{(n)} = 0, \\ \quad ... \\ \rho^2\lambda^{(n)} + N'\lambda' + N''\lambda'' + N'''\lambda''' + ... + N^{(n)}\lambda^{(n)} = 0, \end{cases} \leftrightarrow \quad (i)\begin{cases} [\rho^2\upsilon' + A'\upsilon' + B'\upsilon'' + C'\upsilon''' + ... + N^{(n)}\upsilon^{(n)}]\lambda' \\ +[\rho^2\lambda'' + A''\upsilon' + B''\upsilon'' + C'''\upsilon''' + ... + N''\upsilon^{(n)}]\lambda'' \\ +[\rho^2\lambda''' + A'''\upsilon' + B'''\upsilon'' + C'''\upsilon''' + ... + N'''\upsilon^{(n)}]\lambda''' \\ \quad ... \\ +[\rho^2\upsilon^{(n)} + A^{(n)}\upsilon' + B^{(n)}\upsilon'' + C^{(n)}\upsilon''' + ... + N^{(n)}\upsilon^{(n)}]\lambda^{(n)} = 0, \end{cases}$$

$$(k)\begin{cases} \rho^2\upsilon' + A'\upsilon' + B'\upsilon'' + C'\upsilon''' + ... + N^{(n)}\upsilon^{(n)} = 0, \\ \rho^2\lambda'' + A''\upsilon' + B''\upsilon'' + C'''\upsilon''' + ... + N''\upsilon^{(n)} = 0, \\ \rho^2\lambda''' + A'''\upsilon' + B'''\upsilon'' + C'''\upsilon''' + ... + N'''\upsilon^{(n)} = 0, \\ \quad ... \\ \rho^2\upsilon^{(n)} + A^{(n)}\upsilon' + B^{(n)}\upsilon'' + C^{(n)}\upsilon''' + ... + N^{(n)}\upsilon^{(n)} = 0, \end{cases}$$

Les deux jeux de coefficients *(λ)* et *(v)* sont respectivement déterminés par les systèmes *(c)* et *(k)* entretenant une relation que Lagrange précise de la manière suivante : les systèmes *(c)* et *(k)* sont différents mais l'équation $P=0$ qui permet de trouver $\rho$ par élimination est identique pour les deux systèmes ([40]):

> Et il est bon de remarquer qu'en éliminant de ces équations les quantités *v', v'', v''',...*, on aura une équation finale en $\rho^2$ qui sera nécessairement la même que celle qui résulte des équations (c) par l'évanouissement des quantités $\lambda', \lambda'', \lambda''',...$; ce qui peut se voir aisément à *priori*. [*Ibid.*].

L'équation algébrique $P=0$ permet par conséquent de déterminer les coefficients *v* directement à partir du système *(a)* sans nécessiter le double emploi de la méthode des coefficients indéterminés.

- **L'expression directe des solutions aux conditions initiales à l'aide de l'équation des petites oscillations.**

La condition *(i)* sur les indéterminées ν s'identifie, à un coefficient multiplicatif prêt, à l'équation algébrique $P=0$ :

> Faisons maintenant $\rho=\rho_m$, nous aurons
>
> $$\begin{cases} A'\upsilon'+B'\upsilon''+C'\upsilon'''+...+N^{(n)}\upsilon^{(n)} = -\rho_m^2\upsilon', \\ A''\upsilon'+B''\upsilon''+C''\upsilon'''+...+N''\upsilon^{(n)} = -\rho_m^2\upsilon'', \\ A'''\upsilon'+B'''\upsilon''+C'''\upsilon'''+...+N'''\upsilon^{(n)} = -\rho_m^2\upsilon''', \\ \quad ... \\ A^{(n)}\upsilon'+B^{(n)}\upsilon''+C^{(n)}\upsilon'''+...+N^{(n)}\upsilon^{(n)} = -\rho_m^2\upsilon^{(n)}, \end{cases}$$
>
> et l'équation *(i)* deviendra
>
> $$(\rho^2-\rho^2_m)\,[v'\lambda'+v''\lambda''+v''\lambda'''+....+v^{(n)}\lambda^{(n)}]=0$$
>
> laquelle devra être vraie pour toutes les valeurs de $\rho$ qui satisfont aux équations (e), d'où celle-ci est tirée, on aura en général
>
> $$[v'\lambda'+v''\lambda''+v''\lambda'''+....+v^{(n)}\lambda^{(n)}]=0,$$
>
> lorsque $\rho=\rho_1, \rho_2, \rho_3,...,\rho_n$ excepté $\rho_m$, auquel cas l'équation se vérifie d'elle-même, à cause du facteur $\rho^2-\rho^2_m$.[*Ibid*, p.527].

Les *n-1* conditions *(h)* qui définissent les *(v)* sont ainsi vérifiées pour $\rho=\rho_1,..., \rho_n$, excepté $\rho_m$ :

$$(h)\ v'_m\lambda' + v'_m{}'\lambda'' + ... + v(n)_m\,\lambda(n) = 0$$

Ces conditions s'expriment par une équation polynomiale de degré *n-1* dont les racines coïncident avec *n-1* racines de $P=0$ et qui s'identifie donc, à un coefficient multiplicatif près, à une factorisation de $P$ par les coefficients $(\rho^2-\rho_m^2)$ :

> $$(\rho^2-\rho^2_m)\,[v'_m\lambda'+v''_m\lambda''+v'''_m\lambda'''+....+v^{(n}{}_m{}^{)}\lambda^{(n)}]=\xi\, P=(I-\frac{\rho^2}{\rho_1^2})(I-\frac{\rho^2}{\rho_2^2})....(I-\frac{\rho^2}{\rho_n^2});$$
>
> $$[...]\ \upsilon_m'\lambda'+\upsilon_m''\lambda''+\upsilon_m'''\lambda'''+...+\upsilon_m^{(n)'}\lambda^{(n)}=\chi(I-\frac{\rho^2}{\rho_1^2})(I-\frac{\rho^2}{\rho_2^2})....(I-\frac{\rho^2}{\rho_n^2})$$
>
> en prenant tous les facteurs hormis $(1-\rho^2/\rho_m^2)$. [*Ibid*].

L'expression $Q_m$ $v'_m\lambda'+...+v_m^{(n)}\lambda^{(n)}$ qui permet, via la condition *(h)* de déterminer les $\mu$ s'obtient par conséquent en factorisant l'équation algébrique $P=0$ par $\rho$-$\rho_m$ ([41]). Il est alors possible d'exprimer la solution aux conditions initiales $y^{(s)}$ directement à l'aide de l'équation des petites oscillations ([42]):

$$y^{(s)} = \frac{v_1^{(s)}}{Q_1}\,\theta_1 + \frac{v_2^{(s)}}{Q_2}\,\theta_2 + ... + \frac{v_n^{(s)}}{Q_n}\,\theta_n$$

Dans la pratique, les facteurs $Q_m$ sont déterminés par de simples différentiations de $P$ :

---

[40] En des termes qui nous sont contemporains, les déterminants (ou les polynômes caractéristiques) de deux matrices transposés sont égaux. Cette propriété n'est pas vraie des vecteurs propres de deux matrices transposées : les systèmes *(λ)* et *(v)* ne coïncident pas même si les $\rho$ coïncident.

[41] D'un point de vue contemporain, la nature des coefficients indéterminées $\lambda$ est changeante, ils désignent le plus souvent des polynômes en $\rho$, et permettent donc d'exprimer le polynôme caractéristique (ils correspondent alors à une écriture générale des vecteurs propres du système). La signification change lorsqu'un indice est attribué à $\lambda$ et $\lambda_m$ correspond alors à un nombre (le coefficient du vecteur propre associé à $\rho_m$).

[42] Cette expression donne une expression polynomiale générale des vecteurs propres comme des quotients des mineurs et de facteurs du polynôme caractéristique analogue à l'expression donnée par les colonnes non nulles de la matrice des cofacteurs (voir à ce sujet la formule * dans I.3. et la note *n*°23). Nous reviendrons sur cette expression lors de la description de son héritage chez Cauchy réalisée dans la quatrième partie de cet article.

Prenons les différences de part et d'autre, en faisant varier ρ, et supposons ensuite $\rho=\rho_m$, […]

$$\frac{\chi dP}{d\rho} = - \frac{2}{\rho_m}[v'_m\lambda_m'+v''_m\lambda_m''+...v^{(n)}{}_m\lambda_m{}^{(n)} = \frac{2Q_m}{\rho_m}$$

donc on aura en général,

$$Q = - \frac{1}{2}\chi\rho\,\frac{dP}{d\rho},$$

Ce qui pourra servir à abréger le calcul de la valeur de $Q$ dans plusieurs occasions. [*Ibid*].

- **Application de la méthode à l'exemple des oscillations d'un fil.**

Lagrange illustre sa méthode en développant l'exemple des « oscillations d'un fil fixe par une de ses extrémités, et chargé d'un nombre quelconque de poids » :

> Soit $n$ le nombre de poids, que nous supposerons, pour plus de simplicité, égaux entre eux et également éloignés les uns des autres ; imaginons que le fil ne fasse que des oscillations infiniments petites et dans le même plan : et soient nommées $y', y'', y''',....., y^{(n)}$ les distances des corps à la verticale, à commencer par le plus bas, et a la distance d'un corps à l'autre : on aura […]

$$\frac{d^2y'}{dt^2}+\frac{y'-y''}{a}=0$$

$$\frac{d^2y''}{dt^2}+\frac{-y'+y''-2y'''}{a}=0$$

$$\frac{d^2y'''}{dt^2}+\frac{-2y''+5y'''-3y^{iv}}{a}=0$$

$$\frac{d^2y^{iv}}{dt^2}+\frac{y^{iv}-y^{v}}{a}-3\,\frac{y'''-2y^{iv}+y^{v}}{a}=0$$

……………………………………………..

$$\frac{d^2y^{iv}}{dt^2}+\frac{-(n-1)y^{(n-1)}+(2n-1)y^{(n)}}{a}$$

[*Ibid.*, p.535]

La résolution nécessite la détermination des indéterminés $v$, de l'équation $P=0$, de ses racines $\rho$ et de ses facteurs $Q_m$. Dans ce cas particulier, Lagrange remarque que les *(λ)* et les *(v)* coïncident car les systèmes *(c)* et *(k)* sont identiques : « à l'égard des quantités ν, on les trouvera de la même manière » que les quantités λ. Les coefficients du système *(c),* disposés en miroirs, ne sont pas affectés par le *jeu sur les primes et les indices.* La seule détermination de l'équation $P=0$ permet donc la résolution générale du système :

$$\lambda'' = (1+a\rho^2)\lambda''$$

$$\lambda'''=\frac{-\lambda'+(3+a\rho^2)\lambda''}{2}=(1+2a\rho^2+\frac{a^2\rho^4}{2})\lambda'$$

………………………………………………………

$$\lambda^{iv}= (1+3a\rho^2+\frac{3a^2\rho^4}{2}+\frac{a^3\rho^6}{2.3})\lambda',$$

$$\lambda^{v} = (1+4a\rho^2+\frac{6a^2\rho^4}{2}+\frac{4a^3\rho^6}{2.3}+\frac{a^4\rho^8}{2.3.4})\lambda',$$

et ainsi de suite ; de sorte qu'on aura en général

$$\lambda^{(m)}=[1+(m-1)a\rho^2+\frac{(m-1)(m-2)}{4}a^2\rho^4+\frac{(m-1)(m-2)(m-3)}{4.9}a^3\rho^6+...]\lambda'.$$

Or il est visible que, pour satisfaire à la dernière équation

$$\rho^2\lambda^{(n)}+\frac{-(n-1)\lambda^{(n-1)}+(2n-1)\lambda^{(n)}}{a}=0$$

il faut supposer $\lambda^{(n+1)}=0$,ce qui donne

$$1+na\rho^2+{}^2+\frac{n(n-1)a^2}{4}\rho^4+\frac{n(n-1)(n-2)a^3}{4.9}\rho^6+....=0$$

équation d'où 'on tirera n valeurs de ρ² […].

$$P=1+na\rho^2+\frac{n(n-1)a^2}{4}\rho^4+\frac{n(n-1)(n-2)a^3}{4.9}\rho^6+....$$

d'où l'on tire [...] $Q = -\frac{1}{2}\chi\rho\ \frac{dP}{d\rho}$. [...] Faisons donc ces substitutions dans la dernière formule du n°30 [les expressions des solutions $y^{(s)}$], on aura l'expression générale des quantités $y$, et le problème sera résolu. [*Ibid.*, p.536]

## 3. Le caractère spécifique de la pratique de Lagrange et les rapports remarquables des systèmes mécaniques.

Nous avons vu Lagrange mettre en œuvre un double emploi de la méthode des coefficients indéterminés pour finalement s'affranchir de cette méthode en exprimant les solutions aux conditions initiales directement à partir de l'équation des petites oscillations. Le processus algébrique original élaboré par Lagrange exploite le jeu sur les primes et les indices impliqué par chaque introduction d'un jeu de coefficients indéterminés. Il permet de constater l'identité des équations associées à deux systèmes dont les coefficients sont en miroirs :

$$(c)\begin{cases} \rho^2\lambda'+A'\lambda'+A''\lambda''+A'''\lambda'''+\ldots+A^{(n)}\lambda^{(n)}=0, \\ \rho^2\lambda''+B'\lambda'+B''\lambda''+B'''\lambda'''+\ldots+B^{(n)}\lambda^{(n)}=0, \\ \rho^2\lambda'''+C'\lambda'+C''\lambda''+C'''\lambda'''+\ldots+C^{(n)}\lambda^{(n)}=0, \\ \ldots \\ \rho^2\lambda^{(n)}+N'\lambda'+N''\lambda''+N'''\lambda'''+\ldots+N^{(n)}\lambda^{(n)}=0, \end{cases} \leftrightarrow \quad (i)\begin{cases} [\rho^2\upsilon'+A'\upsilon'+B'\upsilon''+C'\upsilon'''+\ldots+N^{(n)}\upsilon^{(n)}]\lambda' \\ +[\rho^2\lambda''+A''\upsilon'+B''\upsilon''+C'''\upsilon'''+\ldots+N''\upsilon^{(n)}]\lambda'' \\ +[\rho^2\lambda'''+A'''\upsilon'+B'''\upsilon''+C'''\upsilon'''+\ldots+N'''\upsilon^{(n)}]\lambda''' \\ \ldots \\ +[\rho^2\upsilon^{(n)}+A^{(n)}\upsilon'+B^{(n)}\upsilon''+C^{(n)}\upsilon'''+\ldots+N^{(n)}\upsilon^{(n)}]\lambda^{(n)}=0, \end{cases}$$

Dans le cas particulier des « oscillations d'un fil fixe par une de ses extrémités », Lagrange montre que le jeu sur les primes et les indices des coefficients n'affecte pas le système obtenu. Le procédé élaboré pour résoudre le problème des petites oscillations met alors en évidence une propriété des systèmes de la mécanique, la disposition en miroirs des coefficients que Laplace désignera plus tard comme les « rapports remarquables » des systèmes linéaires et dont Lagrange va démontrer qu'elle découle des principes de la dynamique et de la statique. Dans la citation ci-dessous, le terme *(i,j)* provient du développement de Taylor à l'ordre 2 de l'énergie cinétique et *[i,j]* de celui de l'énergie potentielle.

**Notation de [Lagrange 1766]**

$$(a)\begin{cases} \frac{d^2y'}{dt^2}+A'y'+B'y''+C'y'''+\ldots+N'y^{(n)}=0 \\ \frac{d^2y''}{dt^2}+A''y'+B''y''+C''y'''+\ldots+N''y^{(n)}=0 \\ \frac{d^2y'''}{dt^2}+A'''y'+B'''y''+C'''y'''+\ldots+N'''y^{(n)}=0 \\ \ldots \\ \frac{d^2y^{(n)}}{dt^2}+A^{(n)}y'+B^{(n)}y''+C^{(n)}y'''+\ldots+N^{(n)}y^{(n)}=0 \end{cases}$$

**Notation de [Lagrange 1774]**

$$0=(I)\ \frac{d^2\xi}{dt^2}+(I,2)\ \frac{d^2\psi}{dt^2}+(I,3)\ \frac{d^2\varphi}{dt^2}+\&c+[I]\ \xi+[I,2]\ \psi+[I,3]\varphi+\&c,$$

$$0=(2)\ \frac{d^2\psi}{dt^2}+(I,2)\ \frac{d^2\xi}{dt^2}+(2,3)\ \frac{d^2\varphi}{dt^2}+\&c+[2]\psi+[I,2]\xi+[2,3]\varphi+\&c,$$

$$0=(3)\ \frac{d^2\varphi}{dt^2}+(I,3)\ \frac{d^2\xi}{dt^2}+(2,3)\ \frac{d^2\psi}{dt^2}+\&c+[3]\varphi+[I,3]\xi+[2,3]\psi+\&c,$$

*&c.*

$$0=(n)\ \frac{d^2\varphi}{dt^2}+(I,n)\ \frac{d^2\xi}{dt^2}+(2,n)\ \frac{d^2\psi}{dt^2}+\&c+[n]\varphi+[I,n]\xi+[2,n]\psi+\&c,$$

La pratique élaborée par Lagrange pour la résolution du problème des petites oscillations est spécifique par la manière dont elle s'affranchit du recours à la méthode des coefficients indéterminés en exprimant les solutions des systèmes linéaires par les décompositions polynomiales d'une équation particulière, l'équation des petites oscillations. La disposition en miroirs des coefficients des systèmes mécaniques s'avère indissociable du caractère spécifique de l'équation des petites oscillations dès l'origine. Nous allons à présent questionner l'héritage de la pratique élaborée par Lagrange au sein de la discussion sur l'équation des petites oscillations, des travaux de d'Alembert à ceux de Weierstrass (1858).

# IV. La discussion sur les racines de l'équation des petites oscillations (1766-1858).

## 1. La discussion dans un cadre mécanique chez d'Alembert, Lagrange et Laplace.

Comme nous l'avons vu dans la partie précédente, d'Alembert avait abordé, dans son *Traité de dynamique* de 1743, le problème des petites oscillations de deux ou trois corps par une mise en œuvre de la méthode des coefficients indéterminés qui permettait d'associer au système d'équations différentielles une équation polynomiale. La pratique de d'Alembert ne se réduit cependant pas à une application de la méthode des indéterminés et s'avère sous-tendue par une représentation mécanique du problème. Le principe directeur qui guide l'emploi de la méthode des coefficients indéterminés et ramène le problème à la résolution d'une équation algébrique de degré 2, est en effet basé sur une observation mécanique attribuée à Daniel Bernoulli:

> Je me contenterai de dire que l'on remarque aisément dans les valeurs de *x* & de *y* trouvées ci-dessus, la double oscillation que M. Bernouilli a observée dans le mouvement du pendule dont il s'agit ; chacune de ces oscillations est représentée par chacun des deux termes de la valeur de *x* & de celle de *y*. En effet, l'équation du mouvement d'un pendule simple de longueur $\lambda$, est $-ddz=\frac{2azdt^2}{\lambda T^2}$, ou $z=K\times cof.\frac{1}{T}\sqrt{\frac{2a}{\lambda}}$ ; d'où il est facile de voir que les mouvemens des corps *M*, *m* sont composées de deux mouvemens, synchrones chacun à celui d'un pendule simple. [d'Alembert 1758, p.152].

Les mouvements d'un fil lesté de deux masses *m* et *M* peuvent se « représenter » mécaniquement comme « composés » des deux oscillations indépendantes d'un « pendule simple ». De la même manière que « chacune des oscillations est représentée » par *x* et *y,* la méthode des coefficients indéterminés représente algébriquement l'observation de Bernouilli. Le coefficient indéterminé *v* bénéficie d'une double représentation : représentation algébrique d'une part comme solution d'une équation du second degré, représentation mécanique d'autre part comme période propre des petites oscillations d'un « pendule simple ». Les deux « intégrales » *u* et *u'* représentent mécaniquement les oscillations des deux pendules simples et la combinaison algébrique $x+vy=u$, permet d'exprimer le mouvement composé par les variables *x* et *y*. Lorsque Lagrange proposera, comme nous l'avons vu dans la troisième partie, de généraliser le problème de d'Alembert à un système de *n* corps, sa méthode traduira toujours la représentation mécanique selon laquelle le mouvement d'une corde lestée de *n* masses se représente comme la composition des oscillations indépendantes de *n* cordes chargées d'une seule masse ([43]) :

> Ainsi le mouvement des corps sera le même, dans ce cas que s'ils étaient pesans et qu'ils fussent suspendus chacun à un fil de longueur $1/r_m^2$ [ ([44])], la gravité étant prise pour l'unité des forces accélératrices ; d'où l'on voit que le système est susceptible d'autant de différents mouvements isochrones que l'équation $P=0$ a de racines réelles négatives et inégales. [Lagrange 1766, p.534].

Cette représentation mécanique induit une pratique algébrique consistant à ramener l'intégration d'un système de *n* équations différentielles linéaires à coefficients constants à celle de *n* équations indépendantes intégrables séparément. Lorsque, comme nous l'avons vu, Lagrange fait usage de sa méthode de réduction de l'ordre d'une équation différentielle dont des solutions particulières sont connues, les solutions particulières $\xi= Esin(t\sqrt{K}+\varepsilon)$ représentent les oscillations des pendules simples dont la composition mathématique permet d'écrire la solution « générale » du problème. La généralisation impliquée par le passage du problème de d'Alembert, limité à deux ou trois masses, au problème « général » de *n* masses, pose cependant de nouvelles questions : l'équation algébrique, « $P=0$ », à laquelle l'intégration du système différentiel est ramenée, est de degré *n*, elle ne peut donc pas être résolue algébriquement « en général », contrairement aux équations particulières obtenues pour deux ou trois masses. Or le bon usage de la méthode nécessite l'obtention de racines toutes inégales correspondant aux $n$ solutions indépendantes $\xi = E\sin(t\sqrt{K} + \varepsilon)$. Cette difficulté propre à la généralité du problème posé par Lagrange est à l'origine d'une discussion qualitative sur la nature des racines de l'équation des petites oscillations. Le problème mécanique considéré attend en effet une conclusion qualitative sur la nature des oscillations, ces dernières restent-elles bornées ou s'amplifient-elles sans limite ? La représentation mécanique sous jacente à la méthode algébrique associant les racines de l'équation $P=0$ aux périodes des petites oscillations, la stabilité du système est associée à la nature des racines : le système est stable, les oscillations bornées, à condition que les racines soient réelles, négatives et inégales ([45]).

---

[43] Pour la représentation mécanique qui sous tend les méthodes mathématique, voir plus généralement les conclusions faites par d'Alembert [1758, p.152-160], ou Lagrange [1766, p.532-539] en matière de fréquences d'oscillations.

[44] $1/r^2$ désigne une racine de l'équation.

La discussion sur la relation entre la stabilité et la nature des racines se développe chez d'Alembert sur le plan mécanique en mobilisant le principe de conservation des forces vives, « si on veut qu'ils [les corps] arrivent tous en même tems à la verticale », les oscillations sont nécessairement des expressions réelles de $t$ [d'Alembert 1758, p.152]. La nature réelle des oscillations implique alors que les « parties imaginaires » des racines de « l'équation en $v$ » se « détruisent mutuellement » [*Ibidem*, p.165-166]. Mais d'Alembert va plus loin en affirmant que les oscillations devant « toujours être extrêmement petites, par la nature du problème », l'occurrence de racines imaginaires doit être rejetée car elle imposerait la présence de termes exponentielles, croissant à l'infini [*Ibid.*, p.167]. En cas d'occurrence de *racines multiples*, la méthode des indéterminés ne permet pas de ramener le système à des équations indépendantes et d'Alembert évoque, sans le développer, un argument selon lequel il suffirait de faire subir une variation infiniment petite α aux coefficients de l'équation afin de se ramener au cas des racines inégales, la quantité α disparaissant à l'issu du calcul. Lorsqu'il entreprend en 1766 de généraliser la méthode aux oscillations d'un système quelconque de $n$ corps, Lagrange reprend à son compte l'argument de d'Alembert et propose un calcul effectif pour le cas d'une racine double $\rho_1$. Posant $\rho_2=\rho_1+\omega$, $\omega$ étant une « quantité évanouissante » ([46]), Lagrange conclut :

> […] il est évident que les termes de la valeur de $y^{(s)}$ qui répondent aux racines égales contiendront toujours l'angle $t$, et de plus des exponentielles ordinaires si ces racines sont positives, et des sinus et des cosinus si elles sont négatives. [Lagrange 1766, p.529].

Chez Lagrange comme chez d'Alembert, l'occurrence d'une racine multiple est interprétée comme engendrant une solution de la forme $P(t)e^{\delta t}$, $P(t)$ contenant des puissances de l'arc $t$ d'un degré égal à la multiplicité de la racine. La discussion aboutit à une conclusion homogène qui qualifie les solutions associées aux racines imaginaires comme aux racines multiples de mécaniquement incorrectes car impliquant des oscillations croissantes à l'infini ([47]).

Si le système mécanique est stable alors les racines de l'équation algébrique associée sont réelles, négatives et inégales. Les enjeux de cette implication changent de nature lorsque Lagrange aborde, en 1774, le problème des variations séculaires des planètes par une mathématisation similaire à celle des petites oscillations d'un fil traitées en 1766 ([48]). Dans la situation des petites oscillations d'un fil, les présupposés de stabilité du problème mécanique impliquaient des propriétés des racines de l'équation associée- les oscillations sont petites donc les racines sont réelles -, dans le cadre de l'astronomie, la stabilité du système du monde n'est au contraire pas assurée ([49]) :

> Avant de terminer cet Article, nous devons encore remarquer que, quoique nous ayons supposé que les racines *a, b, c,…* de l'équation en *x* soient réelles et inégales, il peut néanmoins arriver qu'il y en ait d'égales ou d'imaginaires ; mais il est facile de résoudre ces cas par les méthodes connues : nous observerons seulement que, dans le cas des racines égales, les valeurs de s, $s_1, s_2, \ldots, u, u_1, u_2, \ldots$ contiendront des arcs de cercle, et que dans celui des racines imaginaires ces valeurs contiendront des exponentielles ordinaires ; de sorte que, dans l'un et l'autre cas, les quantités dont il s'agit croîtront à mesure que *t* croît ; par conséquent la solution précédente cessera d'être exacte au bout d'un certain temps ; mais heureusement ces cas ne paraissent pas avoir lieu dans le Système du monde. [Lagrange 1778, p.665].

[45] Dans cette conclusion, il faut cependant distinguer entre les deux qualités des racines que sont leurs natures (réelles, négatives) et leurs multiplicités. Dans le cas où les racines sont imaginaires ou réelles positives, la méthode s'applique et les oscillations s'expriment par une exponentielle réelle : elles ne sont donc pas bornées. En cas d'occurrence de racines multiples, la méthode ne permet plus l'obtention de $n$ équations indépendantes et l'expression des solutions n'est plus valable.

[46] On retrouvera des raisonnements basés sur l'introduction d'infiniment petits pour traiter le cas des racines multiples chez Cauchy en 1839 et Sylvester en 1882. Comme l'a montré Thomas Hawkins en 1975, un raisonnement par introduction d'infiniment petits pourrait être rendu mathématiquement rigoureux (et par conséquent arriver à des conclusions opposées à celles de Lagrange) en utilisant le théorème de Bolzano Weierstrass sur l'ensemble des matrices orthogonales, fermé et borné dans $M_n$(Ë).

[47] Ce type de solutions peut effectivement intervenir dans la résolution d'un système d'équations différentielles d'ordre $n$ dont la matrice aurait des valeurs propres multiples. Il n'intervient cependant pas dans l'hypothèse symétrique ou plus généralement diagonalisable. L'absence de notion d'indépendance linéaire ne permet pas d'imaginer deux solutions s'exprimant à partir d'une même fonction $e^{\delta t}$ et tout de même « indépendantes ». Il y a ici deux familles indépendantes : celle des fonction $(e^{\delta_i t})$ dans l'espace des fonctions dérivables et celle des vecteurs de Ë$^n$ auquel est isomorphe l'espace des solutions de l'équation. C'est l'indépendance de cette dernière famille qui est nécessaire et qui peut donc être réalisée pour des racines multiples.

La stabilité du système solaire ne pouvant être prise comme hypothèse, la généralisation de la méthode des petites oscillations d'un fil aux oscillations des planètes sur leurs orbites inverse la logique de la discussion sur la nature des racines qui se focalise sur l'implication : les racines sont réelles donc le système est stable. En 1781, Lagrange aborde la question de la stabilité du système solaire par un calcul effectif des solutions de l'équation algébrique associée. Le cas de Saturne et Jupiter est traité séparément en raison de leur éloignement ; l'étude des orbites des planètes Mars, Terre, Vénus et Mercure conduit à une équation du 4e degré que Lagrange résout par application de sa méthode de réduction du degré par substitutions [Lagrange 1781, p.311]. L'application numérique, basée sur les tables astronomiques et trigonométriques, fournit quatre racines réelles distinctes. Pour Lagrange, cette conclusion numérique ne permet cependant pas de lever le « doute » sur la stabilité du système des quatre planètes tant le calcul des valeurs des racines est tributaire des valeurs attribuées aux masses des planètes. Lever le doute nécessiterait une démonstration algébrique de l'inégalité des racines mais la nature du problème - la résolution d'une équation algébrique générale- condamne à abandonner l'ambition d'une démonstration « générale » pour s'orienter vers un « artifice particulier ».

> Mais, comme les racines que nous venons de trouver dépendent des valeurs supposées aux masses des Planètes, on pourrait douter si, en changeant ces valeurs, on ne tomberait peut être pas dans les racines égales ou imaginaires. Pour lever tout à fait ce doute, il faudrait pouvoir démontrer, en général, que, quelles que soient les valeurs des masses, pourvu seulement qu'elles soient positives, les racines de l'équation dont il s'agit sont toujours réelles et inégales. Cela est facile lorsqu'on ne considère à la fois que l'action mutuelle de deux Planètes […] ; mais cette équation se complique et s'élève à mesure que le nombre des Planètes augmente ; c'est pourquoi il devient de plus en plus difficile de juger à priori de la qualité des racines. Cependant il ne parait pas impossible de parvenir, par quelque artifice particulier, à décider cette question d'une manière générale ; et comme c'est un objet également intéressant pour l'analyse et pour l'Astronomie physique, je me propose de m'en occuper. En attendant, je me contenterai de remarquer que, dans le cas présent, les racines trouvées sont trop différentes entre elles pour qu'un petit changement dans les masses adoptées puisse les rendre égales, et encore moins imaginaires. [Lagrange 1781, p.316].

---

[48] En raison de l'effectivité que requiert de la méthode son usage en astronomie, Lagrange développe une nouvelle méthode pour la détermination des constantes des solutions aux conditions initiales et consistant en un procédé d'itération permettant de mener les calculs pour les planètes du système solaire [1778, p.634]. Cette méthode, dont la description dépasse l'objet de ce travail, est parfois désignée comme « méthode de Le Verrier » en algèbre linéaire [Gantmacher 1958]. La demande d'effectivité provenant de l'application de la méthode de Lagrange à l'astronomie se fait particulièrement ressentir pour le cas des variations séculaires qui nécessitent la prise en compte d'observations sur une très longue période [Lagrange, 1778 p.634]. La question de l'extension des prévisions mathématiques des variations séculaires à un temps indéterminé restera très débattue au XIXe siècle pour finalement recevoir une réponse négative chez Poincaré. Voir [Robadey 2006] et [Laskar 1992].

[49] Tout le problème de la stabilité du système solaire se résume donc au calcul des valeurs propres de deux matrices réelles (7×7 à l'époque car Neptune n'était pas encore découverte) *A* et *B*. Si ces valeurs propres sont toutes réelles et distinctes, alors les solutions sont quasi-périodiques, et les excentricités et inclinaisons ne présentent que des variations périodiques autour de leurs valeurs moyennes, mais si l'une des valeurs propres a une partie imaginaire non nulle, on obtient une instabilité exponentielle, une excentricité peut alors devenir très grande et entraîner la possibilité de rencontre de deux planètes.

L'« artifice particulier » évoqué par Lagrange est construit par Laplace qui, dans une suite de mémoires publiés dans les années 1770-1780, se propose de donner une méthode permettant de prévoir les positions des corps célestes pour une durée « extrêmement longue » ([50]).

> Ces valeurs ne peuvent servir que pour un temps limité, après lequel, les excentricités devenant fort grandes, la supposition qu'elles sont peu considérables et d'après laquelle elles ont été trouvées cesse d'être exacte; on ne peut donc étendre à un temps quelconque les résultats obtenus dans cette supposition, qu'autant que l'on est assuré que les racines de l'équation *(k)* sont toutes réelles et inégales ; mais il serait très difficile d'y parvenir par la considération directe de cette équation. [Laplace 1776, p.464]

L'idéal déterministe de Laplace donne une nouvelle importance à la question de la stabilité du système du monde qu'il s'agit de démontrer en « dehors toute hypothèse » numérique sur la masse des planètes. La réalité et l'inégalité des racines de l'équation des inégalités séculaires des planètes est démontrée par Laplace au moyen de l'« artifice particulier » consistant à représenter la stabilité physique, la constance du moment cinétique, par une condition algébrique sur les coefficients du système différentiel. En 1789, Laplace démontre à l'aide des seuls « rapports remarquables » des coefficients disposés en miroirs dans les systèmes mécaniques que l'occurrence de solutions instables est contradictoire à cette condition ([51]).

La démonstration de Laplace semble clore la discussion engagée par d'Alembert et Lagrange sur la nature des racines de l'équation des petites oscillations et, dans le cadre de la mécanique, la discussion sur la stabilité se développera essentiellement par la prise en compte des termes non linéaires dans la mathématisation des inégalités séculaires des planètes. Comme nous allons le voir en abordant les travaux de Cauchy, avant la remise en cause des conclusions de Lagrange sur la relation entre multiplicité des racines et stabilité du système par Weierstrass en 1858 puis Yvon-Villarceau en 1870, la poursuite de la discussion passera par la géométrie analytique et impliquera une disparition des représentations mécaniques associées à la pratique algébrique élaborée par Lagrange.

## 2. La discussion dans un cadre géométrique chez Cauchy.

La discussion des petites oscillations change de nature lorsque, en 1829, la représentation mécanique sous-jacente est supplantée par une représentation issue de la géométrie analytique ([52]). Malgré l'objet annoncé par son titre, « Sur l'équation à l'aide de laquelle on détermine les inégalités séculaires des planètes », le mémoire de Cauchy ne porte pas sur un sujet mécanique mais s'insère dans un ensemble de travaux sur la classification des surfaces du second degré motivés, d'une part, par la charge d'enseignement de géométrie analytique à l'Ecole Polytechnique et, d'autre part, par le développement de la théorie de l'élasticité à partir des recherches de Fresnel sur la double réfraction de la lumière ([53]).

Le mémoire de 1829 est annoncé, trois ans avant sa publication, par une note à l'Académie des sciences dont ni l'intitulé, « Mémoire sur l'équation qui a pour racines les moments d'inertie principaux d'un corps solide et sur diverses équations du même genre » [1826], ni le contenu ne font mention des petites oscillations. La note de 1826 propose pourtant déjà une analogie entre mécanique et géométrie, elle a pour objet la résolution, par une même méthode de recherche d'« axes principaux », du problème mécanique de la rotation d'un solide et de la classification géométrique des « fonctions homogènes du second degré » ([54]). Le problème de la rotation d'un solide fait référence à la formulation analytique donnée par Lagrange en 1775 au résultat mécanique de Leonhard Euler [1758 et 1765] sur l'existence, pour un corps solide en rotation, de trois axes orthogonaux par rapport auxquels le moment d'inertie est nul ([55]). La question de la rotation d'un solide, bien que mettant en jeu un système de trois équations différentielles dont l'intégration passe par la résolution d'une équation cubique, n'avait, chez Lagrange, aucune relation avec le problème des petites oscillations d'un système de corps. D'une part, la considération d'un solide limitait l'étude au cas de trois variables et permettait une démonstration par l'absurde de la réalité des racines de l'équation cubique associée [Lagrange 1788, p.399], d'autre part les représentations sous jacentes au problème des axes principaux de rotation étaient très différentes de celles des petites oscillations. La recherche des axes principaux s'appuyait en effet sur des changements de variables interprétés dans un cadre géométrique comme des changements d'axes orthogonaux [*Ibidem,*.p.395] et cette représentation géométrique porte l'analogie entre le problème de la rotation d'un solide et celui de la recherche des axes principaux des coniques et quadriques explicitée par Cauchy en 1826. Comme nous l'avons vu dans la deuxième partie, c'est à une deuxième analogie portée cette fois par le caractère spécifique de l'équation des petites oscillations, que les trois problèmes de la rotation d'un solide, des petites oscillations et de la classification des surfaces du second degré sont abordés par une même méthode en 1829.

[50] Laplace et Lagrange font leurs calculs à l'ordre 1 et négligent les termes qui dépendraient du carré des masses des planètes. Ces approximations limitent la duré de validité de leurs théorie comme le montrera notamment Le Verrier en 1856. Voir à ce sujet [Laskar 1992, p.184-187].

[51] La preuve de Laplace, basée sur la conservation de l'intégrale première donnant l'énergie mécanique d'un système $L-V=C$, présuppose la forme erronée donnée par Lagrange aux solutions en cas de racines multiples (voir la note 37), elle est donc circulaire pour ce cas de figure. Pour une formulation de la démonstration de Laplace dans le cadre des mathématiques contemporaines, voir [Hawkins 1975, p.15].

[52] Dans son étude de la mathématisation chez Cauchy, A. Dahan Dalmedico a mis en évidence l'influence de représentations propres à la géométrie analytique sur les méthodes algébriques développées par Cauchy pour la théorie de l'élasticité (ellipsoïdes de contraintes et de déformations) et aujourd'hui perçues comme appartenant à l'algèbre linéaire [Dahan Dalmedico 1992, p.296].

[53] La théorie de l'élasticité élaborée par Cauchy nécessite en effet la description des « surfaces représentées par des équations du second degré » pour décrire les variations directionnelles de l'élasticité d'un solide. Voir à ce sujet [Dahan Dalmedico 1992, p.236].

Le véritable objet du mémoire de 1829 est une généralisation à *n* variables d'une méthode élaborée dans le cadre de la géométrie analytique et publiée en 1828 dans un mémoire intitulé « Sur les centres, les plans principaux, et les axes des surfaces du second degré ». Cette généralisation manifeste d'une part un héritage de la pratique spécifique élaborée par Lagrange pour le problème des petites oscillations et, d'autre part, des représentations géométriques propres au problème de la classification des coniques et quadriques. Le théorème énoncé par Cauchy en 1829 ne porte pas sur les systèmes différentiels étudiés par Lagrange mais sur les différentes « formes » revêtues par l'équation d'une surface homogène sujette à des changements d'axes. Parmi ces formes, il en existe une et une seule qui « renferme seulement les carrés » des coordonnées :

> Théorème II. – Etant donnée une fonction homogène et du second degré de plusieurs variables $x,y,z,...,$ on peut toujours leur substituer d'autres variables $\xi$, $\eta$, $\zeta$, ...liées à $x$, $y$, $z$, ... par des équations linéaires tellement choisies que la somme des carrés de $x$, $y$, $z$, ... soit équivalente à la somme des carrés de $\xi$, $\eta$, $\zeta$, ... , et que la fonction donnée de $x$, $y$, $z$... se transforme en une fonction de $\xi$, $\eta$, $\zeta$,... homogène et du second degré mais qui renferme seulement les carrés de $\xi$, $\eta$, $\zeta$,... […] dans le cas particulier où les variables $x$, $y$, $z$, sont au nombre de trois seulement, l'équation *(7)* se réduit à celle qui se représente dans diverses questions de Géométrie et de Mécanique, par exemple, dans la théorie des moments d'inertie ; et le théorème I fournit les règles que j'ai données dans le III[e] volume des *Exercices* comme propres à déterminer les limites des racines de cette équation. Alors aussi les équations *(22)* sont semblables à celles qui existent entre les cosinus des angles que forment trois axes rectangulaires quelconques avec les axes coordonnés, supposés eux-mêmes rectangulaires, et le théorème II correspond à une proposition de Géométrie, savoir que, par le centre d'une surface, on peut mener trois plans perpendiculaires l'un à l'autre, et dont chacun la divise en deux parties symétriques. [Cauchy 1829, p.20].

Comme nous l'avons vu dans la deuxième partie de cet article, le problème est formulé par Cauchy comme celui de la recherche des extremums d'une « fonction réelle homogène et du second degré $s = f(x,y,z,....)$ […] si l'on assujettit ces variables à l'équation de condition $x^2+y^2+z^2 + .... = 1$ » [Cauchy 1829, p.176] ([56]). Par les considérations développées dans les *Leçons sur le Calcul infinitésimal*, les conditions sur les dérivées partielles de $s$ se traduisent par la résolution d'un système linéaire dont l'équation algébrique $S=0$ associée est définie comme « une fonction alternée » (un déterminant) des quantités « comprises dans le Tableau » ([57]):

$$(10)\begin{cases} (A_{xx}-s)x+A_{xy}y+A_{xz}z+...=0, \\ A_{xy}x+(A_{yy}-s)y+A_{yz}z+...=0, \\ A_{xz}x+A_{yz}y+(A_{zz}-s)z+...=0, \\ ......................................, \end{cases}$$

> En opérant ainsi, on trouvera, par exemple, pour $n=2$, *(12)* $S=(A_{xx}-s)(A_{yy}-s)-A_{xy}^2$;pour $n=3$, *(13)* $S=(A_{xx}-s)(A_{yy}-s)(A_{zz}-s)-A_{yz}^2(A_{xx}-s)-A_{xz}^2(A_{yy}-s)-A_{xy}^2(A_{zz}-s)+2A_{xy}A_{xz}A_{yz}$ […] et, généralement, on obtiendra pour $S$ une fonction de $s$, qui sera entière et du degré $n$. [*Ibidem*, p.176].

---

[54] Les axes principaux sont les vecteurs propres des matrices associées aux formes quadratiques dans une base orthonormée.

[55] Pour la physique contemporaine, le théorème des axes principaux d'un solide en rotation revient à énoncer que le mouvement de l'oscillateur spatial se décompose en trois oscillateurs harmoniques indépendants par projection sur trois axes.

A chacune des racines $s_i$ de l'équation $S=0$ correspond un « système de valeurs » $x_i, y_i, z_i,...$ solution du système *(10)*. La caractérisation donnée par Cauchy aux relations entretenues par ces systèmes deux à deux manifeste un héritage de la pratique élaborée par Lagrange en 1766. Comme nous l'avons vu dans la troisième partie, Lagrange avait introduit deux jeux de coefficients indéterminées $(\lambda_i^{(s)})$ et $(v_i^{(s)})$ pour déterminer les solutions $y^{(s)} = \frac{v_1^{(s)}}{Q_1}\theta_1 + \frac{v_2^{(s)}}{Q_2}\theta_2 + ...+\frac{v_n^{(s)}}{Q_n}\theta_n$ du système aux conditions initiales $(\lambda'_i y' + \lambda''_i y'' + .... + \lambda^{(n)}_i y^{(n)} = \theta_i)$. Il avait ensuite déduit des rapports remarquables $A_{xy}=A_{yx}$ du système que les deux jeux $(\lambda_i)$ et $(v_i)$ étaient déterminés par le même système d'équations correspondant au système initial. Les deux jeux de coefficients $\lambda', \lambda'', \lambda'',...$ et $v', v'', v''', v''',...$correspondent donc à deux systèmes de solutions du système *(10)*, notés par Cauchy $x_1, y_1, z_{1,...}$ et $x_2, y_2, z_2,....$ Comme l'avait fait Lagrange avec ses jeux d'indéterminés, Cauchy multiplie la première ligne du système *(10)* par $x_2$, la seconde par $x_1$, etc. puis ajoute les lignes de manière à « éliminer le coefficient $A_{xx}$ » :

Relations entre les $(x_i)$ et $(x_j)$ établies par Cauchy.

$(s_2-s_1)x_1x_2 + A_{xy}(x_2y_2-x_1y_2)+A_{xz}(x_2z_1-x_1z_2)+...=0.$

$A_{xy}(y_2x_1-y_1x_2)+(s_2-s_1)y_1y_2 + A_{yz}(y_2z_1-y_1z_2)+...=0.$

*[...]* etc.

Enfin, si l'on ajoute membre à membre les équations (17), (18), (19), etc., on trouvera

(20) $(x_1x_2+y_1y_2+z_1z_2+...)(s_2-s_1)=0.$

Donc, toutes les fois que les racines $s_1,s_2$ seront inégales entre elles, on aura

(21) $x_1x_2+y_1y_2+z_1z_2+...=0$ ;

Relations entres les $(\lambda_i)$ et $(\lambda_j)$ de Lagrange.

$$(i)\left\{\begin{array}{c} v'_m\lambda'_1+v''_m\lambda''_1+...+v_m^{(n)}\lambda_1^{(n)}=0 \\ v'_m\lambda'_2+v''_m\lambda''_2+...+v_m^{(n)}\lambda_2^{(n)}=0 \\ .... \\ v'_m\lambda'_{m-1}+v''_m\lambda''_{m-1}+...+v_m^{(n)}\lambda_{m-1}^{(n)}=0 \\ v'_m\lambda'_m+v''_m\lambda''_m+...+v_m^{(n)}\lambda_m^{(n)}=Q_m \\ v'_m\lambda'_{m+1}+v''_m\lambda''_{m+1}+...+v_m^{(n)}\lambda_{m+1}^{(n)}=0 \\ .... \\ v'_m\lambda'_n+v''_m\lambda''_n+...+v_m^{(n)}\lambda_n^{(n)}=0 \end{array}\right\}$$

relations qui se résument par l'équation :

$(\rho^2-\rho^2_m)\ [v'\lambda'+v''\lambda''+v'''\lambda'''+.... + v^{(n)}\lambda^{(n)}]=0$

De la relation *(21)*, Cauchy déduit les relations entretenues par deux séries de solutions associées à deux racines distinctes :

$$(22)\begin{cases} x_1^2+y_1^2+z_1^2+...=1, & x_1x_2+y_1y_2+z_1z_2+...=0 & x_1x_n+y_1y_n+z_1z_n+...=0, \\ x_2x_1+y_2y_1+z_2z_1+...=0, & x_2^2+y_2^2+z_2^2+...=1, & x_2x_n+y_2y_n+z_2z_n+...=0, \\ ................................., & ...................................., & ............................, \\ x_nx_1+y_ny_1+z_nz_1+...=0, & x_nx_2+y_ny_2+z_nz_1+...=0, & x_n^2+y_n^2+z_n^2+...=1, \end{cases}$$

Les racines, $s_1, s_2, ...,s_n$ de l'équation $S=0$ permettent d'exprimer l'équation de la surface dans le repère formé de ses axes principaux, $s_1x_1^2+s_2x_2^2+...+s_nx_n^2=0$. Les systèmes $(x_i, y_i, z_i, ...)$ correspondants aux racines $s_i$ de l'équation $P=0$ représentent alors les coordonnées d'axes orthonormés. La pratique algébrique mise en œuvre par Cauchy est indissociable d'une

[56] C'est-à-dire la recherche des extremums de la fonction sur sur ce que nous désignerions comme la sphère unité de Ë$^n$.

[57] A propos du mémoire de 1815 de Cauchy sur le calcul des fonctions alternées envisagé dans le cadre de la genèse du concept de groupe, voir [Dahan Dalmedico 1979]

représentation d'orthogonalité issue de la géométrie analytique et, lorsque ce dernier généralise sa méthode de recherche d'axes principaux à $n$ variables, il interprète les procédés algébriques de la méthode des coefficients indéterminées de Lagrange comme portant sur des changements d'axes orthogonaux ([58]).

Nous avons détaillé dans la troisième partie de cet article le caractère spécifique du rôle donné par Lagrange à l'équation algébrique $S=0$ dans l'expression des oscillations $\lambda_i e^{\rho t}$ du système mécanique des petites oscillations, le coefficient $\rho$ étant une racine de l'équation $S=0$ et $\lambda$ un polynôme en $\rho$ obtenu par une factorisation d'une équation identique, à un facteur multiplicatif près, à l'équation $S=0$ :

$$(\rho^2-\rho^2_m)\,[v'\lambda'+v''\lambda''+v'''\lambda'''+\ldots + v^{(n)}\lambda^{(n)}]=0,$$

L'expression $Q_m=v'\lambda'+v''\lambda''+v'''\lambda'''+\ldots + v^{(n)}\lambda^{(n)}$ s'identifiait à la factorisation de $S=0$ par l'expression $\rho^2-\rho^2_m$, c'est-à-dire à la valeur de $S/(s-s_m)$ pour $s=s_m$ et permettait la formulation des oscillations ($y^{(s)}$) du système mécanique à partir de la seule donnée de l'équation $S=0$. Cauchy reformule les expressions de Lagrange par des relations d'orthogonalité euclidienne et des extractions de sous déterminant $P_{xy}$ du système (10) :

> Soit maintenant $R$ ce que devient la fonction $S$, lorsque, dans le Tableau (11), on supprime tous les termes appartenant à la première colonne horizontale, ainsi qu'à la première colonne verticale ; et $Q$ ce que devient la même fonction, quand on supprime, en outre, les termes renfermés dans les deuxièmes colonnes horizontale et verticale. Enfin, désignons par $P_{uv}$ ce que devient $S$, lorsque l'on supprime dans le Tableau (11) les termes qui appartiennent à la même colonne horizontale que le binôme $A_{uu}-s$ avec ceux qui appartiennent à la même colonne verticale que $A_{vv}-s$ [...] on aura évidemment (23) $R = P_{xx}$. [*Ibid.*, p.178].

« En faisant abstraction de la première » équation du système (10), c'est-à-dire en fixant sa valeur comme arbitraire, le calcul des déterminants permet de déterminer $x, y, z$ par des formules analogue aux expressions $\frac{vi}{Qi}$ de Lagrange (l'expression $X^2+Y^2+Z^2+\ldots$ correspond au terme noté par Lagrange $Q_i^2$ ) ([59]) :

> et l'on conclura des équations (10), en faisant abstraction de la première,
>
> $$(24)\ \frac{x}{P_{xx}} =- \frac{y}{P_{xy}} =- \frac{z}{P_{xz}} =-\ldots$$
>
> Posons d'ailleurs, pour abréger, (25) $X=P_{xx} = R,\ Y=-P_{xy},\ Z=-P_{xz}, \ldots$
>
> La formule (24), combinée avec la formule (3) , donnera
>
> $$(26)\ \frac{x}{X} = \frac{y}{Y} = \frac{z}{Z} = \ldots = \frac{1}{\sqrt{X^2+Y^2+Z^2+\ldots}}$$
>
> et, si l'on désigne par
>
> $$(27)\ X_1, Y_1, Z_1, \ldots, X_2, Y_2, Z_2, \ldots, X_n, Y_n, Z_n, \ldots,$$
>
> les systèmes de valeurs de $X, Y, Z, \ldots$ correspondants aux racines $s_1, s_2, \ldots, s_n$ de l'équation (7),
>
> $$(28)\begin{cases} \frac{x_1}{X_1}=\frac{y_1}{Y_1}=\frac{z_1}{Z_1}=\ldots=\frac{1}{\sqrt{X_1^2+Y_1^2+Z_1^2+\ldots}}, \\ \frac{x_2}{X_2}=\frac{y_2}{Y_2}=\frac{z_2}{Z_2}=\ldots=\frac{1}{\sqrt{X_2^2+Y_2^2+Z_2^2+\ldots}}, \\ \ldots\ldots\ldots\ldots\ldots\ldots\ldots\ldots, \\ \frac{x_n}{X_n}=\frac{y_n}{Y_n}=\frac{z_n}{Z_n}=\ldots=\frac{1}{\sqrt{X_n^2+Y_n^2+Z_{n1}^2+\ldots}}, \end{cases}$$
>
> [*Ibid.*].

Ces dernières formules généralisent les expressions de changements d'axes orthonormés permettant l'écriture d'une conique dans son système d'axes principaux, obtenues en 1828 par un usage du calcul des « fonctions alternées » analogue à celui que nous avons vu à l'œuvre dans la reformulation par Cauchy des procédés de Lagrange. La question de la nature des racines de l'équation $S=0$ se pose alors de façon nouvelle. L'écriture d'une fonction homogène de $n$ variables comme somme de carrés, nécessite l'utilisation de formules de changements d'axes *(28)* dans lesquelles le terme de droite peut donner « une forme 0/0 » en cas d'occurrence de racine multiple dans l'équation caractéristique $S=0$ ce qui se formule, dans le calcul des fonctions alternés, comme l'occurrence d'une racine commune $s_i$ entre l'équation $S=0$ et l'une des équations $R=0$, $Q=0$, associées aux sous déterminant $X_i$ extraits du système $S$. Pour une racine double par exemple, le facteur $(s-s_1)^2$ dans l'expression de $S$ ne serait pas compensé par les facteurs $(s-s_1)$ du mineur $X_i$ et l'expression $x_1 = \dfrac{X_i}{\sqrt{X_1^2+Y_1^2+Z_1^2+..}} = \dfrac{X_i}{\frac{|S|}{s-s_i}}$ pourrait prendre « une forme 0/0 ». Cauchy, qui critique les méthodes qui tendent « à faire attribuer aux formules algébriques une étendue indéfinie, tandis que […] la plupart de ces formules subsistent uniquement sous certaines conditions », précise la limite de validité de sa méthode :

> D'après ce qui a été dit ci-dessus, il ne peut rester de doutes sur l'exactitude du théorème I, si ce n'est dans le cas où quelques valeurs de $s$ vérifieraient à la fois les deux équations (36) $S=0$, $R=0$, $Q=0$, …prises consécutivement. [*Ibid*, p.235].

[58] En des termes qui nous sont contemporains, les relations d'orthogonalité duale entre les $v$ et les $\lambda$ données par Lagrange donnent lieu à des relations d'orthogonalité euclidienne chez Cauchy.

[59] Si $B(s)$ est la matrice des mineurs de la matrice caractéristique $A-sI$, alors $B(s)(A-sI)=\Delta(s)I$ et si $D$ est un vecteur colonne, $(A-sI)X=D \Leftrightarrow X=B(s)D$. Les coordonnées des vecteurs propres sont obtenues par les colonnes non nulles de la matrice adjointe $Ad(A-sI)$ de la matrice caractéristique, c'est-à-dire la matrice des cofacteurs de $A-sI$. Pour une formulation contemporaine de ce procédé voir l'expression (*) donnée en I.3. et la note n°23.

Les « formules algébriques » de la méthode « générale » ne « subsistent » pas à la « condition » de multiplicité des racines. Le problème posé par l'occurrence de racines multiples se manifeste alors comme un cas singulier limitant l'« étendue » d'une méthode algébrique et nécessitant pas conséquent un raisonnement spécifique que Cauchy mène par l'introduction d'infiniment petits ([60]). Les événements politiques de 1830 sont proches et ce n'est qu'à son retour d'exil que Cauchy propose, dans les années 1839-1840, une nouvelle méthode, indépendante de la nature des racines, et basée sur le calcul des résidus [Cauchy 1839b]. Dès son élaboration par Cauchy en 1826, le calcul des résidus se présentait comme une méthode homogène de résolution des équations différentielles ordinaires à coefficients constants. Dans le « Mémoire sur l'analogie des puissances et des différences et sur l'intégration des équations linéaires » [1825], Cauchy avait publié une solution des équations différentielles linéaires d'ordre $n$ à coefficients constants, de la forme $\frac{d^n y}{dx^n} + a_1 \frac{d^{n-1} y}{dx^{n-1}} + \ldots + a_{n-1}\frac{dy}{dx} + a_n y = 0$, par l'emploi d'une « analogie » entre factorisation symbolique de l'équation différentielle et factorisation du polynôme $F(r)= r^n+a_1r^{n-1}+\ldots+a_{n-1}r+a_n$. Ce calcul symbolique perd cependant toute validité en cas d'occurrence de racines multiples dans l'équation $F(r)=0$ ([61]). Au contraire, le calcul des résidus permet de donner une résolution homogène indépendante de la nature des racines de l'équation $F(r)=0$ sous la forme $y = \Re es \frac{\Phi(r)e^{rx}}{((F(r))}$, $\Phi(r)$ étant un polynôme arbitraire ([62]). A l'occasion de deux mémoires, publiés en 1839-1840 dans un contexte de travaux sur la propagation des ondes, l'élasticité et la théorie de la lumière, Cauchy revient sur sa méthode de 1829 et donne le nom d'équation caractéristique à l'équation algébrique au cœur des analogies formelles entre divers problèmes de physique mathématique et de géométrie analytique. A cette époque, les préoccupations de Cauchy ont changé et le recours à la géométrie qui caractérisait l'approche de 1829 fait place à de nouvelles

---

[60] Nous avons déjà vu un tel raisonnement à l'œuvre chez d'Alembert et Lagrange dans la partie III de cet article.

[61] Cauchy écrit l'équation différentielle sous la forme $(D-r_1)(D-r_2)\ldots(D-r_n)y = -f(x)/a_0$, où les $r_i$ sont les racines réelles du polynôme $F$. Alors, $(D-r_1)y_{n-1}=f(x)/a_0$, $(D-r_2)y_{n-2} = y_{n-1}$, …, $(D-r_n)y= y_1$, et comme $(D-r)y=f(x)$ implique $y = e^{rx} \int e^{-rx}f(x)dx$, Cauchy obtient la formule suivante qui nécessite l'occurrence de racines $r_i$ toutes distinctes :

$$y = \frac{e^{r_n x}}{a_0}\int e^{(r_n - r_{n-1})x}\left(\int e^{(r_{n-1} - r_{n-2})x}\left(\ldots\int e^{(r_2 - r_1)x} f(x)dx\ldots\right)dx\right)dx$$

Pour tout complément, consulter [Dahan Dalmedico 1992, p.197].

[62] C'est-à-dire par l'intégrale curviligne $\frac{1}{2\pi i}\int \frac{\phi(r)e^{rx}}{F(r)}dr$ prise autour des racines de $F$.

méthodes développées pour le traitement de problèmes de physique mathématique comme ceux qui étaient à l'origine de la discussion des petites oscillations et auxquels le calcul de résidus donne une résolution homogène ([63]). Le qualificatif homogène s'oppose ici aux cas singuliers qui « encombrent » les méthodes de l'algèbre en limitant l'« étendue » de la « subsistance » des formules algébriques. Le renouvellement de la signification de la généralité qui accompagne cette exigence d'homogénéité, va porter la discussion des petites oscillations de Cauchy à Weierstrass, en passant par Carl Jacobi, Sylvester, Johann Dirichlet,Carl Borchardt et Hermite ([64]).

### 3. La discussion de Weierstrass sur la transformation des couples de fonctions homogènes.

Comme nous l'avons déjà mentionné, la propriété énoncée par Lagrange en 1766 sur la relation entre la stabilité des systèmes mécaniques et la nature des racines de l'équation des petites oscillations ne sera pas contestée avant la parution, à près d'un siècle de distance, d'un mémoire de Weierstrass sur les couples de fonctions homogènes du second degré ([65] ):

> Nous avons donc démontré le théorème suivant.
> Soit $\Phi$, $\Psi$ des fonctions entières homogènes du deuxième degré de $n$ variables $x_1, x_2, ...x_n$, à coefficients réels et dont la première conserve un signe constant et ne s'annule pas identiquement pour les valeurs réelles de $x_1, x_2, ...x_n$. Le déterminant de la fonction $s\Phi$-$\Psi$ est une fonction entière de degré $n$ de la grandeur variable $s$ qui ne s'annule que pour des valeurs réelles. Soit $s_1, s_2, ..., s_m$ ces valeurs, alors le déterminant ne contient que les facteurs suivant, $(s-s_1)^{\lambda 1}$ $(s-s_2)^{\lambda 2}$ $...(s-s_m)^{\lambda m}$, où $\lambda_1, \lambda_2, ...\lambda_m$ sont des nombres positifs dont la somme est $n$ ; il existe alors des fonctions homogènes parfaitement déterminées, $\theta_1, \theta_2,...., \theta_m$, de degré deux et de variables $x_1, x_2, ...,x_n$, telles que $\Phi$ et $\Psi$ s'expriment sous la forme :
>
> $$\Phi = \theta_1 + \theta_2 + ... + \theta_m$$
> $$\Psi = s_1\theta_1 + s_2\theta_2 + ... + s_m\theta_m$$
>
> où $\theta_\mu$, ou $-\theta_\mu$ selon que $\Phi$ reste positive ou négative, représente une somme quadratique de $\lambda_\mu$ fonctions linéaires réelles des grandeurs $x_1, x_2, ...., x_n$. [Weierstrass 1858, p.242, *traduction F.B.*]

La nature des racines caractéristiques n'intervient pas dans cet énoncé qui donne donc une résolution à la fois *homogène* et *algébrique* des différents problèmes abordés par Lagrange, Laplace et Cauchy, désormais insérés dans une même théorie des couples de fonctions homogènes. Quelle que soit la multiplicité des racines de l'équation caractéristique $f(s)=0$, la stabilité mécanique du système des petites oscillations et la possibilité de transformer les équations des surfaces du second degré en sommes de carrés par des changement d'axes orthogonaux sont assurées par la présence d'une « circonstance remarquable » (« bemerkenswerther ») propre aux fonctions $\Phi$ définies positives (ou négatives) et selon laquelle une racine $s_i$ de multiplicité $\mu$-$1$ de toutes les équations obtenues par les mineurs $f(s)_{\alpha\beta}=0$, sera nécessairement une racine de multiplicité $\mu$ de l'équation $f(s)=0$.

> [Le problème] a été complètement résolu par Cauchy, Jacobi etc. pour le cas où l'on ne trouve aucune grandeur égale parmi $s_1, s_2, ..., s_n$. Il n'est pas encore résolu en revanche dans les circonstances exceptionnelles où les racines de l'équation $f(s)= 0$ ne sont pas différentes l'une de l'autre, la difficulté qui se présente alors aurait déjà du être éclaircie et je propose de l'examiner attentivement plus en détail. Je ne pensais pas initialement qu'une solution serait possible sans des discussions spécifiques aux nombreux cas différents qui peuvent se produire. Il me fallait espérer que la résolution du problème soit susceptible d'une méthode indifférente à la multiplicité des grandeurs $s_1, s_2, ..., s_n$. [*Ibidem*, p. 233, *traduction F.B.*].

---

[63]Cette évolution des méthodes de Cauchy a été mise en évidence par A. Dahan Dalmedico pour le cas de la théorie de l'élasticité : à une première approche des années 1820-1830, basée sur une représentation géométrique des ellipsoïdes de pressions, succède une seconde méthode, basée sur un point de vue moléculaire, et développée par Cauchy à partir de 1839. « Dans le cadre de la deuxième théorie moléculaire de Cauchy, l'appareil conceptuel qui avait favorisé la figuration géométrique des pressions et des déformations locales, en particulier par le biais des quadriques de pressions et de déformations est quelque peu perdu de vue » [Dahan Dalmedico 1992, p.285], « le calcul des résidus et le calcul symbolique sur les opérateurs différentiels prend une place prédominante » [*Ibidem*, p. 211].

[64] Nous renvoyons pour ces différents développements de la discussion à l'étude détaillée de [Hawkins 1977] qui a notamment émis l'hypothèse qu'un mémoire de Borchardt [1846], pourrait être à l'origine de l'intérêt de Weierstrass pour la question. Hawkins fait référence à une lettre de Weierstrass, citée par Dugac [1972, p.158] mettant en évidence la relation privilégiée entre les deux hommes dans les années 1850.

[65] Hawkins [1977] a suggéré que l'attention de Weierstrass sur la question de la stabilité des systèmes mécaniques pourrait avoir été suscitée par un article de Dirichlet [1846], publié en appendice à la troisième édition de la mécanique analytique [1853], qui reprend la démonstration de Lagrange selon lequel un état d'équilibre d'un système est stable si la fonction de potentiel $V$ y atteint un minimum strict.

Le succès de Weierstrass à élaborer « une méthode indifférente à la multiplicité des grandeurs $s_1, s_2, ..., s_n$ » tient à l'expression qu'il donne au problème comme relevant d'une question de « transformation » d'une paire de fonctions $\Phi$, $\Psi$ et qui sous tend une pratique qui combine les variables $x_i$, $y_i$ en « formes » linéaires. Partant de l'écriture du couple de formes $s\Phi$-$\Psi$, Weierstrass exprime les variables $x_\alpha$ par les formes linéaires $\Phi_\alpha$ de manière à représenter l'écriture de $\Phi$ et $\Psi$ comme combinaison linéaire des formes quadratiques $\Phi_\alpha\Phi_\beta$ ([66]) :

On pose, si $\alpha, \beta, \gamma$, sont des nombres appartenant à $1,2,..., n$

$$\frac{1}{2}\ \frac{d\Phi}{dx\alpha} = \Phi\alpha,\ \frac{1}{2}\ \frac{\delta\Psi}{\delta x\alpha} = \Psi_\alpha$$

[...]

$$(5)\ \Phi = \Sigma_\alpha \Phi_\alpha x_\alpha$$
$$\Psi = \Sigma_\alpha \Psi_\alpha x_\alpha$$

[…] et si on sort $x_1, x_2, ..x_n$ de

$$s\Phi_1 - \Psi_1,\ s\Phi_2 - \Psi_2, ...,\ s\Phi_n - \Psi_n$$

On obtient :

$$(1)\ x_\alpha = \Sigma_\beta \{ \frac{f(s)_{\alpha\beta}}{f(s)} (s\Phi_\beta - \Psi_\beta) \}$$

où $f(s)_{\alpha\beta}$ est une fonction entière de degré inférieur au $(n-1)^s$ […]

$$(2)\ f(s)_{\alpha\beta} = f(s)_{\beta\alpha}$$

[*Ibid.*, *traduction F.B.*].

La variable $x_\alpha$ ne contenant pas de puissance de $s$, un développement « en puissances décroissantes » de l'expression $\Sigma_\beta \{ \frac{f(s)_{\alpha\beta}}{f(s)} (s\Phi_\beta - \Psi_\beta) \}$ présente uniquement un terme constant. D'un autre côté, le développement de $\frac{f(s)_{\alpha\beta}}{f(s)}$ ne contient que des termes de degré inférieurs à $s^{-1}$ car $deg(f(s)) = n$ (le déterminant de $\Phi$ n'est pas nul). Ne peuvent donc subsister dans le développement que des termes provenant du développement de $s\Phi_\beta$ (la notation $[...]_{s-1}$ représente le coefficient de $s^{-1}$ dans le développement) :

$$(3)\ x_\alpha = \Sigma_\beta [\frac{f(s)_{\alpha\beta}}{f(s)}]_{s-1}\Phi_\beta$$

et comme les coefficients de $s^{-1}$ disparaissent :

$$(4)\ \Sigma[f\frac{f(s)_{\alpha\beta}}{f(s)}\ _{s-1}\Psi_\beta = \Sigma[s\frac{f(s)_{\alpha\beta}}{f(s)}]_{s-1}\Phi_\beta$$

Ces deux identités, combinées à (4) et (2) permettent d'écrire :

$$(6)\ \Phi = \Sigma [\frac{f(s)_{\alpha\beta}}{f(s)})]_{s-1}\Phi_\alpha\Phi_\beta,$$

$$(7)\ \Psi = \Sigma [\frac{f(s)_{\alpha\beta}}{f(s)}]_{s-1}\Psi_\alpha\Phi_\beta = \Phi = \Sigma [\frac{f(s)_{\alpha\beta}}{f(s)}]_{s-1}\Phi_\alpha\Psi_\beta,$$

[…]

$$(8)\ \Psi = \Sigma[s\frac{f(s)_{\alpha\beta}}{f(s)}]_{s-1}\Phi_\alpha\Phi_\beta$$

La formule *(8)* exprime les fonctions $\Phi$ et $\Psi$ comme des sommes des mêmes formes quadratiques $\Phi_\alpha\Phi_\beta$. Pour obtenir la transformation désirée des fonctions $\Phi$ et $\Psi$, il ne reste qu'à réorganiser l'écriture par regroupement des variables en groupes de formes $\theta_i$ et $\Theta_i$ pour lesquelles $\Theta_i = s\theta_i$.

---

[66] Les $\Phi_\alpha$ sont des formes linéaires qui constituent une base du dual $E^*$. Il s'agit donc d'exprimer les variables $x_\alpha$ dans cette base duale de manière à avoir une écriture des formes quadratiques comme somme de produits de formes linéaires $\Phi = \Sigma k_{s1}\Phi_\alpha\Phi_\beta$.

La relation $x_\alpha=\Sigma_\beta\{ \frac{f(s)_{\alpha\beta}}{f(s)} (s\Phi_\beta- \Psi_\beta)$ exprime les $x_\alpha$ à l'aide de factorisations $f(s)_{\alpha\beta}$ du déterminant $f(s)$ ([67]). L'expression $\frac{f(s)_{\alpha\beta}}{f(s)}$, analogue au $\frac{Pxx}{s}$ de Cauchy, manifeste un héritage des pratiques spécifiques élaborées par Lagrange et Cauchy, depuis le mémoire publié par ce dernier en 1829, elle s'interprète comme le quotient d'un mineur par le déterminant $f(s)$ dont il est extrait et la circonstance remarquable énoncée par Weierstrass implique que cette expression est toujours définie, même en l'occurrence de racines multiples de $f(s)=0$ ([68]). Pour réorganiser l'écriture par groupes de variables $\theta_i$ et $\Theta_i$ pour lesquelles $\Theta_i=s\theta_i$, Weierstrass décompose les fractions $\frac{f(s)_{\alpha\beta}}{f(s)}$ en éléments simples à l'aide des pôles donnés par les racines $s_1, ..., s_m$ de l'équation $f(s) =0$ ([69]).

$$[\frac{f(s)_{\alpha\beta}}{f(s)}]_{s-1} = \Sigma_\mu [\frac{f(s)_{\alpha\beta}}{f(s)}]_{(s-s\mu)-1}$$

Les variables sont structurées en groupes caractérisés ([70]), pour chaque racine $s_\mu$, par les termes comportant son résidu $[\frac{f(s)_{\alpha\beta}}{f(s)}]_{(s-s\mu)-1}$ ([71]):

$$(10)\ \Sigma_{\alpha\beta} \{ \frac{f(s)_{\alpha\beta}}{f(s)} \}_{(s-s\mu)-1}\, \Phi_\alpha\Phi_\beta = \theta_\mu$$

$$\Sigma_{\alpha\beta} \{ \frac{sf(s)_{\alpha\beta}}{f(s)} \}_{(s-s\mu)-1}\, \Phi_\alpha\Phi_\beta = \Theta_\mu$$

de sorte que

$$\Phi = \theta_1 + .... + \theta_m,$$
$$\Psi = \Theta_1 + .... + \Theta_m$$

La circonstance remarquable caractéristique des formes définies positives permet d'établir que :

$$s\frac{f(s)_{\alpha\beta}}{f(s)} = s_\mu \frac{f(s)_{\alpha\beta}}{f(s)} + (s-s_\mu)\frac{f(s)_{\alpha\beta}}{f(s)} \text{ et } \frac{s-s_\mu}{f(s)} \text{ ne devient pas } \infty \text{ si } s=s_\mu$$

Par conséquent, $\{ s\frac{f(s)_{\alpha\beta}}{f(s)} \}_{(s-s\mu)-1} = \frac{f(s_\mu)_{\alpha\beta}}{f'(s_\mu)}$, ce qui permet de déduire :

---

[67] D'un point de vue qui nous est contemporain, si on note $B(s)$ la matrice des cofacteurs de la matrice caractéristique $sA-B$, alors

$$(sA-B)B(s) =f(s) \Rightarrow \forall\, \alpha\, f(s)x_\alpha = \Sigma_\beta\, b_{\alpha\beta}(s)\, [sA-B)_{\beta\alpha}$$

[68] La circonstance remarquable est démontrée de la manière suivante. Si on développe le quotient $f(s)_{\alpha\gamma}/f(s)$ en puissances croissantes de $(s-s_\mu) : f(s)_{\alpha\gamma}/f(s) = h_{\alpha\gamma}(s-s_\mu)^l +...$, alors la théorie du déterminant implique :

$$\frac{d(\frac{f(s)_{\gamma\gamma}}{f(s)})}{ds} = -\sum_{\alpha\beta} A_{\alpha\beta} \frac{f(s)_{\alpha\gamma}}{f(s)} \frac{f(s)_{\beta\gamma}}{f(s)} \text{ et } lh_{\gamma\gamma}(s-s_\mu)^{l-1} + ... = -\sum_{\alpha\beta} A_{\alpha\beta} h_{\alpha\gamma} h_{\beta\gamma} (s-s_\mu)^{2l} + ...$$

$\Phi$ est définie et ne s'annule donc pas pour des valeurs réelles des variables $x_i$, non toutes nulles, le coefficient de $(s-s_\mu)^{2l}$ ne peut pas être nul et $2l \geq l-1$ donc $l \geq -1$. Chaque terme $f(s)_{\alpha\gamma}$ doit par conséquent être divisible par $(s-s_\mu)^{\lambda-1}$.

[69] Cet examen des pôles des quotients par Weierstrass pourrait, selon Thomas Hawkins, manifester une influence de la résolution donnée par Cauchy en 1839 dans le cadre du calcul des résidus [Hawkins 1975, p.131].

[70] L'exemple suivant permet d'illustrer la notation de Weierstrass :

$$\frac{x+1}{x^2+x-2} = \frac{1}{3} [ \frac{2}{x-1} + \frac{1}{x+2}] \text{ signifie } [\frac{x+1}{x^2+x-2}]_{s-1} = 2/3 + 1/3$$

[71] Cette méthode est analogue à la manière dont ont exprime aujourd'hui les exponentielles de matrices intervenant dans la résolution des systèmes différentiels. S'il s'agit par exemple de calculer l'exponentielle de

$$M = \begin{pmatrix} 1 & 4 & -2 \\ 0 & 6 & -3 \\ -1 & 4 & 0 \end{pmatrix}$$

Le polynôme caractéristique est: $P_M= -(x-2)^2(x-3)$. On décompose $1/F$ en éléments simples, $\frac{1}{F} = \frac{1}{X-3} - \frac{X-1}{(X-2)^2}$

Le théorème de Bezout permet de déterminer les projecteurs sur les sous espaces caractéristiques, $(X-2)^2 -(X-1)(X-3)=1$

$$U_1 = (-X-1) \text{ et } U_2=1,\ P_1 = U_1Q_1 = -(X-1)(X-3) \text{ et } P_2 = (X-2)^2$$

Les projecteurs sont $p_1 = P_1(M)$ et $p_2 = P_2(M)$, $exp(M) = e^\lambda (I_3 + 1/1! (M-2I_3) p_1 {}^{+ e\lambda 2} (M-I_3)p_2$

$$\exp(M) = \begin{pmatrix} -6e^2+3e^3 & -4e^2+4e^3 & 10e^2-6e^3 \\ -6e^2+3e^3 & -3e^2+4e^3 & 9e^2-6e^3 \\ -7e^2+3e^3 & -4e^2+4e^3 & 11e^2-6e^3 \end{pmatrix}$$

$$(13).\ \Phi = \theta_1 + .... + \theta_m,$$
$$\Psi = s_1\,\theta_1 + .... + s_m\,\theta_m,$$

Les regroupements des variables par la structure du développement en éléments simples de la fraction $\frac{f(s)_{\alpha\beta}}{f(s)}$ donnent donc la transformation désirée pour le couple de fonctions $(\Phi, \Psi)$, pour tout $\mu$ de $\{1,2,...m\}$, $\Theta_\mu = s_\mu\theta_\mu$.

## Conclusion.

Après avoir spécifié un moment de référence nous avons constitué un corpus regroupant en un ensemble cohérent des textes publiés sur la période 1766-1874. Nous avons d'abord associé l'identité de ce corpus au caractère spécifique d'une équation algébrique, l'équation des petites oscillations. En étudiant l'origine du corpus dans les travaux de Lagrange et ses évolutions chez Laplace, Cauchy et Weierstrass nous avons ensuite mis en évidence la permanence d'une pratique consistant à exprimer les solutions d'un système linéaire par des factorisations polynomiales de l'équation caractéristique de ce système, c'est-à-dire par les expressions polynomiales (*) déjà mentionnées qui mettent en jeu les expressions $\frac{\frac{P_{1i}}{S}(x)}{x - s_j}$ ainsi que les quotients successifs $\frac{\Delta_{i+1}}{\Delta_i}$ des mineurs principaux de $S$. Nous avons vu comment la pratique spécifique de jeu sur les primes et les indices des systèmes linéaires mise en œuvre par Lagrange en 1766 était à l'origine de l'une des spécificités de l'équation des petites oscillations, l'occurrence de « rapports remarquables » en miroirs dans les coefficients des systèmes linéaires. En suivant la discussion sur les racines de l'équation des petites oscillations, nous avons vu un héritage de la pratique élaborée par Lagrange se manifester chez Laplace, Cauchy et Weierstrass et, en même temps, l'identité de l'équation des petites oscillations évoluer. Caractérisée par une *pratique indissociable d'une équation,* l'identité de notre corpus présente donc une nature *algébrique* qu'il nous faut à présent préciser.

A l'exception de Jordan et, comme nous allons le voir, Darboux en 1874, aucun des auteurs du corpus ne se réfère à une théorie algébrique. L'identité algébrique de la discussion ne peut donc se caractériser, avant 1874, comme une identité théorique, elle se place davantage du côté des pratiques que des méthodes. A la permanence d'une pratique polynomiale de résolution des systèmes linéaires il faut en effet opposer la diversité des méthodes élaborées dans les cadres théoriques distincts au travers desquels se déploie la discussion. La méthode élaborée par Lagrange pour ramener l'intégration d'un système à celle de $n$ équations indépendantes était, comme nous l'avons vu, indissociable d'une représentation mécanique, cette méthode s'avère très différente de celle que nous avons vu à l'œuvre chez Cauchy en 1829 et qui était sous tendue par des représentations géométriques de changements d'axes orthogonaux. Dans le cadre du corpus, une même pratique algébrique est donc employée aux seins de méthodes et cadres théoriques différents dans lesquels elle adopte des représentations diverses. Le point de vue selon lequel les procédés algébriques ne portent pas de représentations propres mais prennent leurs « significations » relativement aux théories dans lesquels ils sont employés est celui que Kronecker oppose en 1874 à l'organisation algébrique que veut donner Jordan à la théorie des formes bilinéaires. Pour Kronecker, dans la théorie des formes bilinéaires comme dans d'autres cadres théoriques, le « travail algébrique » n'a qu'une « place relative », il « est effectué au service d'autres disciplines mathématiques dont il reçoit ses fins et dont dépendent ses objectifs », il doit par conséquent « être seulement considéré comme le moyen et non comme le but de la recherche » [Kronecker, 1874, 367, traduction F.B.]. Chez Kronecker comme chez la plupart des auteurs du corpus, le travail algébrique n'est pas une fin en soi mais sert des objectifs issus de contextes divers et de disciplines comme la mécanique, la géométrie ou l'arithmétique. Si, parallèlement à cette variabilité des contextes, une identité algébrique se manifeste dans la discussion sur l'équation des petites oscillations, cette identité ne tient pas à la perception d'une théorie algébrique sous jacente mais à une perspective historique. C'est en effet par des références à des travaux plus anciens que les auteurs du corpus comme Laplace, Cauchy ou Weierstrass donnent à leurs travaux une identité dépassant les contextes dans lesquels ceux-ci s'inscrivent. Cauchy par exemple, donne à son mémoire de 1829 une identité dépassant le contexte géométrique de la recherche des axes principaux des coniques et quadriques par sa référence à Lagrange qui identifie la pratique algébrique dont Cauchy s'inspire et à l'« équation à l'aide de laquelle on détermine les inégalités séculaires des planètes » qui se rapporte implicitement aux travaux de Laplace et à la question de la nature des racines caractéristiques. Avant qu'une identité théorique comme celle de la théorie des formes bilinéaires ne permette d'identifier le corpus à l'histoire d'une théorie (comme Kronecker qui se réfère à l'« histoire de la théorie des faisceaux de formes quadratiques »), c'est une identité historique qui permet aux auteurs d'identifier leurs travaux à une discussion de nature algébrique.

De l'origine de la discussion chez Lagrange à sa fin chez Weierstrass et Jordan, la revendication de généralité est la principale motivation qui amène des auteurs à insérer leurs propres travaux dans une perspective historique plus large. Pour Lagrange, la capacité de sa méthode à résoudre un problème mécanique pour un système quelconque de corps lui confère un caractère général. La nécessité d'obtenir des racines caractéristiques distinctes ne réduit pas la généralité polynomiale car ce cas de figure est considéré comme incompatible avec la situation des petites oscillations. L'intervention de Cauchy dans la discussion est d'abord motivée par la généralisation à $n$ variables que permet le problème des petites oscillations d'une méthode élaborée pour deux ou trois variables dans un cadre géométrique. Les « formules algébriques » généralisant les expressions de changements d'axes orthonormés ne « subsistent » cependant pas à la « condition » de multiplicité des racines caractéristiques conduisant à la forme $\frac{0}{0}$. Pour Cauchy, le problème posé par l'occurrence de racines multiples manifeste un cas singulier limitant l'« étendue » d'une formule et nécessitant un traitement spécifique initialement mené par introduction d'infiniment petits. C'est par opposition à de tels cas singuliers qui « encombrent » les formules de l'algèbre que Cauchy introduit dès 1826 le calcul des résidus à l'aide duquel il donnera en 1839 une résolution homogène au problème des petites oscillations. Le renouvellement de la signification de la généralité qui accompagne cette exigence d'homogénéité va porter la discussion des petites oscillations de Cauchy à Weierstrass, en passant par Jacobi, Borchardt, Hermite, Dirichlet et Sylvester. La question se focalise alors sur le problème de la multiplicité des racines, d'abord abordé par des méthodes différentielles par Sturm, Cauchy et Sylvester. Ce dernier, après une première série de travaux consacrée à la nature des racines algébriques (1840-1842) articule la méthode différentielle à la théorie du déterminant et ses représentations géométriques en termes de « rectangle symétrique » par rapport à sa « diagonale ». Ses travaux qui, avec ceux d'Arthur Cayley, sont à l'origine de la théorie des invariants et voient l'introduction des termes « matrices » et « mineurs » (1850-1852) ([72]), sont investis dans un cadre arithmétique en relation avec les préoccupations d'Hermite sur les formes quadratiques et la décomposition en quatre carrés (1853-1856).

La question de la généralité est indissociable de l'identité algébrique du corpus mais tout en donnant, comme nous l'avons vu, une conclusion algébrique, générale et homogène à la discussion, Weierstrass s'appuie en 1858 sur des procédés de « transformation » de « formes » linéaires qui manifestent l'orientation arithmétique qu'ont pris les travaux sur l'équation des petites oscillations dans les années 1850. En réponse à un mémoire de Borchardt intitulé « Neue Eigenschaft der Gleichung mit deren Hülfe man die saecularen Storungen der Planeten bestimmt » [1846], Hermite adresse à Borchardt une lettre que celui-ci publiera dans son journal en 1857 sous le titre « Sur l'invariabilité du nombre des carrés positifs et des carrés négatifs dans la transformation des polynômes homogènes du second degré ». Hermite interprète les déterminants employés par Cauchy en 1829 comme des invariants caractérisant les classes d'équivalences des couples de formes quadratiques :

> Dans le cas où vous le jugeriez convenable, vous pourriez publier la démonstration suivante, du principe découvert par Jacobi, et employé par lui à la démonstration des belles formules pour les conditions de réalité des racines des équations algébriques, que vous avez données dans votre Mémoire sur l'équation à l'aide de laquelle, etc. Rien d'ailleurs n'est plus simple que d'établir ce principe que j'énoncerai ainsi : *Quelque substitution réelle que l'on emploie pour réduire un polynôme homogène du second degré à une somme de carrés, le nombre des coefficients de ces carrés qui auront un signe donné sera toujours le même.* [Hermite 1857, p.429].

[72] Au sujet des origines de la théorie des invariants, consulter [Parshall 1989 et 2006].

Dans un cadre arithmétique, le problème est présenté sous une forme semblable à la loi d'inertie, « réduire un polynôme homogène du second degré à une somme de carrés », et abordé par les méthodes propres à la théorie des formes quadratiques. On recherche les invariants qui caractérisent les coefficients des formes canonique*s* permettant de caractériser les classes d'équivalences des formes. Comme le formulera clairement Darboux qui, en 1874 généralisera les procédés élaborés par Hermite en 1857 pour élaborer une nouvelle démonstration du théorème de Weierstrass de 1858, l'invariant du couple de formes ($f(x_1,...,x_n)$, $x_1^2+...+x_n^2$) qu'est le déterminant $|f(x_1, x_2,...,x_n)-\lambda(x_1^2+x_2^2+...+x_n^2)|$ ne permet pas de caractériser les différentes classes d'équivalences et il faut donc aller au-delà du déterminant en considérant ses « mineurs » (annexe 3) ([73]). Contrairement à la loi d'inertie de la théorie arithmétique des formes quadratiques, le problème porte ici sur la transformation simultanée d'un couple de formes envisagé comme un *polynôme* de formes quadratiques $f(x_1, x_2,...,x_n)-\lambda(x_1^2+x_2^2+...+x_n^2)$. En raison de son caractère polynomial, le problème s'insère pour Darboux dans « la théorie algébrique des formes quadratiques ». La mathématique des formes est-elle une algèbre ou une arithmétique ? C'est sur cette question que se développe la querelle entre Jordan et Kronecker de 1874. Pour Kronecker qui revendique l'héritage de Gauss, la théorie des formes bilinéaires est de nature arithmétique et ses méthodes doivent par conséquent relever d'un calcul d'invariants obtenus par des procédés effectifs de calculs de p.g.c.d. ([74]). Jordan donne au contraire une organisation algébrique à la théorie qu'il structure par l'action des groupes de substitutions (substitutions linéaires, orthogonales etc.) et une méthode de « réduction » des « formes » à des « formes canoniques simples ».

Comme nous l'avons vu dans la première partie, le corpus de la discussion avait joué le rôle d'une histoire commune au travers de laquelle de premières identités s'étaient manifestées entre les deux théorèmes de Jordan et Weierstrass. L'étude que nous avons menée dans cet article nous permet d'aller plus loin. Elle montre que la discussion sur l'équation des petites oscillations s'identifie plus profondément à une pratique algébrique commune à Jordan et Weierstrass et consistant à transformer la forme des systèmes linéaires par des décompositions de la forme polynomiale d'une équation algébrique. Contrairement à la théorie arithmétique des formes quadratiques dans laquelle le terme « forme » bénéficie d'une définition mathématique explicite en relation avec la notion de relation d'équivalence, les significations des termes « formes » et « transformations » restent le plus souvent implicites et évoluent considérablement dans le cadre algébrique de la discussion. Dans sa Mécanique Analytique de 1788, Lagrange distinguait ce qu'on « était en droit d'attendre de la Dynamique », donner les équations fondamentales du mouvement, et la « charge » incombant « au calcul intégral » de donner à ces équations une « forme simple » propre à leurs intégrations.

[73] Darboux voit notamment la notion de diviseur élémentaire de Weierstrass en germe dans les invariants introduits par Sylvester en 1851 pour caractériser les intersections de coniques et quadriques. C'est à l'occasion de ces travaux que Sylvester introduit les termes de « mineurs » d'une « matrice ». Voir à ce sujet [Brechenmacher 2006d].

[74] C'est en effet dans ses travaux fondateurs sur l'arithmétique des formes quadratiques que les notions de formes et de classes d'équivalences avaient été définies par Gauss pour les fonctions homogènes à coefficients entiers ou rationnels. Lors de la querelle de 1874, Kronecker introduit ce que l'on désignerait aujourd'hui comme les facteurs invariants d'une matrice $A$ sur un anneau principal. Si $D(\lambda)=|A-\lambda I|$, dans la suite des p.g.c.d. des mineurs successifs $\Delta_r(\lambda), \Delta_{r-1}(\lambda), ..., \Delta_1(\lambda)$, chaque polynôme est divisible par le précédent et si les quotients correspondants sont désignés par $i_1(\lambda), i_2(\lambda),...i_r(\lambda)$ et appelés les polynômes invariants de la matrice $A(\lambda)$ la décomposition de ces polynômes en facteur irréductibles distincts sur $K$ donne les diviseurs élémentaires de la matrice $A(\lambda)$ sur le corps $K$.

$$i_1(\lambda)=[\varphi_1(\lambda)]^{s_1}...[\varphi_s(\lambda)]^{c_s},\ i_2(\lambda)=[\varphi_1(\lambda)]^{d_1}...[\varphi_s(\lambda)]^{d_s},\ \dots,\ i_r(\lambda)=[\varphi_1(\lambda)]^{l_1}...[\varphi_s(\lambda)]^{l_s}\ (c_k>d_k>....>l_k)$$

> Cependant comme ces équations peuvent avoir différentes formes plus ou moins simples, & surtout plus ou moins propres pour l'intégration, il n'est pas indifférent sous quelle forme elles se présentent d'abord ; & c'est peut-être un des principaux avantages de notre méthode, de fournir toujours les équations de chaque problème sous la forme la plus simple relativement aux variables qu'on y emploie, & de mettre en état de juger d'avance quelles sont les variables dont l'emploi peut en faciliter le plus l'intégration. [Lagrange 1788].

Le choix entre les différentes formes sous lesquelles se présentent les équations était normé par un critère de simplicité lié au calcul intégral. Pour le problème des petites oscillations l'existence d'une forme intégrable du système était garantie par la représentation mécanique du problème selon laquelle le mouvement d'une corde chargée de *n* masses se décompose mécaniquement en oscillations de *n* cordes chargées d'une seule masse. Si la méthode de Lagrange consistait bien à ramener l'intégration d'un système d'équations différentielles à celle de *n* équations indépendantes, le changement de la *forme* du système ne s'appuyait pas sur un procédé de transformation mais sur la résolution d'une équation algébrique donnant les paramètres mécaniques du système. Chez Cauchy la notion de « transformation d'une fonction homogène » était indissociable d'une représentation géométrique de transformation des axes de coordonnées d'une conique. Si la question qui semble aujourd'hui essentielle, et que l'on formulerait comme la caractérisation de la forme d'une matrice non nécessairement diagonalisable car présentant des valeurs propres multiples, n'avait pas été abordée par Lagrange, Laplace ou Cauchy il ne s'agissait pas là d'un défaut de généralité comme le sanctionne Kronecker en 1874 mais d'une question étrangère à la manière dont étaient envisagés les termes « formes » et « transformations » d'un système. Cette question était au contraire naturelle dans le cadre arithmétique des travaux d'Hermite et Weierstrass des années 1850 qui ouvraient la voie aux premiers travaux d'Elwyn Christoffel sur les formes bilinéaires d'une part ([75]), dans le cadre algébrique des travaux de Jordan sur les groupes de substitutions d'autre part. Dans la résolution générale que donne ce dernier au problème de l'intégration des systèmes d'équations linéaires à coefficients constants, la décomposition en facteurs irréductibles du polynôme caractéristique permet de regrouper les variables en « un certain nombre de séries » et de « ramener » le système d'équations à une « suite » de « formes simples » dont l'intégration est connue (annexe 2) [Jordan 1871, p.787]. Aux diverses significations associées au terme « formes » chez Lagrange, Laplace ou Cauchy, répond chez Hermite, Weierstrass, Kronecker, Jordan ou Darboux une mathématique dont la « forme », la « transformation » est l'objet. Cette mathématique qui permet de poser, en toute généralité, la question de la caractérisation de la forme d'un système linéaire est-elle une algèbre ou une arithmétique ? Comme nous le détaillons dans un article consacré à la querelle de 1874 [Brechenmacher 200 ?], si Jordan et Kronecker se réfèrent à une histoire commune dans laquelle ils identifient une pratique spécifique, consistant à aborder les transformations des couples de formes bilinéaires $(A,B)$ par la décomposition de la forme polynomiale de l'équation $|A+sB|=0$, les deux géomètres opposent la généralité de deux méthodes (calculs d'invariants, réductions canoniques), deux cadres théoriques (formes bilinéaires, substitutions), deux disciplines (arithmétique, algèbre) dans lesquels s'insèrent cette pratique partagée par les deux savants. Comme l'illustre cette controverse, la manière dont une mathématique des formes s'élabore au XIX[e] siècle, en manifestant une tension entre arithmétique et algèbre et la rencontre de différents héritages comme ceux portés par la théorie des formes quadratiques et la pratique algébrique spécifique à la discussion sur l'équation à l'aide de laquelle on détermine les inégalités séculaires des planètes, permet d'enrichir l'histoire de l'algèbre linéaire souvent abordée par l'intermédiaire de l'histoire d'une théorie, d'une notion, d'un mode de raisonnement et la problématique de l'émergence de structures.

[75] Dans deux mémoires successifs, consacrés à ce que nous désignerions comme le cas hermitien dans lequel les racines sont réelles comme dans le cas symétrique, et publiés en 1864, Christoffel s'appuie sur les résultats de Weierstrass de 1858 pour généraliser des travaux mécaniques de Clebsch [1860] et des résultats d'[Hermite 1854, 1855, 1856] sur les résidus biquadratiques, les propriétés arithmétiques des nombres $a+bi$ ($a$ et $b$ entiers) et la décomposition d'un entier en somme de quatre carrés. Le premier mémoire se présente comme un développement mathématique du second, consacré à des problèmes de petites oscillations dans le cadre de la théorie de Cauchy selon laquelle la lumière correspond aux petites vibrations de molécules ponctuelles d'éther soumises à des forces attractives et répulsives. Le procédé nécessite l'intégration d'un système d'équations différentielles linéaires dont les coefficients sont des constantes complexes [Mawhin 1981] présentant des « propriétés remarquables car ils apparaissent comme des extensions de ces équations qui interviennent dans la théorie des perturbations séculaires des planètes et dans tant d'autres recherches » [Clebsch 1860, p.326] et généralisant le cas quadratique à ce que Christoffel décrit comme des « fonctions bilinéaires » en référence à un mémoire de 1857 dans lequel Jacobi généralisait la loi d'inertie des formes quadratiques aux « fonctions bilinéaires » $\Sigma a_{ij}x_iy_j$.

## Annexe 1. Représentation graphique simplifiée du corpus formant la discussion sur l'équation à l'aide de laquelle on détermine les inégalités séculaires des planètes.

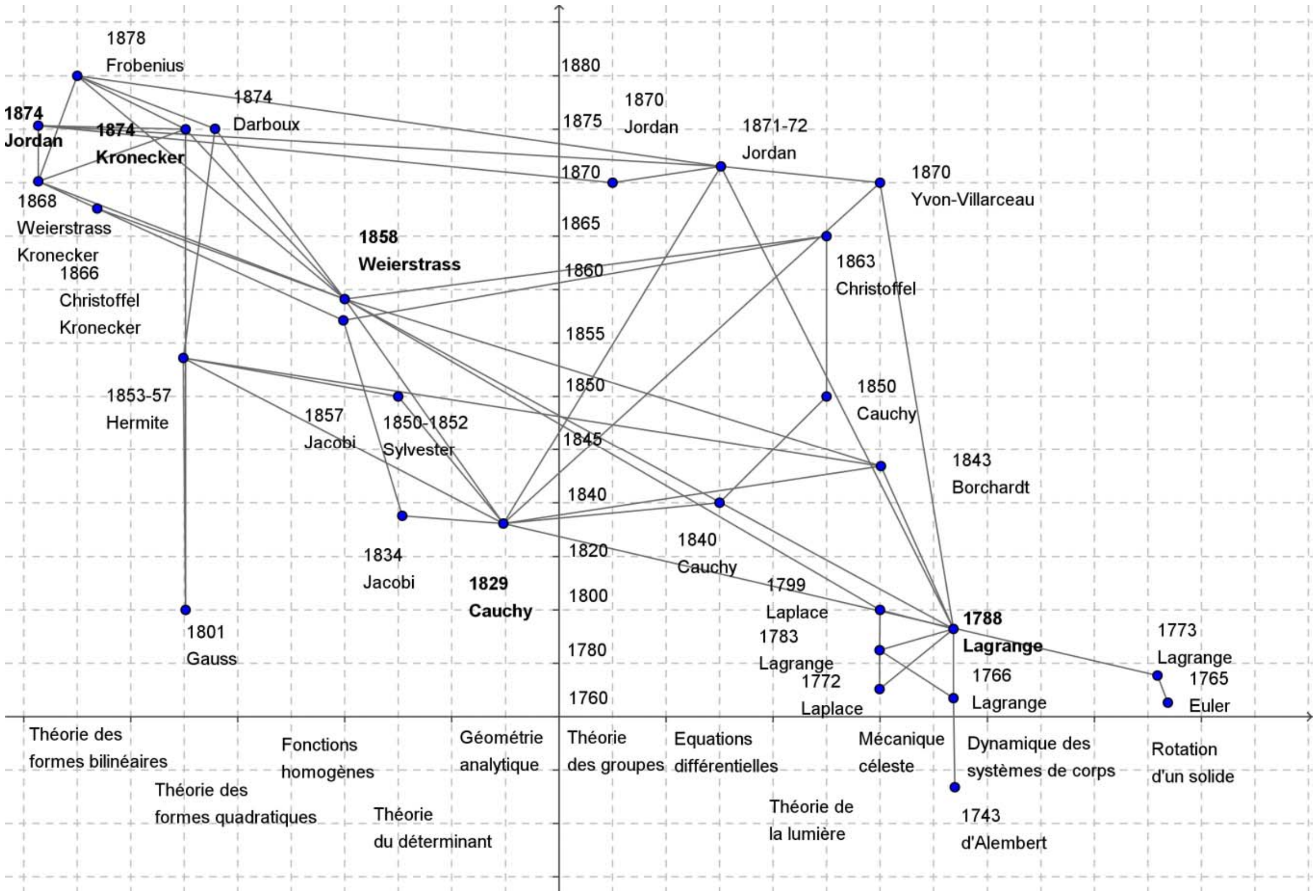

## Annexe 2. La réponse de Jordan à Yvon-Villarceau (1871).

Dans une des séances de l'hiver dernier, M. Yvon Villarceau a signalé une lacune dans le procédé généralement indiqué pour la solution d'un système d'équations différentielles linéaires à coefficients constants […]. On sait en effet, que l'intégration de ce système dépend de l'équation caractéristique […] mais le cas où cette équation a des racines égales présente une légère difficulté. On connaît en gros le moyen de la résoudre ; mais on n'a pas donné, que nous sachions, une analyse complète et embrassant tous les cas de la question […].

$$(I) \qquad \begin{aligned} \frac{dx_1}{dt} &= a_1x_1 + \ldots + l_1x_n, \ldots, \\ \frac{dx_2}{dt} &= a_2x_1 + \ldots + l_2x_n, \ldots, \\ &\ldots\ldots \\ \frac{dx_n}{dt} &= a_nx_1 + \ldots + l_nx_n, \end{aligned}$$

Ce problème peut cependant se résoudre très simplement par un procédé identique à celui dont nous nous sommes servi, dans notre Traité des substitutions, pour ramener une substitution linéaire quelconque à sa forme canonique. Nous allons ramener de même le système (I) à une forme canonique qui puisse s'intégrer immédiatement […]. Les équations différentielles prendront la forme

$$(5)\begin{cases} \frac{dy_1}{dt}=\sigma y_1,..,\frac{dy_v}{dt}=\sigma y_v \\ \frac{dx_{v+1}}{dt}=a'_1 x_{v+1}+...+k'_1 x_n + fonct.(y_1,...,y_v), \\ .... \\ \frac{dx_n}{dt}=a'_1 x_{v+1}+...+k'_{n-v} x_n + fonct.(y_1,...,y_v), \end{cases}$$

et l'on aura $\Delta=(\sigma-s)^v \Delta'$, $\Delta'$ désignant le déterminant

$$\begin{vmatrix} a'_1-s & ... & k'_1 \\ ... & ... & ... \\ a'_{n-v} & ... & k'_{n-v}-s \end{vmatrix}$$

Poursuivant ainsi, on voit que les variables indépendantes peuvent être choisies de telle sorte qu'aux $\mu$ racines égales à $\sigma$ que possède l'équation $\Delta=0$ correspondent $\mu$ variables nouvelles formant un certain nombre de séries contenant respectivement $r, r',...$ variables, $r+r'+...$ étant égal à $\mu$, et les variables d'une même série étant liées par une suite de relations de la forme

$$(6)\quad \frac{dy_1}{dt}=\sigma y_1,\ \frac{dz_1}{dt}=\sigma z_1+y,\ \frac{du_1}{dt}=\sigma u_1,+z_{1,\,...}\ \frac{dw_1}{dt}=\sigma w_1+v_1$$

[…] Soit $r$ le nombre de variables de la série $y_1,...,w_1$ ; le système des équations (6) aura évidemment pour intégrales le système suivant :

$$w_1=e^{\sigma t}\psi(t),\ v_1=e^{\sigma t}\psi'(t),...,\ y_1=e^{\sigma t}\psi^{r-1}(t),$$

$\psi(t)$ étant une fonction entière arbitraire du degré $r-1$. [Jordan 1871, p.787].

## Annexe 3. Extrait du « Mémoire sur la théorie algébrique des formes quadratiques » de Darboux (1874).

Notation : $\Delta, \Delta_1, \Delta_2, ..., \Delta_{n-1}, 1$ sont les mineurs du déterminant caractéristique. Les fonctions $\Phi_p$ sont des invariants définis par Darboux, à partir des mineurs de la forme quadratique.

Considérons une forme quadratique homogène à $n$ variables

$$f(x_1,x_2,...,x_n),$$

que nous écrirons aussi de la manière suivante :

$$\Sigma_i \Sigma_j a_{ij}x_ix_j,$$

[…] Nous supposerons, suivant l'usage $a_{ji}=a_{ij}$; par suite un terme rectangle se trouvera deux fois dans la somme précédente avec les mêmes coefficients, et l'expression développée de la forme sera

$$a_{11}x_1^2+a_{22}x_2^2+....+a_{nn}x_n^2+2a_{12}x_1x_2+...$$

[…] Cela posé, réduisons, par les procédés connus, $f$ à une somme de carrés de la forme suivante :

$$f=\varepsilon_1(x_1+a_2x_2+...+a_nx_n)^2+\varepsilon_2(x_2+\beta_3x_3+\beta_nx_n)^2+...+\varepsilon_{n-1}(x_{n-1}+\lambda_nx_n)^2+\varepsilon_nx_n^2$$

on trouvera, comme l'a prouvé M. Hermite,

$$\varepsilon_1=\Delta_{n-1},\ ...,\ \varepsilon_{n-i}=\frac{\Delta_i}{\Delta_{i-1}},\ \varepsilon_n=\frac{\Delta}{\Delta_1}$$

et par suite on aura

$$f=\Delta_{n-1}X_1^2+\frac{\Delta_{n-2}}{\Delta_{n-1}}X_2^2+...+\frac{\Delta}{\Delta_1}X_n^2$$

$X_1, X_2,...,X_n$ étant des fonctions linéaires des $n$ variables $x_1, x_2, ..., x_n$. Il résulte de cette décomposition de $f$ que le nombre des carrés positifs de la forme est égal à celui des permanences de la suite

$$\Delta, \Delta_1, \Delta_2, ..., \Delta_{n-1}, 1\ ;$$

quant aux fonctions linéaires $X_i$ elles ont été mises sous la forme de déterminants par divers géomètres et notamment par M. Weierstrass. […] Le nombre de carrés positifs de la forme ne peut changer que si $\lambda$ passe par une racine de l'équation obtenue en égalant l'invariant à zéro, et dans ce cas le nombre de carrés positifs de la forme ne peut varier d'une quantité supérieure à l'ordre de multiplicité de la racine considérée. C'est le théorème qu'il s'agissait d'établir ; mais il donne lieu à une remarque essentielle : c'est que, dans son énoncé, on pourrait entendre par ordre de multiplicité d'une racine, non plus le nombre des dérivées, mais le nombre des fonctions $\Phi,\ldots,\Phi_n$ qu'annule cette racine, quels que soient les arbitraires figurant dans ces fonctions. Ainsi une racine multiple pourra être considéré comme simple si elle n'annule pas tous les mineurs du premier ordre ; comme double si, annulant tous les mineurs du premier ordre, elle n'annule pas tous ceux du second et ainsi de suite. Faisons quelques applications de ce théorème. Soit d'abord la forme quadratique

$$f(x_1, x_2,\ldots,x_n)-\lambda(x_1^2+x_2^2+\ldots+x_n^2).$$

L'équation qu'on obtient en égalant l'invariant à zéro a été traitée d'abord par Cauchy, puis par une foule de géomètres : MM. Borchardt, Sylvester, etc. On déduit la réalité de ses racines du théorème que nous venons d'établir [….]. De plus, si elle a des racines multiples, une racine d'ordre *p* devra annuler tous les mineurs d'ordre *p-1* de l'invariant. La méthode de M. Sylvester seule permet de démontrer ce dernier point et d'écarter la difficulté relative aux racines multiples, qui se présente dans la méthode de Cauchy et de M. Borchardt. [Darboux 1874, p.367].

## Bibliographie.